\documentclass[10pt]{amsart}
\usepackage{amsmath,amsfonts,amsthm}
\usepackage [dvips]{graphicx}
\usepackage{amsmath}
\usepackage{graphicx,psfrag}  
\usepackage[psamsfonts]{amssymb}
\usepackage{amscd}
\usepackage[all,cmtip]{xy}
\title[The Weil-Petersson geodesic flow]{The Weil-Petersson geodesic flow is ergodic %\\(Preliminary version)
}

\author{K. Burns, H. Masur and  A. Wilkinson}
\thanks{Each  author partially 
supported by the NSF}
\address{Keith Burns \\ Department of Mathematics \\Northwestern University\\ 2033 Sheridan Rd.\\ Evanston, IL 60208.  USA.  burns@math.northwestern.edu}
\address{Howard Masur \\ Department of Mathematics \\University of Chicago\\ 5734 S. University Ave.\\
Chicago, Illinois 60637\\ USA.  masur@math.uchicago.edu}
\address{Amie Wilkinson \\ Department of Mathematics \\Northwestern University\\ 2033 Sheridan Rd.\\ Evanston, IL 60208.  USA.  wilkinso@math.northwestern.edu}
\date{\today}
\theoremstyle{plain}

\newtheorem{mainthm}{Theorem}
\newtheorem{theorem}{Theorem}[section]
\newtheorem{proposition}[theorem]{Proposition}
\newtheorem{lemma}[theorem]{Lemma}
\newtheorem{corollary}[theorem]{Corollary}

\def\Teich{\hbox{Teich} }
\def\mod{\hbox{mod} }

\def\inj{\operatorname{inj} }

\def\title{\em}

\def\bar{\overline}
\def\bd{\operatorname{\bf d}}
\def\bexp{\operatorname{\bf exp}}

\def\Sat{\operatorname{Sat}}

\def\loc{\text{loc}}

\def\uell{\underline\ell}
\def\oell{\overline\ell}
\def\Re{\mbox{Re}}
\def\Im{\hbox{Im}}

\def\cX{\mathcal X} 
\def\cY{\mathcal Y}
\def\cD{\mathcal{D}}
\def\cW{\mathcal{W}}

\def\cB{\mathcal{B}}
\def\cD{\mathcal{D}}

\def\cJ{\mathcal{J}}
\def\cM{\mathcal{M}}
\def\cR{\mathcal{R}}

\def\U{\mathcal{U}}

\def\C{\mathcal{C}}
\def\cC{\mathcal{C}}

\def\M{\mathcal{M}}
\def\T{\mathcal{T}}
\def\cT{\mathcal{T}}

\def\ok{\overline k}

\def\grad{\operatorname{grad}}

\def\Vol{\hbox{Vol}}
\def\cH{\mathcal{H}}

\def\cU{\mathcal{U}}

\def\Proj{\hbox{Proj}}
\def\transverse{\,\raise2pt\hbox to1em{\hfil$\top$\hfil}\hskip -1em \hbox
to1em{\hfil$\cap$\hfil}\,} 

\newcommand\PGL{{\operatorname{PGL}}}
\newcommand\SL{{\operatorname{SL}}}
\newcommand\PSL{{\operatorname{PSL}}}
\newcommand\MCG{{\operatorname{MCG}}}

\newcommand\R{\mathbb R}
\newcommand\RR{{\mathbb R}}
\newcommand\CC{{\mathbb C}}
\newcommand\HH{{\mathbb H}}
\newcommand\LL{{\mathbb L}}

\newcommand\ZZ{{\mathbb Z}}
 
\newlength{\figboxwidth} 

\def\la{\langle}
\def\ra{\rangle}

\begin{document}

\begin{abstract}{We prove that the geodesic flow for the Weil-Petersson metric on the moduli space of Riemann surfaces is ergodic (and in fact Bernoulli) and has finite, positive metric entropy.}
\end{abstract}

\maketitle

%\tableofcontents

\section*{Introduction}

This paper is about the dynamical properties of the Weil-Petersson geodesic flow for the moduli space of Riemann surfaces.  Our main result is that this flow is ergodic: any invariant set must have volume zero or full volume.  Ergodicity implies that a randomly chosen, unit speed Weil-Petersson geodesic in moduli space becomes equidistributed over time. What is more, the tangent vectors to such a geodesic also become equidistributed in the space of all unit tangent vectors to moduli space.

To state our result more precisely and to put it in context, we first review the basic setup from Teichm\"uller theory. Let $S$ be a surface of genus $g\geq 0$ with $n\geq 0$ punctures, and let $\cM(S)$ be the moduli space of conformal structures on $S$, up to conformal equivalence.  Assume  that $3g+n\geq 4$, which implies that  in each conformal class there is  complete hyperbolic metric.  Then $\cM(S)$ has the alternate description of the moduli space of hyperbolic structures on $S$, up to isometry.  The orbifold universal cover of $\cM(S)$ is the Teichm\"uller space $\Teich(S)$ of marked conformal structures on $S$. 

It is a classical result
due to Fricke and Klein that  $\Teich(S)$ is homeomorphic to a ball of dimension $6g-6+2n$.
Teichm\"uller space carries a natural complex structure  via a special embedding of $\Teich(S)$ into a complex representation variety $QF(S)$, called {\em quasifuchsian space}. 
Under this map, called the {\em Bers embedding}, the image of $\Teich(S)$ sits as a complex subvariety (indeed there is a biholomorphic equivalence $QF(S) \cong \Teich(S)\times \Teich(S)$).  
The orbifold fundamental group of $\cM(S)$ is the mapping class group $\MCG(S)$ of orientation preserving homeomorphisms of $S$ modulo isotopy. The mapping class group acts holomorphically on $\Teich(S)$. The stabilizer of each point is finite, which  gives $\cM(S)$ the structure of a complex orbifold.

A naturally defined and well-studied metric on Teichm\"uller space, and the focus of this paper, is the {\em Weil-Petersson metric}  $g_{WP}$, which is the K\"ahler metric induced by the Weil-Petersson symplectic form $\omega_{WP}$ and the almost complex structure $J$ on $\Teich(S)$:
%\footnote{}  
$$g_{WP}(v,w) = \omega_{WP}(v,Jw).$$
%One then defines the Weil-Petersson (WP for short) metric $g_{WP}$ on $T(\cT_{g,n})$ by $g_{WP}(v,w) = \omega_{WP}(v,iw)$.  
We refer to the Weil-Petersson metric as the WP metric, for short.
The WP metric is invariant under $\MCG(S)$ and so descends to a metric on $\cM(S)$.  It has   finite volume determined by the volume form $|\omega_{WP}^{\wedge 3g -  3 + n}|$.

A striking feature of the WP metric is its intimate connections with hyperbolic geometry, among them:
\begin{itemize}
\item the hyperbolic length of a closed geodesic (for a fixed free homotopy class on $S$) is a convex function along WP geodesics in $\Teich(S)$  \cite{WolpertBehavior};
 \item in Fenchel-Nielsen coordinates $(\ell_i,\tau_i)_{i=1}^{3g-3+n}$ on $\Teich(S)$, the WP symplectic form $\omega_{WP}$
has the simple expression \cite{WolFN}
$$\omega_{WP} =\frac12 \sum_{i=1}^{3g-3+n} d\ell_i \wedge d\tau_i.$$ 
\item the growth of the hyperbolic lengths of simple closed curves on $S$ is related to the WP volume of $\cM(S)$ \cite{Mi}; and
\item the WP metric has a formulation in terms of dynamical invariants of the geodesic flow on hyperbolic surfaces \cite{Bridgeman, CTM}.
\end{itemize}

The Weil-Petersson metric has several notable features that make it an interesting geometric object of study in its own right.
The WP metric is negatively curved, but incomplete. The sectional curvatures are neither bounded away from $0$ (except in the simplest cases of $(g,n) = (1,1)$ and $(0,4)$), nor bounded away from $-\infty$.   The  WP geodesic flow thus presents a naturally-occurring example of a singular hyperbolic dynamical system, for which one might hope to reproduce the known properties of
the geodesic flow for a compact, negatively curved manifold, such as: ergodicity, equidistribution of closed orbits, exponentially fast mixing and decay of correlations, and central limit theorem.

We summarize the previous literature on the WP geodesic flow.  Wolpert \cite{WolpertGeometry} showed that the geodesic flow is defined for all time on a full volume subset of the the unit tangent bundle $T^1 \Teich(S)$ and thus descends to
a volume-preserving flow on the finite volume quotient $\cM^1(S) := T^1\Teich(S)/\MCG(S)$.
Pollicott, Weiss and Wolpert \cite{PWW} proved   in the case $(g,n) = (1,1)$ that the geodesic flow is transitive  on $\cM^1(S)$  and periodic orbits are dense in $\cM^1(S)$ \cite{PWW}. Brock, Masur and Minsky \cite{BMM} proved transitivity and denseness of
periodic orbits for arbitrary $(g,n)$ and also showed
that the topological entropy of the geodesic flow is infinite (that is, unbounded on compact invariant sets).
Hamenst\"adt \cite{Ham} proved a measure-theoretic version of density of closed orbits:  the set of invariant Borel probability measures for the WP geodesic flow that are supported on a closed orbit is dense in the 
space of all ergodic invariant probability measures.

In this paper, we prove:
\begin{mainthm}\label{t=main} Let $S$ be a Riemann surface of genus $g\geq 0$, with $n\geq 0$ punctures.  Assume that $3g+n\geq 4$.
The Weil-Petersson geodesic flow on $\cM^1(S)$ is ergodic (and in fact Bernoulli) with respect to WP volume and has finite, positive measure-theoretic entropy.
\end{mainthm}

The Bernoulli property means that the time-$1$ map of the geodesic flow is abstractly isomorphic (as a measure-preserving system) to a Bernoulli process on a finite alphabet.  In particular it is mixing of all orders.   An interesting open question is to determine the rate of mixing of this flow.

Our basic approach to proving Theorem~\ref{t=main} is as follows.  The WP geodesic flow $\varphi_t$ preserves a finite probability volume $m$ on $\cM^1(S)$, and one can show using properties of the WP metric that $\log\|D\varphi_1\|$ is integrable with respect to the measure $m$. The Multiplicative Ergodic Theorem of Oseledec (cf. \cite[Theorem S.2.9]{KatokHasselblatt})  then implies that there is a full volume subset $\Omega\subset \cM^1(S)$ such that for every $v\in \Omega$ and every nonzero tangent vector $\xi\in  T_v\cM^1(S)$, the limit
$$
\lambda(\xi) := \lim_{t\to\infty} \frac{1}{t} \log\|D_v\varphi_t(\xi)\| 
$$
exists and is finite.  The real  number $\lambda(\xi)$ is called the {\em (forward) Lyapunov exponent} of $\varphi_t$ at $\xi$.  Observe that if $\xi$ is in the line bundle $\RR\dot\varphi(v)$ 
tangent to the orbits of the flow, then $\lambda(\xi)=0$.  We say that $\varphi_t$ is {\em nonuniformly hyperbolic} if for almost every $v\in \Omega$ and every $\xi\in  T_v\cM^1(S) \setminus \RR\dot\varphi(v)$, the Lyapunov exponent  $\lambda(\xi)$ is nonzero.

Using the fact that the WP sectional curvatures are negative, we establish that the WP geodesic flow is nonuniformly hyperbolic.    Nonuniform hyperbolicity is the starting point for a rich ergodic theory of volume-preserving diffeomorphisms and flows, developed first by Pesin for closed manifolds, and expanded by Sinai, Katok-Strelcyn, Chernov and others to systems with singularities, such as the WP geodesic flow.  The basic argument for establishing ergodicity of such systems originates with Eberhard Hopf and his proof of ergodicity for geodesic flows for closed, negatively curved surfaces \cite{Ho39}. His method
was to study the Birkhoff averages of continuous functions along leaves
of the stable and unstable foliations of the flow.  This type
of argument has been used since then in increasingly general contexts,
and has come to be known as the Hopf Argument. 

The core of the Hopf Argument is very simple.  Suppose that $\psi_t$
is a $C^\infty$ flow defined on a full measure subset $\Omega$  of a Riemannian manifold $V$, 
preserving a finite volume on $V$. For any $x\in \Omega$  one defines the stable and unstable sets: 
$$\cW^s(x) = \{x'\in \Omega : \lim_{t\to \infty} d(\psi_t(x), \psi_t(x')) = 0\}\quad \hbox{ and }\, \cW^u(x) = \{x'\in \Omega : \lim_{t\to -\infty} d(\psi_t(x), \psi_t(x')) = 0\}.$$
The stable (respectively unstable) sets partition $\Omega$ into measurable subsets. 

The first step in the Hopf Argument is to observe that for any continuous function $f\colon V\to \RR$ with compact support, the forward and backward upper Birkhoff averages
$$
f^s = \limsup_{T\to +\infty} \frac{1}{T}\int_0^T f \circ \psi_t\, dt\quad \hbox{ and }\,f^u = \limsup_{T\to -\infty} \frac{1}{T} \int_0^T f \circ \psi_t\, dt
$$
have the property that $f^s$
is constant on any stable set $\cW^s(x)$ and $f^u$ 
is constant on any unstable set $\cW^u(x)$.
Both functions $f^s$ and $f^u$ are evidently invariant under the flow $\psi_t$, and
the Birkhoff and von Neumann Ergodic Theorems (cf. \cite[Theorem 4.1.2 and Proposition 4.1.3]{KatokHasselblatt} imply that $f^s = f^u$ almost everywhere.
To show that $\psi_t$ is ergodic it suffices to show that $f^s$ is constant almost everywhere, for every continuous $f$ with compact support. The fundamental idea is to use the properties of the equivalence relation generated by the stable sets, the unstable sets, and the flow to conclude that $f^s = f^u$ must be constant.

In the next step in the Hopf Argument, one assumes some form of hyperbolicity of the flow, which will imply that the stable and unstable sets are in fact smooth manifolds.   In the original context of Hopf's argument,  $V=\Omega = T^1S$ is the unit tangent bundle of a compact, negatively curved surface $S$ and $\psi_t$ is the geodesic flow.  In this setting, the stable and unstable sets have a particularly nice description.  For almost every unit vector $v$, the stable and unstable Busemann  functions $b_v^s$ and $b_v^u$ are globally defined $C^\infty$ functions.  The stable and unstable sets are the orthogonal vectors to the level sets of these functions or equivalently the gradients of these functions on the level sets.  They are $C^\infty$, globally defined,  and for $\ast\in\{s,u\}$, the collection $$\cW^\ast:= \{\cW^\ast(v) : v\in T^1S\}$$ 
defines a $C^1$ foliation of $T^1S$. At each point $v\in T^1S$, the tangent space $T_v T^1S$ is spanned by the tangents to $\cW^s(v), \cW^u(v)$ and the direction $\dot\psi(v)$ of the flow.  A local argument in $C^1$ charts using Fubini's theorem shows that any $\psi_t$-invariant function that is almost everywhere constant along leaves of $\cW^s$ and $\cW^u$ must be locally almost everywhere constant, and hence globally almost everywhere constant, since $T^1S$ is connected.  In particular
the function $f^s$ is constant for any continuous, compactly supported $f$, and so $\psi_t$ is ergodic.

Hopf's original argument does not generalize immediately to geodesic flows for higher dimensional compact, negatively curved manifolds.  In this higher-dimensional setting, the stable and unstable foliations $\cW^s$ and $\cW^u$ exist, again arise from the   level sets of Busemann functions, and have $C^\infty$ leaves.  In general, however they fail to be $C^1$ foliations (except when the curvature is $1/4$-pinched) and so the argument using Fubini's theorem in local $C^1$ charts fails.

In the late 1960's Anosov \cite{An67} overcame this obstacle by proving that for any compact, negatively curved manifold, the foliations $\cW^s$ and $\cW^u$ are {\em absolutely continuous.} Absolute continuity, a strictly weaker property than $C^1$, is sufficient to carry out a Fubini-type argument to show that any $\psi_t$-invariant function almost everywhere constant along leaves of $\cW^s$ and $\cW^u$ is locally constant. See Section~\ref{s=generalsetup} for a more detailed discussion of absolute continuity.  Anosov thereby proved that the geodesic flow for any compact manifold of negative sectional curvatures is ergodic.

There is an extensive literature devoted to extending the Hopf Argument beyond the uniformly hyperbolic setting of geodesic flows on compact negatively curved  manifolds.  For smooth flows defined everywhere on compact manifolds, Pesin \cite{Pesin} introduced an ergodic theory of nonuniformly hyperbolic systems. In short, Pesin theory shows that if $\psi_t\colon V\to V$ preserves a finite volume and is nonuniformly hyperbolic, then almost everywhere the stable and unstable sets are smooth manifolds.  The family of stable manifolds is measurable and absolutely continuous in a suitable sense. 

 From Pesin theory, one deduces that a nonuniformly hyperbolic diffeomorphism of a compact manifold  has countably many ergodic components of positive measure.  More information about the flow can be used in some contexts to deduce ergodicity.  The obstruction to using the full Hopf Argument in this setting is that stable manifolds are defined only almost everywhere, and they may be arbitrarily small in diameter, with poorly controlled curvatures, etc.  

In a somewhat different direction than Pesin theory, Sinai \cite{Sinai} introduced methods for proving ergodicity of hyperbolic flows with singularities and applied them in his study of the $n$-body problem of celestial mechanics. Here the flow $\psi_t$ locally resembles the geodesic flow for a compact, negatively curved manifold, but globally encounters discontinuities and places where the norms of the derivatives $\|D\psi_t\|$ and $\|D^2\psi_t\|$  become unbounded.

 Introducing new techniques in the Hopf argument, Sinai was able to show that for several important classes of systems, including some billiards and flows connected to the $n$-body system, ergodicity holds.  These arguments have since been generalized to much larger classes of singular hyperbolic systems and
singular nonuniformly hyperbolic systems. 

In the singular nonuniformly hyperbolic setting, all aspects of Hopf's argument require careful revisiting.  The mere existence of local stable manifolds is a delicate matter and depends in a strong way on the  growth of the derivative of $\psi_t$ near the singularities. To give a sense of how delicate these issues can be,  we remark that:
\begin{itemize}
\item for compact surfaces of {\em nonpositive} curvature and genus $g\geq 2$, it is unknown whether the geodesic flow is always ergodic (even though it is always transitive);
\item there exist complete, finite volume surfaces of pinched negative curvature (but unbounded derivative of curvature) whose stable foliations are not even H\"older continuous \cite{BallmannBrinBurns};
\item for $C^1$ nonuniformly hyperbolic systems that are not $C^2$, stable sets can fail to be manifolds \cite{Pugh};
\item nonuniformly hyperbolic systems on compact manifolds can fail to be ergodic and can even have infinitely many ergodic components with positive measure \cite{DHP}. 
\end{itemize}
A general result providing for the existence and absolute continuity of local stable and unstable manifolds for singular, nonuniformly hyperbolic systems was proved by Katok-Strelcyn \cite{KS}.  We will use this work in an important way in this paper.

Returning to the context of the present paper, the WP geodesic flow is a singular, nonuniformly hyperbolic system.  To prove that it is ergodic, the first step is to verify the Katok-Strelcyn conditions to establish existence and absolute continuity of local stable and unstable manifolds. 
In particular, one needs to control the norm of the first two derivatives of the geodesic flow in a neighborhood of the boundary of  $\cM^1(S)$. 

To control the first derivative, we use the asymptotic expansions of Wolpert for the WP curvature and covariant derivative found in \cite{WolpertGeometry, WolpertExtension, WolpertUnderstanding}, combined with a careful analysis of the solutions to the WP Jacobi  equations.  This is the content of Theorem~\ref{t=firstderivest}. The precise estimates obtained by Wolpert appear to be essential for these calculations.

Since Wolpert's expansions of the WP metric are only to second order, and we need third order control to estimate the second derivative of the flow, we borrow ideas of McMullen in \cite{CTMKahler}. There is a {\em nonholomorphic} (in fact totally real) embedding of  $\Teich(S)$ into quasifuchsian space  $QF(S)$, under which the WP symplectic form has a holomorphic extension.  This holomorphic form is the derivative of a one-form that is bounded in 
the Teichm\"uller metric.  Using the Cauchy Integral Formula and a comparison formula between Teichm\"uller and WP metrics, one can then obtain bounds on all derivatives of the WP metric. This is the content of Proposition~\ref{p=derivcurvbounds}.  These bounds are adequate to control the second derivative of the geodesic flow, using the  bounds on the first derivative already obtained.

Once the conditions of \cite{KS} have been verified,  we are guaranteed the almost everywhere existence of absolutely continuous families $\cW^s$ and $\cW^u$  of  local stable and unstable manifolds.  Nonetheless these stable and unstable manifolds may not have uniform size.  At this point, we use negative curvature and another key property of the WP metric called {\em geodesic convexity} to show that in fact $\cW^s$ and $\cW^u$ have well-controlled uniform size.  

As a by-product of our arguments, we obtain that the WP Busemann function is $C^\infty$ for almost every tangent direction to $\Teich(S)$   (see Proposition~\ref{p=achoro}).  The local geometry of $\cW^s$ and $\cW^u$ is sufficiently nice that Hopf's original argument can be used with small modifications.  In particular, none of the more complicated local ergodicity arguments, such as the ``Hopf chains" developed by Sinai, are necessary. We also obtain positive, finite entropy of the WP flow using results of Katok-Strelcyn and Ledrappier-Strelcyn in \cite{KS}.

The paper does not quite follow the structure of this outline.  Rather than restricting to the special case of the WP metric, we instead develop an abstract criterion for ergodicity of the geodesic flow for an incomplete, negatively curved manifold.  This has the advantage of clarifying the  issues involved and also might allow for further applications. This is carried out in Section~\ref{s=generalsetup}, which may be read independently of the rest of the paper.  The remainder of the paper is devoted to setting up and verifying the conditions in Section~\ref{s=generalsetup} in the case of the Weil-Petersson metric.  

We remark that 
Pollicott and Weiss \cite{PW} gave a fairly complete outline of how to prove ergodicity for the Weil-Petersson metric in the cases $(g,n) = (1,1)$ and $(0,4)$.  
They say in the paper that the missing ingredients are the
bounds on the first and second derivatives of the geodesic flow, 
which are two of the major steps accomplished in this paper in the case of general $(g,n)$.

\subsection{The case of the punctured torus}

Several interesting features of the WP metric are already present in the simplest cases
 $(g,n) = (1,1)$  and $(0,4)$, where $S$ is the once-punctured torus or the four times punctured sphere.  In these cases, $\Teich(S)$ is the upper half space $\HH$ and  
%the Teichm\"uller metric $d_T$ is precisely  the hyperbolic metric, (cf. \cite{FarbMargalit}) 
$\cM(S)$ is the classical moduli space of elliptic curves $\HH/\PSL(2,\ZZ)$, which is a sphere with one puncture and two cone singularities of order $2$ and $3$.  

The mapping class group
$\MCG(S)$ is the modular group $\SL(2,\ZZ)$.  Due to the presence of torsion elements
in $\PSL(2,\ZZ)$, the space $\cM(S)$ is not a manifold, but the finite branched cover $\HH/\Gamma[k]$, for $k\geq 3$ is a manifold
\cite{Serre},  where $\Gamma[k]$ is the level-$k$ congruence subgroup
$$
\Gamma[k] = \{A\in \PSL(2,\ZZ)\,\mid\, A\equiv I \quad \mod \, k\}.
$$
The tangent bundle to $\Teich(S)$ is canonically identified with $\PGL(2,\RR)$.
% and the unit tangent bundle in the Teichm\"uller metric with $\PSL(2,\RR)$.

There are global coordinates $(\ell,\tau)$ in $\Teich(S)$, the so-called {\em Fenchel-Nielsen coordinates}, which have the asymptotic (first-order) expansions
$$\ell(z) \sim \frac{1}{\Im(z)}, \quad\text{ and } \tau(z) \sim \frac{\Re(z)}{\Im(z)},\quad\text{ as } \, \Im(z)\to\infty,
$$
and the WP form  has the first-order asymptotic expansion
$$
\omega_{WP} =  \frac{1}{2} d\ell\wedge d\tau \, \sim \, \frac{1}{\Im(z)^3} dz\wedge d\bar z ,\quad\text{ as } \,\Im(z) \to\infty.
$$
Since the complex structure on $\Teich(S)$ is the standard one on $\HH$, we obtain the expansion
$$
g_{WP}^2  \, \sim \, \frac{|dz|^2}{\Im(z)^3}.
$$
%By contrast,
%$$g_{T}^2 = \frac{|dz|^2}{\Im(z)^2} .$$
A neighborhood of the cusp in $\cM(S)$ is formed by taking the quotient of the points above the line $\Im(z) =\Im(z_0)$, for $\Im(z_0)$  sufficiently large, by the mapping class element $z\mapsto z+1$.  A model for this neighborhood is  the  surface of revolution for the curve $\{y = x^3 : x>0\}$ about the $x$-axis.    

From the form of the metric one can see the  incompleteness: a vertical  ray to the cusp at infinity starting at $\Im z=y_0$ has length $\sim 2y_0^{-1/2}\sim 2\ell^{1/2}$.  Moreover the curvature $K$ satisfies  $K\sim -\frac{3}{2\ell}\to-\infty$ as $\Im(z)\to \infty$.  These precise rates of divergence for the minimum sectional curvature hold as well in higher genus and will be crucial to our investigations.

Pollicott and Weiss \cite{PW} studied the model case of a negatively curved surface whose singularities coincide with a  surface of revolution for a polynomial and proved ergodicity of the geodesic flow in this case. 

\bigskip

\noindent{\bf Acknowledgments. }The authors express their appreciation to Scott Wolpert and Curt McMullen for many helpful conversations during the time this paper was being written. We also thank   Nikolai Chernov,  Benson Farb and Carlangelo Liverani for useful discussions and Ursula Hamenst\"adt for bringing our attention to the problem.

\section{Background on Teichm\"uller theory, Quasifuchsian space,  and Weil-Petersson geometry}

Much of the discussion in this section is based on McMullen's paper \cite{CTMKahler}.  Useful background can be found in \cite{Nag} and the course notes \cite{CTMnotes}.

\subsection{Riemann surfaces and tensors of type $(r,s)$}

We begin with some preliminary facts about Riemann surfaces. A Riemann surface is a topological surface equipped with an atlas of charts into $\CC$ with holomorphic transition maps. Suppose that $X$ is a Riemann surface of genus $g$ with $n$ punctures.  
We assume that $3g +n \geq 4$. Uniformization implies that $X$ is conformally equivalent to a quotient $\HH/\Gamma$, where $\HH$ denotes the upper half plane, and $\Gamma$ is a discrete subgroup of $PSL(2,\RR)$.  The hyperbolic metric $\tilde\rho$ on $\HH$ given by 
$$
\tilde\rho(z) =  \frac{|dz|}{\Im z}
$$
descends to a metric $\rho$ on $\HH/\Gamma$ of finite area, which is the unique Riemannian metric of constant curvature $-1$ on $X$ that induces the same conformal structure.

Denote by $\kappa$ the holomorphic cotangent bundle  and by $\kappa^{-1}$ the holomorphic tangent 
bundle of $X$, both of which are holomorphic complex line bundles over $X$.
For $r$ an integer, we denote by $\kappa^r$ the $|r|$-fold complex tensor product $\otimes^{|r|}\kappa$, if $r \geq 0$, and  $\otimes^{|r|}\kappa^{-1}$
if $r<0$.

A {\em tensor of type $(r,s)$} on $X$ is a section of the complex line bundle $\kappa^r\otimes \bar\kappa^s$ over $X$.
This leads to the construction of $L^p$ norms on the space of measurable $(r,s)$ tensors, defined as follows; for $\psi$ an $(r,s)$ tensor, and $p\geq 1$ we define:
$$
\|\psi\|_p := \left( \int_X  \rho^{2- p(r+s)} |\psi|^p  \right)^{1/p}\qquad \|\psi\|_\infty := \hbox{ess}\sup_X \rho^{-(r+s)}|\psi|.
$$
These norms will give rise to the Teichm\"uller ($p=1$) and WP ($p=2$) metrics on Teichm\"uller space, which we now define.

\subsection{Teichm\"uller and Moduli spaces}

A {\em marked complex structure} is a Riemann surface $X$ together with a homeomorphism $f\colon S\to X$, where $S$ is a fixed Riemann surface.
Given a marking surface $S$ of genus $g$ with $n$ punctures, we define the 
{\em Teichm\"uller space $\Teich(S)$} to be the set of equivalence classes of marked complex structures
$f\colon S\to X$, where $f_1\colon S\to X_1$ and $f_2\colon S\to {X_2}$ are equivalent if there is a conformal  map
$h\colon X_1 \to {X_2}$ isotopic to $f_2 f_1^{-1}$.  

  Uniformization gives an identification of $\Teich(S)$ with an open component of the representation variety 
of homomorphisms from $\pi_1(S)$ into the real Lie group $PSL(2,\RR)$, modulo conjugacy; this
identification gives $\Teich(S)$ a real analytic structure.  $\Teich(S)$ also carries a compatible complex analytic structure, which we shall
describe a little later.

The {\em mapping class group} $\MCG(S)$ is the set of equivalence classes of orientation preserving diffeomorphisms of $S$ modulo
isotopy, which forms a group under composition.  $\MCG(S)$ acts properly by diffeomorphisms of $\Teich(S)$ via precomposition with
the marking homeomorphisms $f\colon S\to X$; the quotient $\cM(S) = \Teich(S)/\MCG(S)$ is easily seen to be the moduli space of 
Riemann surfaces homeomorphic to $S$, modulo conformal equivalence.  The $\MCG(S)$-stabilizer of any 
point in $\Teich(S)$ is finite. In denoting an element of $\Teich(S)$, we will often omit the marking given by the equivalence class
of maps $f\colon S\to X$ and refer only to the target Riemann surface $X$.
  We do this because the tangent space and cotangent spaces at any point do not depend on the marking, but only on the target $X$.

We review the definition of the Weil-Petersson norms 
on the tangent and cotangent spaces $T_X\Teich(S)$ and $T^\ast_X\Teich(S)$
at a point $X\in \Teich(S)$.
%These definitions are made by identifying the cotangent and tangent spaces
%with appropriate spaces of differentials on $X$, as we now explain.  
An {\em integrable meromorphic quadratic differential} on $X$ is a  tensor of type $(2,0)$ 
that has a local representation
of the form $q(z) dz^2$ where $q(z)$ is holomorphic on $X$ and has at most simple poles at the punctures.   We define $Q(X)$ to be the vector space of integrable meromorphic quadratic differentials $\phi$ on $X$.

A {\em Beltrami differential} on $X$ is a measurable tensor of type $(-1,1)$, which has a local representation
of the form $b(z) d\bar z/ dz $.  Note that the product of a Beltrami differential with a quadratic
differential is a $(1,1)$-tensor. Let $M(X)$ be the vector space of all measurable Beltrami differentials $\mu$ on $X$ with the property that
$\int_X |\phi\mu| < \infty$, for every $\phi\in Q(X)$. We then have a natural complex pairing of the space $M(X)$ with $Q(X)$ given by  
\begin{equation}
\label{eq:pair}
\langle \phi,\mu \rangle = \int_X \phi\mu \quad \hbox{for } \phi \in Q(X),\, \mu\in M(X).
\end{equation}
In view of the fact that elements of $Q(X)$ have finite
$L^p$ norm for every $1\leq p\leq \infty$, it follows that elements of 
$M(X)$ are precisely those Beltrami differentials $\mu$ on $X$ of finite $L^q$ norm, for
$1\leq q \leq \infty$.

We have the fundamental isomorphisms of vector spaces:
$$T_X\Teich(S) \cong M(X)/Q(X)^\perp
\quad\hbox{ and }\quad 
T^\ast_X\Teich(S) \cong Q(X),
$$
where $Q(X)^\perp = \{\mu\in M(X) : \langle \mu, \phi \rangle = 0,\, \forall \phi\in Q(X)\}$.

Having described these identifications, we now can define the WP norm.
The {\em Weil-Petersson metric} on $T^\ast _X\Teich(S)$ is defined  by the $L^2$ norm:
$$\|\phi\|_{WP} = \|\phi\|_2 = \left(\int_X \rho^{-2} |\phi|^2\right)^{1/2}.
$$
Note that the definition of the WP metric involves both conformal and hyperbolic data from $X$; this feature
makes the WP metric somewhat tricky to work with.  On the other hand, the hyperbolic input from the metric 
$\rho$ leads to the delicate and beautiful connections between the WP metric and hyperbolic geometry and dynamics
discussed in the introduction.

The WP norm on the tangent space $T_X\Teich(S)$ is induced by the pairing (\ref{eq:pair}) via the formulae:
$$
 \|v\|_{WP} = \sup_{\phi\in Q(X),\,\,\|\phi\|_{WP}=1 } \Re (\langle\phi, \mu\rangle ),
$$
for any $\mu \in M(X)$ representing the tangent vector $v\in T_X\Teich(S)$.  

%By duality of the appropriate $L^p$ norms, it follows that
%the Teichm\"uller norm is induced by the $L^\infty$ norm on $M(X)$ and the WP norm is induced by the $L^2$ norm:
%$$
%\|v\|_T = \|\mu\|_\infty ; \quad \|v\|_{WP} = \|\mu\|_2 = \left(\int_X \rho^{2} |\mu|^2\right)^{1/2} ,
%$$
%for any $\mu \in M(X)$ representing $v\in T_X\Teich(S)$. The Teichm\"uller and WP norms are both invariant under the  %derivative action of $\MCG(S)$.

%A simple comparison between the WP and Teichm\"uller norms follows from the 
%Cauchy-Schwarz and Gauss-Bonnet Theorems (see Prop 2.4 in \cite{CTMKahler}):
%$$\|v \|_{WP}\leq |2\pi \chi(S)|^{1/2} \|v\|_T,
%$$
%for any $X\in \Teich(S)$ and $v\in T_X\Teich(S)$.
%We prove in Section~\ref{s=comparison} a complementary bound: there exists $C>0$ such that
%$$\|v\|_{WP}\geq C\uell(X)^{-1/2}\| v \|_{T},$$
%for any $v \in T_X\Teich(S)$, where $\uell(X)$ denotes the hyperbolic length of the shortest closed
%geodesic on $X$, also known as the {\em systole} of $X$.
%\footnote{Since this paper was written, Wolpert has computed a precise bound on the value of $C$ in this estimate}. 
%This lower bound will be used to estimate from above the derivatives of the WP metric.

\subsection{The bundle of projective structures on $S$}\label{ss=projectivebundle}

A {\em projective structure} on a surface $X$ is an atlas of charts into 
$\CC$ whose overlaps are M\"obius transformations (elements of $PSL(2,\CC)$); note that
a projective structure determines a unique complex structure. Fix as above a Riemann surface $S$ of
genus $g$ with $n$ punctures.
A {\em marked projective structure} is a homeomorphism $f\colon S\to X$, where $X$
is endowed with a projective structure; we say that two marked structures $f_1\colon S\to {X_1}$ and $f_2\colon S\to {X_2}$
are equivalent if there is a projective isomorphism from ${X_1}$ to ${X_2}$ homotopic to $f_2f_1^{-1}$.  
Denote by $\Proj(S)$ the space of equivalence classes of projective structures marked by $S$.

It is a classical fact that $\Proj(S)$ has the structure of a complex manifold  that
arises from its embedding into the representation variety of homomorphisms from $\pi_1(S)$ into
$PSL(2,\CC)$, modulo conjugacy (see \cite{Hubbard}). The map that assigns to each marked projective structure the 
compatible marked conformal structure defines a fibration $\pi\colon \Proj(S)\to \Teich(S)$.  The fiber $\Proj_X(S)$ over $X$ is
an affine space modelled on $Q(X)$. In particular there is a well-defined difference $\beta_1-\beta_2\in Q(X)$, for $\beta_1,\beta_2\in \Proj(S)$, which defines a holomorphic map  from $\Proj(S)\times \Proj(S)$ to $Q(X)$.

%In the next subsection we describe how to define a complex structure on $\Teich(S)$ that
%gives $\Proj(S)$ the structure of a holomorphic fiber bundle, and a holomorphic identification
%of $\Proj(S)$ with the cotangent bundle $T^\ast \Teich(S)$. We shall also describe how to construct a holomorphic section %$\sigma_{QF}$ 
%of this bundle that arises from a process called {\em quasifuchsian uniformization}.
%For now we note that
%the Fuchsian uniformization $X\cong \HH/\Gamma$ defines a section
%$\sigma_F\colon \Teich(S)\to \Proj(S)$.   This section is real analytic but not
%holomorphic. 

\subsection{Quasifuchsian space}\label{ss=quasifuchsian}
Let $S=\HH/\Gamma$ be a hyperbolic Riemann surface with $\Gamma < PSL(2,\RR)$, and denote by $\bar S$ the hyperbolic Riemann surface
$\LL/\Gamma$, where $\LL$ is the lower half plane. Since $\Gamma$ is a Fuchsian group, it acts on the Riemann sphere
$\hat\CC$ fixing $\HH$, $\LL$ and the real axis/circle at infinity $\RR_\infty = \hat\CC \setminus (\HH \cup \LL)$.
Following McMullen \cite{CTMKahler}, we define {\em quasifuchsian space} $QF(S)$ to be the product:  
$$QF(S)=\Teich(S)\times \Teich(\bar S).$$

Then $QF(S)$ parametrizes marked quasifuchsian groups equivalent to $\Gamma(S)$. A quasifuchsian group is a Kleinian group 
$\Gamma(X,Y)$  with a domain of discontinuity $\Omega(X,Y)$ consisting of two components whose quotients by $\Gamma(X,Y)$ are $X$ and $Y$ respectively.
% more precisely we have:
%\begin{theorem}[Bers Simultaneous Uniformization]\label{t=Bers} 
%For each pair $([f\colon S\to {X}], [g\colon \overline{S} \to {Y}]) \in QF(S)$, there exist a 
%quasiconformal homeomorphism $\phi\colon \hat\CC \to \hat\CC$, a Kleinian group $\Gamma({X},{Y})<PSL(2,\CC)$,
%and a homomorphism $h\colon \Gamma \to \Gamma({X},{Y})$ (all unique up to the action of $PSL(2,\CC)$ by conjugation)
%with the following properties:
%\begin{enumerate} 
%\item for each $\gamma\in \Gamma$, we have $\phi\circ \gamma = h(\gamma)\circ \phi$;
%\item $\phi$ sends $(\HH\cup \LL, \RR_\infty)$ to $(\Omega({X},{Y}),\Lambda({X},{Y}))$, 
%where $\Omega({X},{Y}) $ is the domain of discontinuity of $\Gamma$ and $\Lambda({X},{Y})\subset \hat\CC$ is the limit sets %and is a  quasicircle;
%\tem there is a conformal isomorphism $\Omega({X},{Y})/\Gamma({X},{Y}) \cong ({X}\cup {Y})$ making 
%the following diagram commute:

%\[
%\begin{xy}
%\xymatrixcolsep{3pc}
%\xymatrixrowsep{2pc}
%\xymatrix{ (\HH\cup \LL)\ar[r]^{\phi}\ar[d]_{\mod \,\, \Gamma} & \Omega({X},{Y})\ar[d]^{\mod \,\,\Gamma({X},{Y})} \\
%(S\cup\overline S)\ar[r]^{f\cup g} &  ({X}\cup {Y}) \\
%}
%\end{xy}
%\]
%\end{enumerate}

We thus have a  ``quasifuchsian uniformization" map
$$
\sigma\colon \Teich(S)\times \Teich(\overline S) \to \Proj(S)\times \Proj(\overline S)
$$
that sends $(X,Y)$ to the projective structures on $X$ and $Y$ inherited from $\Omega(X,Y)$ from the action of $\Gamma(X,Y)$. The map $\sigma$ is a section of the bundle $ \Proj(S)\times \Proj(\overline S) \to QF(S)$.
We write:
$$\sigma(X,Y) = (\sigma_{QF}(X,Y), \overline\sigma_{QF}(X,Y)).
$$
We define the {\em Fuchsian locus $F(S)$} to 
be the image of $\Teich(S)$ under the antidiagonal embedding $\hat\alpha(X) = (X,\overline{X})\in QF(S)$.  
%Observe that $\sigma_{QF}(X,\bar X) = \sigma_F(X)$.

The complex structure on $\Teich(S)$ is then defined via the {\em Bers embedding:} fixing $X\in \Teich(S)$, we define
$\beta_X\colon \Teich(S) \to Q(X)$ by 
$$
\beta_X(Y) = \sigma_{QF}(X,\bar Y) - \sigma_F(X).
$$
The map $\beta_X$ is an embedding, and the pullback of the complex structure on $Q(X)$ gives a complex structure on $\Teich(S)$ that
is independent of $X$ (that is, two different $X$s give isomorphic structures).  Recall that $Q(X)$ is a Banach space when endowed with any
$L^p$ norm.  
%A fundamental result of Nehari implies that the image of $\Teich(S)$ under $\beta_X$ is a bounded subset of
%$Q(X)$ {\em in the $L^\infty$ norm}, where the bound is independent of $X$.  This fact is used in a crucial way 
%in McMullen's paper and in ours.

We have defined a complex structure on $\Teich(S)$, which induces a conjugate complex structure on $\Teich(\bar S)$.  The
complex structure on $QF(S)$ is defined to be the product complex structure.  The Fuchsian locus $F(S)$ is then a totally real
submanifold of $QF(S)$. It can be checked that the fibration $\Proj(S)\to \Teich(S)$ is holomorphic with respect to these
structures.  
%Furthermore, the naturally-defined quasifuchsian uniformization section $\sigma\colon QF(S)\to \Proj(S)\times \Proj(\bar S)$ is holomorphic. 
 Hence, for a fixed
$Y\in \Teich(S)$, the map $X\mapsto \sigma_{QF}(X,\bar Y)$ gives a holomorphic section of $\Proj(S)$ over $\Teich(S)$; this
section gives an isomorphism between  the cotangent bundle $T^\ast\Teich(S)$ and an open subset of $\Proj(S)$. 

We will use the quasifuchsian uniformization section $\sigma$ in a crucial way to estimate higher derivatives of the $WP$ metric in Section~\ref{s=higherorder}.  We record here the properties that we will use.

\begin{theorem} The holomorphic section $\sigma$ satisfies the following properties: 
\label{thm:quasi}
\begin{enumerate}
\item $\sigma_{QF}(X,\overline X) = \sigma_F(X)$;
%where $\sigma_F$ is the Fuchsian uniformization section defined above; 
\item for any $Y, Z \in \Teich(\bar S)$, the map $X\mapsto \sigma_{QF}(X, Y) - \sigma_{QF}(X, Z)$ defines a bounded holomorphic $1$-form on $\Teich(S)$ in the $L^\infty$ norm;
%\item there exists $C\geq 1 $ such that for any $X\in \Teich(S)$ and $Y, Z \in \Teich(\bar S)$,
%$$\|\sigma_{QF}(X,Y) - \sigma_{QF}(X,Z)\|_T \leq C,
%$$
%where $\|\cdot\|_T$ is the Teichm\"uller norm on $T^\ast_X \Teich(S)$.
%\end{itemize}
\item for each $Z\in \Teich(\bar{S})$ the $1$-form 
%$\theta_{WP}$ be
%the $(1, 0)$-form on Teich(S) given by
%\begin{eqnarray*}
$\theta_{WP}(X)  =  \sigma_F(X) -  \sigma_{QF}(X, Z )=-\beta_X(\bar{Z})$
satisfies 
$d(i\theta_{WP}) =\omega_{WP}$.
\end{enumerate}
\end{theorem}
The boundedness of the $1$-form in (2) follows from Nehari's bound  (see Theorems 2.2 in \cite{CTMKahler}).
% and a simple estimate which bounds the Teichm\"uller  $L^1$ norm in terms of the $L^\infty$ norm in Nehari's bound. 
The last statement 
is due to  McMullen [\cite{CTMKahler}, Theorem 7.1]
%\label{t=ctmmain} Fix $Z\in\Teich(\bar{S})$, and let 
%\end{theorem}

\subsection{Fenchel-Nielsen coordinates}

Continue to denote by $S$ a marked Riemann surface of genus $g$ with $n$ punctures.  We define here a natural system of global coordinates
on $\Teich(S)$, called Fenchel-Nielsen coordinates, in which the K\"ahler form $\omega_{WP}$ takes a  simple form. 

Recall that a curve in $S$ is {\em nonperipheral} if it is not homotopic to a  loop surrounding a single puncture. 
A {\em pants decomposition} of  $S$ is a collection $P$ of $3g-3+n$ pairwise disjoint, homotopically nontrivial, nonperipheral and homotopically distinct simple closed curves. The complement of these curves is a collection of surfaces called {\em pairs of pants}.  Topologically, a pair of pants is a three-times punctured sphere.   A pair of pants has one of three types of conformal structure depending on whether each puncture is locally modelled on the punctured plane or on the complement of a closed disk in the plane, in which case we say that the boundary component is a circle.  A pair of pants with $j$ boundary circles  has a $j$-dimensional space of hyperbolic structures, parametrized by the 
hyperbolic lengths of the boundary circles.

We introduce notation that will  be used throughout the paper.  If $f\colon S\to X$ is a marked Riemann surface and $\alpha$ is a homotopically nontrivial, nonperipheral, simple closed curve in $S$, we denote by  $\ell_\alpha(X)$ the hyperbolic length in $X$ of the unique geodesic in the homotopy class of $f_\ast[\alpha]$. This geodesic length function is intimately connected with the WP metric and is used to define Fenchel-Nielsen coordinates.

Fix a pair of pants decomposition $P=\{\alpha_1,\ldots,\alpha_{3g-3+n}\}$ of $S$. The Fenchel-Nielsen coordinates 
$$(\ell_\alpha,\tau_\alpha)_{\alpha\in P}\colon \Teich(S) \to  (\RR_{>0}\times \RR)^{3g-3+n}$$
determined by $P$  are defined as follows. For $f\colon S\to X$ a marked Riemann surface and
$\alpha\in P$, we define $\ell_\alpha(X)$ to be the geodesic length as above and
$\tau_\alpha(X)$ to be the twist parameter, which records the 
relative displacement in how the pairs of pants are glued together 
along $\alpha$ to obtain the hyperbolic metric on 
$X$; more precisely, a full Dehn twist about the curve $\alpha$ changes $\tau_\alpha$ by the amount $\ell_{\alpha}$.  
One must adopt a convention for how this relative displacement $\tau$ is measured, as it is intrinsically
only well-defined up to a constant, but this does not introduce any serious issues. These give global coordinates on $\Teich(S)$ a fact which shows that $\Teich(S)$ is homeomorphic to $\R^{6g-6+2n}$.
 
The Fenchel-Nielsen coordinates are natural with respect to the WP metric. Wolpert \cite{WolFN}  proved that for any pants decomposition $P$, we have
$\omega_{WP} = \frac12\sum_{\alpha\in P} d\ell_\alpha \wedge d\tau_\alpha$.  An ingredient in the proof  of this
formula is the important  fact that the vector field $\partial/\partial \tau_\alpha$, which generates
the Dehn twist flow about $\alpha$, is the symplectic gradient of the Hamiltonian function $\frac12\ell_\alpha$:
$$
\frac12 d\ell_\alpha =  \omega_{WP}( \cdot , \frac{\partial}{\partial \tau_\alpha}),
$$
or equivalently
$$ \grad\ell_\alpha  = - 2 J \frac{\partial}{\partial \tau_\alpha}.
$$
This fundamental relationship is the starting point for many of Wolpert's deep asymptotic expansions for the WP metric,
which we discuss in more detail in Section~\ref{s=firstderivWP}.

\subsection{The Deligne-Mumford compactification of moduli space}\label{ss=dmcompact}

As mentioned earlier,  $\Teich(S)$ is incomplete with respect to the WP distance \cite{Wol1}.  
This occurs precisely because  it is possible to shrink a simple closed curve $\alpha$  to a point and leave Teichm\"uller space along a WP geodesic  in finite time --- 
indeed, the time it takes  is on the order of $\ell_\alpha^{1/2}$. 
%It turns out that pinching one or more closed geodesics to a point is 
%the {\em only} way to leave every compact set of $\Teich(S)$ in finite time along a WP geodesic;
This fact allows one to prove   \cite{Mas} that the  completion of $\Teich(S)$ is the  {\em augmented Teichm\"uller space},   denoted $\bar\Teich(S)$.  The mapping
class group $\MCG(S)$ acts on $\bar\Teich(S)$ and the quotient $\bar\cM(S)$ is the {\em Deligne-Mumford compactification} of the  moduli space $\cM(S)$ and gives the completion on the quotient.

Augmented Teichm\"uller space $\bar\Teich(S)$ is obtained by adjoining lower-dimensional Teichm\"uller spaces of noded Riemann surfaces, which gives it the structure of a stratified space.  The combinatorics of this stratification are encoded by a symplicial complex $\C(S)$ called the 
{\em curve complex}.  We review this construction here.

We first define the curve complex $\C(S)$, which is a $3g-4+n$ dimensional simplicial complex.  The 
vertices of $\C(S)$ are homotopy classes of homotopically nontrivial, nonperipheral simple closed curves on $S$.  
We join two vertices by an edge if the corresponding pair of curves has disjoint representatives. More generally,
a $k$ simplex $\sigma\in \C(S)$ consists of $k+1$  distinct vertices that have disjoint representatives. 
We note that in the sporadic cases of the punctured torus $(g,n) = (1,1)$ and $4$-times punctured sphere $(g,n)=(0,4)$,
$\C(S)$ is just an infinite discrete set of vertices, since there do not exist disjoint  homotopically distinct curves on the underlying
surface $S$.
  Except in these sporadic cases,  $\C(S)$ is a connected locally infinite complex.\footnote{In  the sporadic cases there is more than one possible definition of $\C(S)$; in another, very standard definition in these cases, 
one adds edges joining curves that intersect minimally (once in the case of the torus and twice in the case of the sphere). The resulting $1$-complex is the Farey graph in both cases.}   
Note that a maximal simplex in $\C(S)$ defines a  pants decomposition of $S$.
The mapping class group $\MCG(S)$ acts on $\C(S)$.

 A {\em noded Riemann surface} is a complex space with at most isolated
singularities, called nodes, each possessing a neighborhood biholomorphic
to a neighborhood of $(0, 0)$ in the curve $$\{(z,w)\in \CC^2: zw=0\}.$$
Removing the nodes of a noded Riemann surface $Y$ yields a (possibly
disconnected) punctured Riemann surface, which we will usually denote by $\hat Y$.
The components of $\hat Y$ are called the {\em pieces} of
$Y$. 

Given a simplex $\sigma \in C(S)$, a {\em marked noded Riemann surface} 
with nodes corresponding to  $\sigma$ is a noded Riemann surface $X_\sigma$ equipped with a continuous
mapping $f \colon S \to X_\sigma$
so that $f\vert_{S\setminus \sigma}$ is a homeomorphism to $\hat X_\sigma$. 
Two marked noded Riemann surfaces $[f_1 \colon S\to X^1_\sigma]$ and $[f_2 \colon S\to X^2_\sigma]$ are equivalent 
if there is a biholomorphic node preserving map $h\colon X^1_\sigma\to X^2_\sigma$   such that  $f_1\circ h$ is isotopic to $f_2$.
We denote by $\cT_\sigma$ the set of equivalence classes
$[f\colon S\to X_\sigma ]$ with nodes at $\sigma$.  We adopt the convention that when $\sigma=\emptyset$ then $\cT$
is the Teichm\"uller space $\Teich(S)$ of unnoded surfaces. Then  the augmented Teichm\"uller space 
$$
\bar\Teich(S)= \cT \cup \bigcup_{\sigma \in \
C(S)} \cT_\sigma.$$
(The space $\bar\Teich(S)$  should not be confused with $\Teich(\bar S)$, which was introduced in $\S$\ref{ss=quasifuchsian}.)  

\medskip

\noindent{\bf Notational convention.} If  the topological type of the surface $S$ is fixed,   $\bar\cT$ will  denote the augmented space $\bar\Teich(S)$.  We also denote by $\partial \cT$ the boundary $\bar\cT\setminus \cT$.  We denote by $\pi:T\cT\to \cT$ the natural projection.
As with the elements of $\Teich(S)$,  we will frequently
abuse notation and omit the marking when referring to an element of $\bar\cT$;
% hence
%``let $X_\sigma\in \bar\Teich(S)$" is shorthand for ``let  $[f\colon S\to X_\sigma]$ be an element
%of $\bar\Teich(S)$."

\medskip

To describe a neighborhood of a point $[f\colon S\to X_\sigma]$ in $\bar\Teich(S)$, 
we give coordinates adapted to the simplex $\sigma$. For any such $\sigma$, let
$P$ be the a maximal simplex in $\C(S)$  (pants decomposition) containing $\sigma$,
and let $(\ell_\alpha,\tau_\alpha)_{\alpha\in P} $ be the corresponding Fenchel-Nielsen coordinates on $\Teich(S)$.
Then the {\em extended Fenchel-Nielsen coordinates} for $P$ are obtained by allowing
the lengths $\ell_\alpha$ to range in $\RR_{\geq 0}$ and taking the quotient by identifying
$(0,t)$ with $(0,t')$ in each $\RR$ factor corresponding to the curves in $\sigma$. 
 
This also defines a topology on $\bar\Teich(S)$.  We note that the space is not locally compact. A neighborhood of a noded surface allows for the twists  $\tau_\alpha$ corresponding to the curves $\alpha\in \sigma$ to be arbitrary real numbers.

\section{Background on the geodesic flow}\label{s=background}

Let $M$ be a  Riemannian manifold. As usual $\langle v,w\rangle$  denotes the inner product of two vectors and $\nabla$ is the Levi-Civita connection defined by the Riemannian metric. It is the unique connection that is symmetric and compatible with the metric.
%: if $X$, $Y$ and $Z$ are vector fields, and $[\cdot,\cdot]$ denotes the Lie bracket, we have
 %$$
 %\nabla_XY - \nabla_YX = [X,Y]
 %$$
%and
 %$$
 %\nabla_X\langle Y,Z \rangle = \langle \nabla_XY, Z \rangle + 
   %                           \langle Y, \nabla_XZ \rangle.
 %$$

The covariant derivative along a curve $t \mapsto c(t)$ in $M$ is denoted by $D_c$, $\frac D{dt}$ or simply $'$ if it is not necessary to specify the curve; if $V(t)$ is a vector field along $c$ that extends to a vector field $\widehat V$ on $M$, we have 
 $$
 V'(t) = \nabla_{\dot c(t)} \widehat V.
 $$
Given a smooth map $(s,t) \mapsto \alpha(s,t)$, we let $\frac D{\partial s}$ denote covariant differentiation along a curve of the form $s \mapsto \alpha(s,t)$ for a fixed $t$. Similarly  $\frac D{\partial t}$ denotes
covariant differentiation along a curve of the form $t \mapsto \alpha(s,t)$ for a fixed $s$. The symmetry of the Levi-Civita connection means that 
 $$
 \frac D{\partial s}\frac{\partial\alpha}{\partial t}(s,t) = \frac D{\partial t}\frac{\partial\alpha}{\partial s}(s,t)
 $$
 for all $s$ and $t$.

 The curve  $c$ is a geodesic if it satisfies the equation $D_c \dot c(t) = 0$. Since this equation is a first order ODE in the variables $(c,\dot c)$, a geodesic is uniquely determined by its initial tangent vector. Geodesics have constant speed, since  we have $\dfrac d{dt}\langle \dot c(t), \dot c(t)\rangle = 2\langle\dot D_c\dot c(t), \dot c(t) \rangle = 0$ if $c$ is a geodesic.
 
The Riemannian curvature tensor $R$ is defined by
$$
R(A,B)C = (\nabla_A\nabla_B - \nabla_B\nabla_A - \nabla_{[A,B]})C.
$$
The sectional curvature of the $2$-plane spanned by vectors $A$ and $B$ is defined by
$$
K(A,B) = \frac{\langle R(A,B)B,A \rangle}{\|A\wedge B\|^2}.
$$

The action of the Levi-Civita connection extends to covectors and tensors in such a way that the product rule holds. In particular
$$
(\nabla_W R)(X,Y)Z = \nabla_W\left(R(X,Y)Z\right) - R(\nabla_WX,Y)Z - R(X,\nabla_WY)Z - R(X,Y)\nabla_WZ.
$$

Similarly the second derivative $\nabla^2_{X,Y}T$ of a tensor $T$ is defined by the product rule formula
 $$
 \nabla_X(\nabla_YT) = \nabla^2_{X,Y}T + \nabla_{\nabla_X Y}T.
 $$
 We will use this later in the case $T = R$. If $T$ is a vector field $Z$,
a short calculation using the symmetry of the Levi-Civita connection yields
 $$
 \nabla^2_{X,Y}Z - \nabla^2_{Y,X}Z = R(X,Y)Z.
 $$

\subsection{Vertical and horizontal subspaces and the Sasaki metric} 

The tangent bundle $TTM$ to $TM$ may be viewed as a bundle over $M$
in three natural ways shown in the following commutative diagram:
\[
\begin{xy}
\xymatrixcolsep{4pc}
\xymatrixrowsep{3pc}
\xymatrix{
TTM\ar[r]^{D\pi_M}\ar[d]^\kappa\ar[dr]^{\pi_{TM}\circ\pi_M}&TM\ar[d]^{\pi_M} \\
TM\ar[r]^{\pi_M}&M \\
}
\end{xy}
\]
The first is via the composition of  the natural bundle projections $\pi_{TM} :TTM\to TM$ and $\pi_M : TM \to M$. The second is via the composition of the derivative map $D\pi_M : TTM \to TM$ with $\pi_M$.
The third involves a map $\kappa:  TTM \to TM$, often called the {\em connector map}, which is  determined by the Levi-Civita connection. If $\xi \in TTM$ is tangent at $t = 0$ to a curve $t \mapsto V(t)$ in $TM$  and $c(t) =  \pi_M( V(t))$ is the curve of footpoints of the vectors $V(t)$, then
 $$
 \kappa(\xi)  = D_cV(0).
 $$

The  {\em vertical subbundle} is the subbundle $ker(D\pi_M)$.  It is naturally identified with $TM$ via the map $\kappa$.  The  {\em horizontal subbundle} is the subbundle $ker(\kappa)$.  It is naturally identified with $TM$ via the map $D\pi_M$ and is transverse to the vertical  subbundle. If $v \in T_pM$, we may
identify $T_vTM$ with $T_pM \times T_pM$ via the map $D\pi_{M} \times \kappa : TTM \to TM\times TM$.

Each element of $T_vTM$ can thus be represented uniquely by a pair
$(v_1,v_2)$ with $v_1\in T_pM$ and $v_2\in T_pM$.  
Put another way, every element $\xi$ of $T_vTM$ is tangent to a curve $V\colon (-1,1) \to TM$ with $V(0) = v$.   Let
$c = \pi_M\circ V\colon (-1,1) \to M$ be the curve of basepoints of $V$ in $M$.  Then $\xi$ is represented
by the pair
 $$
 (\dot c(0),  D_cV(0))  \in T_pM\times T_pM.
 $$
These coordinates on the fibers of $TTM$ restrict to coordinates on $TT^1M$. 

Regarding $TTM$ as a bundle over $M$ in this way
gives rise to a natural Riemannian metric on $TM$, called the {\em Sasaki metric}.  In this metric, the
inner product of two elements $(v_1,w_1)$ and $(v_2,w_2)$ of $T_vTM$ is defined:
$$
\langle (v_1,w_1) , (v_2,w_2) \rangle_{Sas} = \langle v_1 , v_2 \rangle  + \langle w_1 , w_2 \rangle .
$$
This metric is induced by a symplectic form $\omega$ on $TTM$; for vectors $(v_1,w_1)$ and 
$(v_2,w_2)$ in $T_vTM$, we have:
$$
\omega((v_1,w_1) , (v_2,w_2)) = \langle v_1,w_2 \rangle - \langle w_1,v_2 \rangle.
$$
This symplectic form is the pull back of the canonical symplectic form on the cotangent bundle $T^*M$ by the map from $TM$ to $T^*M$ induced by identifying a vector $v \in T_pM$ with the linear function $\langle v, \cdot\rangle$ on $T_pM$.

Sasaki \cite{Sasaki} showed that the fibers of the tangent bundle are totally geodesic submanifolds of $TTM$ with the Sasaki metric. A  parallel vector field along a geodesic of $M$ (viewed as a curve in $TM$) is a geodesic of the Sasaki metric. Such a geodesic is orthogonal to the fibers of $TM$. If $v \in T_pM$ and $v' \in T_{p'}M$, we can join them by first parallel translating $v$ along a geodesic from $p$ to $p'$ to obtain $w \in T_{p'}M$ and then moving from $w$ to $v'$ along  a line in  $T_{p'}M$. If $v'$ is close to $v$, we can choose the geodesic so that its length is $d(p,p')$. 
%Then $d_{Sas}(v,v')$ is the length of the hypotenuse of a small right angled triangle whose other two sides have lengths $d(p,p')$ and $\|w - v'\|$.
%Topogonov's comparison theorem \cite[Theorem 2.2]{CheegerEbin} tells us that if $-k^2$ is a lower bound for the curvature of the Sasaki metric in a neighborhood of $v$, then $d_{Sas}(v,v')$ is less than the length of the hypotenuse of a right triangle in a space form of constant curvature $-k^2$ with sides of length  $d(p,p')$ and $\|w - v'\|$
%meeting at the right angle. 
It follows easily from Topogonov's comparison theorem \cite[Theorem 2.2]{CheegerEbin} that
 $$
 d_{Sas}(v,v') \asymp   d(p,p') + \|w - v'\|,
 $$
 as $v' \to v$, where the rate of convergence is controlled by the curvatures of the Sasaki metric in a neighborhood of $v$. The notation $a\asymp b$, here and in the rest of the paper, means that the ratios
$a/b$ and $b/a$ are bounded from above by a constant.  In this case the constant is $2$.

\subsection{The geodesic flow and and Jacobi fields}\label{SSJfields}
For $v \in TM$ let $\gamma_v$ denote the  unique geodesic $\gamma_v$ satisfying $\dot\gamma_v(0) = v$. The geodesic flow $\varphi_t : TM\to TM$ is defined by
 $$
 \varphi_t(v) = \dot\gamma_v(t),
 $$
 wherever this is well-defined. The geodesic flow is always defined locally.  Since the geodesic flow is Hamiltonian, it  preserves a natural volume form on $T^1M$ called the Liouville  volume form.  When the integral of this form is finite, it induces a unique probability measure on $T^1M$ called the {\em Liouville measure} or {\em Liouville volume}.

Consider now a one-parameter family of geodesics, that is a map $\alpha:(-1,1)^2\to M$ with the property that
$\alpha(s,\cdot)$ is a geodesic for each $s\in (-1,1)$.  Denote by $J(t)$ the vector field  
$$
J(t) = \frac{\partial \alpha}{\partial s}(0,t)
$$
along the geodesic $\gamma(t) = \alpha(0,t)$.
Then $J$ satisfies the {\em Jacobi equation:}
$$
J'' + R(J,\dot\gamma)\dot\gamma = 0,
$$
in which $'$ denotes covariant differentiation along $\gamma$.
Since this  is a second order linear ODE,
the pair of 
vectors $(J(0), J'(0)) \in T_{\gamma(0)}M \times T_{\gamma(0)}M$ uniquely determines the vectors
$J(t)$ and $J'(t)$ along $\gamma(t)$. A vector field $J$ along a geodesic $\gamma$ satisfying the Jacobi equation is called a {\em Jacobi field}.

The pair $(J(0),J'(0))$ corresponds in the manner described above to the tangent vector at $s= 0$ to the 
curve $V(s) = \frac{\partial \alpha}{\partial t}(s,0)$. To see this, note that $V(s)$ is a vector field along the curve $c(s)=\alpha(s,0)$, so $V'(0)$ corresponds to the pair
 $$
  (\dot c(0), D_c\frac{\partial \alpha}{\partial t}(s,0))
    = (J(0), \frac D{\partial s}\frac{\partial \alpha}{\partial t}(s,0))
    = (J(0), \frac D{\partial t}\frac{\partial \alpha}{\partial s}(s,0))
    = (J(0),J'(0)).
  $$
  In the same way one sees that $(J(t),J'(t))$ corresponds to the tangent vector at $s = 0$ to the curve
$s \mapsto  \frac{\partial \alpha}{\partial t}(s,t) = \varphi_t \circ V(s)$, which is $D\varphi_t(V'(0))$.
  
To summarize the preceding discussion, there is a 
one-one correspondence between elements of $T_vTM$ and Jacobi fields along the geodesic $\gamma$ with $\dot\gamma(0) = v$.   Note that the pair $(J(t), J'(t))$ defines a section of $TTM$ over $\gamma(t)$.  We have the following key proposition:
\begin{proposition} \label{prop:key}
The image of the tangent vector $(v_1,v_2)\in T_vTM$ under the derivative of the
geodesic flow $D_v\varphi_t$ is the tangent vector $(J(t),J'(t))\in T_{\varphi_t(v)}TM$, where $J$ is the unique
Jacobi field along $\gamma$ satisfying $J(0) = v_1$ and $J'(0) = v_2$.
\end{proposition}

%The {\em Wronskian} of two Jacobi fields $J_1$ and $J_2$ is defined by
 %$$
% W(J_1,J_2) = \langle J_1, J_2'\rangle - \langle J_1', J_2 \rangle.
% $$
%It is constant since the Jacobi equation and one of the symmetries of the curvature tensor give
 %\begin{align*}
% W(J_1,J_2)' &= \langle J_1', J_2'\rangle + \langle J_1, J_2''\rangle -
  %              \langle J_1'', J_2 \rangle - \langle J_1', J_2' \rangle \\
    %  &= \langle J_1, R(J_2,\dot\gamma)\dot\gamma \rangle - 
    %                 \langle  R(J_1,\dot\gamma)\dot\gamma, J_2 \rangle \\
   %  &= \langle R(J_2,\dot\gamma)\dot\gamma, J_1 \rangle - 
    %                 \langle  R(J_1,\dot\gamma)\dot\gamma, J_2 \rangle \\
   %  &= \langle R(J_2,\dot\gamma)\dot\gamma, J_1 \rangle - 
   %                  \langle  R(J_2,\dot\gamma)\dot\gamma, J_1 \rangle \\
    %  &= 0.
% \end{align*}

Any vector field of the form $J(t) = (a + bt)\dot\gamma(t)$ is a Jacobi field, since in that case $R(J,\dot\gamma) = 0$ and the Jacobi equation reduces to $J'' = 0$, which holds since $\dot\gamma' = 0$. Conversely, any Jacobi field that is always tangent to $\gamma$ must have this form.  Computing the Wronskian of the Jacobi field $\dot\gamma$ and an arbitrary Jacobi field $J$ shows that $\langle J', \dot\gamma \rangle$ is constant.
It follows that if $J'(t_0) \perp \dot\gamma(t_0)$ for some $t_0$, then $J'(t) \perp \dot\gamma(t)$ for all $t$. Similarly if $J(t_0) \perp \dot\gamma(t_0)$ and $J'(t_0) \perp \dot\gamma(t_0)$ for some $t_0$, then $J(t) \perp \dot\gamma(t)$ and $J'(t) \perp \dot\gamma(t)$ for all $t$; in this  case we call $J$ a {\em perpendicular Jacobi field}.

An easy consequence of the above discussion is that any Jacobi field $J$ along a geodesic $\gamma$ can be expressed uniquely as $J = J_\parallel + J_\perp$, where $J_\parallel$ is a Jacobi field tangent to 
$\gamma$ and $J_\perp$ is a perpendicular Jacobi field.

\subsection{Matrix Jacobi and Riccati equations}\label{SSMJRE}

Choose an orthonormal basis $e_1=\dot\gamma(0),e_2, \ldots, e_n$ at $0$ for the tangent space at $\gamma(0)$ and parallel transport the basis  along $\gamma(t)$.  Let $\cR(t)$ be the matrix whose entries are
$$
\cR_{jk}(t)=\langle R(e_j(t),e_1(t))e_1(t),e_k(t)\rangle.
$$
Any Jacobi field can be written in terms of the basis as $J(t)=\sum_{k=1}^n y^k e_k(t)$ and the Jacobi equation can be written as 
$$
\frac{d^2 y^k}{dt^2}(t) + \sum_j y^j(t)\cR_{jk}(t ) =0.
$$
A solution is determined by values and derivatives at $0$ of the  $y^k$. 

 Let $\cJ(t)$ denote any matrix of solutions to the Jacobi equation. 
%Define the Wronskian of $\cJ$ to be the matrix $$W(\cJ,\cJ) =
%(\cJ^*)'\cJ - \cJ^*\cJ',$$ where $ *$  denotes transpose.
%Note that by the Jacobi equation, 
%\begin{equation}
%\label{eq:Wronskian}
%W'(\cJ,\cJ)=(\cJ^*)''\cJ-\cJ^*\cJ''=(\cR\cJ)^*\cJ-
%\cJ^*\cR\cJ=0,
%\end{equation}  
%since $\cR$ is symmetric.
When the matrix $\cJ$ is nonsingular, we can define 
$$
U=\cJ'\cJ^{-1}.
$$  
Then  $U$ satisfies the 
 {\em matrix Riccati equation}:
\begin{equation} \label{MRE}
U'  + U^2 + \cR   = 0,
\end{equation}
where $\cR$ is the matrix above. A standard calculation using the Wronskian shows that the operator  $U=\cJ'\cJ^{-1}$  is symmetric if and only if for any two columns $J_i,J_j$ of $\cJ$, we have $$\omega_{\RR^{2n}}((J_i,J_i'), (J_j,J_j'))=0,$$ where $\omega_{\RR^{2n}}$ is the standard symplectic form on $\RR^{n}$.

%\begin{lemma} 
%The operator  $U=\cJ'\cJ^{-1}$  is symmetric if and only if $$W(\cJ(t),\cJ(t))=0,$$ and  this is equivalent to the statement that for any two columns $J_i,J_j$ of $\cJ$, we have
%$$\omega((J_i,J_i'), (J_j,J_j'))=0,$$ where $\omega$ is the symplectic form on $T^1M$.
%\end{lemma}
%\begin{proof}
%The first statement follows from the definition of $W$ and the second follows from the definition of $\omega$.
%\end{proof}

\subsection{Perpendicular Jacobi fields and invariant subbundles}

There are two natural subbundles of $TTM$ that are invariant under the derivative $D\varphi_t$ of the geodesic flow, the first containing the second.  The first is the tangent bundle $TT^1M$ to the unit tangent bundle of $M$.  Under the natural identification $T_vTM \cong T_xM\times T_xM$, for $v\in T_x^1M$,
the subspace $T_vT^1M$ is the set of all pairs $(w_0,w_1)$ such that $\la v, w_1\ra = 0$.
To see this, note that if $\alpha(s,t)$ is a variation of geodesics 
generating the Jacobi field $J$ along the geodesic $\gamma$, with $\dot\gamma(0)=v$ and $\|{\partial\alpha}/{\partial t}(s,t)\| = 1$ for all $s,t$, then 
$$
0 = \left. \frac{D}{\partial s}\left\|\frac{\partial \alpha}{\partial t} \right\|^2\right\vert_{(0,0)} =  
\left. 2\left\la \frac{D^2}{\partial s\partial t} \alpha, \frac{\partial \alpha}{\partial t}\right\ra\right\vert_{(0,0)} =
\left. 2\left\la \frac{D^2}{\partial t\partial s} \alpha, \frac{\partial \alpha}{\partial t}\right\ra\right\vert_{(0,0)} = 2\la J'(0), \dot\gamma(0) \ra.
$$
The $D\varphi_t$-invariance of $TT^1M$  follows from the $\varphi_t$-invariance of $T^1M$. 
It is reflected in the fact, noted  at the end of Section~\ref{SSJfields}, that $\la J'(t), \dot \gamma\ra$ is constant for any Jacobi field $J$ along a geodesic 
$\gamma$.

The second natural invariant subbundle is the orthogonal complement $\dot\varphi^\perp$ in
$TT^1M$ to the vector field $\dot\varphi$ generating the geodesic flow. Under the natural identification $T_vTM \cong T_xM\times T_xM$, for $v\in T_x^1M$, the vector $\dot\varphi(v)$ is  $(v,0)$, and 
the subspace $\dot\varphi^\perp(v)$  is the set of all pairs $(w_0,w_1)$ such that $\la v, w_0\ra = \la v, w_1\ra = 0$. The  $D\varphi_t$-invariance
of $\dot\varphi^\perp$ follows from the observation, made at the end of section ~\ref{SSJfields}, that a Jacobi field $J$ with $J(t_0) \perp \dot\gamma(t_0)$ and $J'(t_0) \perp \dot\gamma(t_0)$  for some $t_0$ is perpendicular to $\gamma$ for all $t$.

To summarize, the space of all perpendicular Jacobi fields along $\gamma$ 
corresponds to the orthogonal complement to the direction of the geodesic flow $\dot\varphi(v)$ at the point $v = \dot\gamma(0) \in T^1M$.  To estimate the norm of the derivative $D\varphi_t$ on $TT^1M$, it suffices to restrict attention to vectors in the invariant subspace $\dot\varphi^\perp$; that is, it suffices to estimate the growth of perpendicular Jacobi fields along geodesics.

\subsection{Consequences of negative curvature and unstable Jacobi fields}

If the sectional curvatures of the Riemannian metric are negative along $\gamma$, then it follows from the Jacobi equation that $\langle J'', J \rangle > 0$, for any Jacobi field with the property that $J(t)$ and $\dot\gamma(t)$ are linearly independent. This has the following consequence; for a proof, see \cite{Eb2}.

\begin{lemma}
\label{lem:negcon}
 If the sectional curvatures are negative along $\gamma$, then
the functions $\|J(t)\|$ and $\|J(t)\|^2$ are strictly convex, for any nontrivial perpendicular Jacobi field $J$ along $\gamma$.
\end{lemma}

We also have the following results from \cite[Section 1.10]{Eb}. 
Let $\gamma: (-\infty, a] \to M$ be a geodesic ray along which the sectional curvatures of the Riemannian metric are always negative. Then,
for each $w \in \dot\gamma(a)^\perp$, there is a unique perpendicular Jacobi field $J_{+,w}$ along $\gamma$ such that $J_{+,w}(a) = w$ and 
$$
\|J_{+,w}(t)\| \leq \|w\| \quad\text{for all $t \leq a$.}
$$
Since $\|J_{+,w}(t)\|$ is a strictly convex function of $t$ by 
Lemma~\ref{lem:negcon}, $\|J_{+,w}(t)\|$ must be strictly increasing for $t \leq a$.
In fact $J_{+,w} = \lim_{\tau \to -\infty} J_{+,w,\tau}$, where $J_{+,w,\tau}$ is the Jacobi field such that $J_{+,w,\tau}(a) = v$ and $J_{+,w,\tau}(\tau) = 0$.
 We call $J_{+,w}$ an {\em unstable} Jacobi field.

For each $t \leq a$, there is a linear map $U_+(t): \dot\gamma(t)^\perp \to  \dot\gamma(t)^\perp$ such that
 $$
 J_+'(t) = U_+(t)(J_+(t))
 $$
 for every unstable Jacobi field $J_+$. A Jacobi field along $\gamma$ is unstable if and only if it satisfies $J' = U_+ J$.   
\begin{proposition} \label{prop:unstable}
The operators $U_+(t)$ are symmetric and positive definite. They satisfy the matrix Riccati equation (\ref{MRE}). Thus
$$
U_+'  + U_+^2 + \mathcal R   = 0 .
$$
In other words, for any vector $w \in \dot\gamma(t)^\perp$, we have: 
$$\langle w,U_+'(w)\rangle=-\langle R(w,\dot\gamma)\dot\gamma,w\rangle-\langle w,U_+ ^2(w)\rangle.
$$
\end{proposition}

 We call $U_+$ the {\em unstable solution} of the Riccati equation along the ray $\gamma$. If $v \in T^1M$ is a vector such that $\gamma_v(t)$ is defined for all $t < 0$, then we define $U_+(v)$ to be the operator $U_+(0)$ associated to the ray $\gamma_v: (-\infty,0] \to M$.

If $\gamma$ is a geodesic in a complete Riemannian manifold with negative curvature, the unstable Jacobi fields 
along $\gamma$ are obtained by varying $\gamma$ through geodesics $\beta$ such that $d(\beta(t), \gamma(t)) \leq d(\beta(0), \gamma(0))$ for  $t < 0$. These geodesics are orthogonal to a family of immersed hypersurfaces whose lifts to the universal cover of $M$ are called horospheres. The operators $U_+(t)$ are the second fundamental forms of horospheres.

There is an analogous definition of {\em stable} Jacobi fields and the {\em stable solution} of the Riccati equation along a ray $\gamma: [a,\infty)$.
If $\gamma : (-\infty, \infty) \to M$ is a complete geodesic, the unstable Jacobi fields along $\gamma$ are the stable Jacobi fields along the geodesic $t \mapsto \gamma(-t)$. We define $U_-(v)$ analogously to $U_+(v)$; it is symmetric and negative definite. The norm of a stable Jacobi field $J(t)$ defined on a ray $\gamma:[a,\infty) \to M$ is strictly decreasing for $t \geq a$.

Let
$$
 \text{ $\mathcal D = \{ v \in T^1M: \gamma_v(t)$ is defined for all 
           $t\}$. }
$$
If $v \in \mathcal D$, both $U_+(v)$ and $U_-(v)$ exist. This allows us to define a splitting of the $2n -1$ dimensional space $T_vT^1M$ as the direct sum of a one dimensional space $E^0(v)$ and two spaces $E^u(v)$ and $E^s(v)$ each of dimension $n-1$. The space $E^0(v)$ is $\RR\dot\varphi(v)$, and we will have $E^u(v) \oplus E^s(v) = \dot\varphi(v)^\perp$. In our usual coordinates, $E^0(v)$ is spanned by $(v,0)$ while
 $$
 E^u(v) = \{ (w,U_+(v)w) : w \in v^\perp \} \quad\text{and}\quad
 E^s(v) = \{ (w,U_-(v)w) : w \in v^\perp \}.
 $$
 The splitting at $v$ is mapped to the splitting at $\varphi_t(v)$ by $D\varphi_t$. 

The next proposition shows that while the splitting $T_{\mathcal D}T^1M = E^u\oplus E^0\oplus E^s$ is defined only over the set ${\mathcal D}$,  the geometry of this splitting is locally uniformly controlled.  

\begin{proposition} \label{prop:ctssplit} There exists a continuous function $\delta\colon T^1M\to \RR_{>0}$ such that for all $v\in {\mathcal D}$, if $(w,w') \in E^u(v)$, then
$$
\la w, w'  \ra \geq \delta(v) \|(w,w')\|_{Sas}^2,
$$
and if $(w,w') \in E^s(v)$, then
$$
\la w, w'  \ra \leq - \delta(v) \|(w,w')\|_{Sas}^2.
$$
%The splitting $E^0 \oplus E^u \oplus E^s$  extends continuously to a splitting over 
%$\overline{\mathcal D} \cap \hbox{\rm Int}(T^1M)$.

\end{proposition}

\begin{proof}  It will suffice to show that the functions
$$
\delta^u(v) = \inf_{(w,w')\in  E^u(v)\setminus \{0\}} \frac{\la w, w'\ra}{ \|(w,w')\|_{Sas}^2}\quad\hbox{ and }\,
\delta^s(v) = \inf_{(w,w')\in  E^s(v)\setminus \{0\}} -\frac{\la w, w'\ra}{ \|(w,w')\|_{Sas}^2}
$$
are locally uniformly  bounded away from $0$ for  $v\in {\mathcal D}$. We prove the statement for $\delta^s$.

Suppose that $\delta^s$ is not locally bounded away from $0$.  Then  there would be $v \in {\mathcal D}$, a sequence of vectors $v_n$ in $\mathcal D$ with $\lim_{n \to \infty} v_n =v $ and a sequence $\xi_n \in E^s(v_n)$ such that 
$\xi_n$ converges to a vector $\xi = (w,w')$  with $\la w, w'\ra = 0$.
By renormalizing we may assume that
$\|\xi_n\|_{Sas} = \|\xi \|_{Sas} = 1$ for each  $n$.

Since $v\in {\mathcal D}$, there exists $\delta > 0$ such that $\gamma_v(t)$ is defined for $|t| < \delta$. 
Let $J$ be the Jacobi field along the geodesic $\gamma_v$ determined by $\xi$, and let $J_n$
be the (stable) Jacobi field along $\gamma_{v_n}$ defined by $\xi_n$. 
Then $(\|J\|^2)'(0) = 2 \la w, w'\ra = 0$. On the other hand,
since $\xi_n \to \xi$ and $\|J_n(t)\|$ is a decreasing function of $t$ for each $n$, we see that $\|J\|$ is nonincreasing on $(-\delta,\delta)$. It follows from this and the
strict convexity of $\|J\|^2$ given by Lemma~\ref{lem:negcon} that the function $\|J\|^2$ cannot have a critical point in the interval $(-\delta,\delta)$.
\end{proof}

This proposition has the following corollary, which will be used for the Hopf argument in Section~\ref{s=generalsetup}.

\begin{corollary}\label{c=cones} Let  $\delta\colon T^1M\to \RR_{>0}$ be the function given by Proposition~\ref{prop:ctssplit}.  The continuous conefields 
$$
\cC^u(v) = \{(w,w')\in \dot\varphi^\perp(v) : \la w, w'  \ra \geq \delta(v) \|(w,w')\|_{Sas}^2\}
$$ 
and
$$
\cC^s(v) = \{(w,w')\in \dot\varphi^\perp(v) : \la w, w'  \ra \leq - \delta(v) \|(w,w')\|_{Sas}^2\}
$$
defined for $v\in  T^1M$ intersect only at the origin, and satisfy
$$
E^u(v) \subset \cC^u(v)\quad\hbox{and}\quad E^s(v) \subset \cC^s(v),
$$
 for all $v\in \cD$.
\end{corollary}

\section{A general criterion for ergodicity of the geodesic flow}\label{s=generalsetup}

In this section we establish a general criterion for ergodicity of the geodesic flow on a negatively curved manifold, not necessarily complete.  In the sections that follow we will verify that the hypotheses of our criterion hold for a quotient of Teichmuller space in the WP metric that is a finite branched cover of moduli space.

If $R$ is the curvature
tensor of a Riemannian metric on a manifold $M$, then for $x\in M$ we define 
$$
\| R_x\| = \sup_{v_1,v_2, v_3 \in T^1_xN} \| R_x(v_1, v_2)v_3 \|, \quad
\| \nabla R_x\| = \sup_{v_1,v_2,v_3,v_4 \in T^1_xN}
\|\nabla_{v_1}R_x(v_2, v_3)v_4 \|,
$$
and
$$
\| \nabla^2 R_x\| = \sup_{v_1,\ldots, v_5 \in T^1_xM}
\| \nabla^2_{v_1,v_2} R_x(v_3,v_4)v_5\| ,
$$
where $\nabla^2 R$ is the second covariant derivative of the curvature tensor:
$
\nabla^2_{X,Y} R = \nabla_X \nabla_Y R  -\nabla_{\nabla_XY}  R
$.

Let $M$ be a contractible Riemannian manifold, negatively curved, possibly incomplete.
Let $\Gamma$ be a group that acts freely and  properly discontinuously on $M$ by isometries,
and denote by $N$ the quotient manifold $N = M/\Gamma$.   We denote by $d$ both the path metric
on $M$ and the quotient metric on $N$, which is just the path metric for the induced Riemannian metric on $N$. The quotient map $p \colon M \to N$ is a covering map and a local isometry. 

Recall that the completion 
$\bar X$ of a metric space $(X,d)$ is the set of all Cauchy sequences $\langle x_n \rangle$ in $X$ modulo the equivalence relation:
$$\langle x_n \rangle \sim \langle y_n \rangle \quad\iff\quad \lim_{n\to\infty} d(x_n,y_n) = 0,
$$
with the induced metric
$
d(\langle x_n \rangle, \langle y_n \rangle) =  \lim_{n\to\infty} d(x_n,y_n)
$.
Let $\bar M$ be the metric completion of $M$ and let  $\bar N$
be the completion of $N$. Let $\partial N = \bar N\setminus N$. We will use $d$ to denote the metric on all of these spaces.

Consider the following additional six assumptions on $M$ and $N$:
\begin{itemize}
\item[I.] $M$  is  a {\em  geodesically convex:} for every $p , p'\in M$ there is a unique geodesic segment  in $M$ connecting $p$ to $p'$.  
\item[II.] {\em $\bar N$ is compact.} 
 \item[III.] $\partial N$ is {\em volumetrically cusplike:}  there exist constants $C>1$ and $\nu>0$ such that:
$$
\Vol \left( \{p\in N \,:\,  d(p,\partial N) < \rho\} \right) \leq C \rho^{2+\nu},
$$
 for every $\rho>0$.
\end{itemize}
For the final three assumptions we assume there exist constants $C>1$ and $\beta>0$ such that:
\begin{itemize}
 \item[IV.]  {\em $N$ has controlled curvature:} for all  $x\in N$, the curvature tensor  $R$  satisfies
$$\max\{\|R_x\| ,\|\nabla R_x\|, \|\nabla^2 R_x\| \} \leq  C  d(x,\partial N) ^{-\beta}.$$
 \item[V.]  {\em  $N$ has controlled injectivity radius:}  for every $x\in N$,
$$
\text{inj}(x)  \geq  C^{-1} d(x,\partial N) ^{\beta}.
$$
\item[VI.] {\em   The derivative of the geodesic flow is controlled:}  for every infinite geodesic $\gamma$ in $N$ and every $t\in [0,1]$:
$$\|D_{\dot\gamma(0)} \varphi_t \| \leq C d\left (\gamma\left([-t,t]\right),  \partial N\right )^{-\beta};
$$
\end{itemize}
Note that if II. and  III. hold, then $N$ has finite volume. In this case, we denote by $m$ the Riemannian volume (measure) on $N$, normalized so that $m(T^1N)=1$.  The main result in this section is:

\begin{theorem}\label{t=generalergodicity} Under assumptions I.-VI., the geodesic flow $\varphi_t$ on $T^1N$  is $m-$a.e. defined for all time $t$.  It is nonuniformly hyperbolic and ergodic (and in fact Bernoulli). The entropy $h(\varphi_t)$  of $\varphi_t$ is positive and finite, in fact equal to the  sum of the positive Lyapunov exponents of $\varphi_t$ with respect to $m$, counted with multiplicity.
\end{theorem}

\begin{remark} It seems that Assumption II. (compactness of $\bar N$) can be relaxed to the assumption that $N$ has finite diameter, but we have not verified all of the details.  We also remark that in applying Theorem~\ref{t=generalergodicity}, verifying Assumptions IV.-VI. is where the work lies.  In the case of the WP metric, assumptions I.-III. are either already known or follow in a straightforward way from known results.
\end{remark}

\begin{proof}[Proof of Theorem~\ref{t=generalergodicity}] We first establish 
several properties of $M$ that can be proved from assumptions I.-III. alone.  
The first such property is  $CAT(0)$.
A metric space $X$  is $CAT(0)$ if it is a geodesic space and and every geodesic triangle in $X$ satisfies the $CAT(0)$ inequality with the comparison Euclidean triangle (see \cite[p.159]{BridsonHaefliger}).

\begin{lemma} If I. holds, then $M$ and  $\bar M$ are both $CAT(0)$ spaces.
\end{lemma}

\begin{proof} 
The fact that $M$ is $CAT(0)$ follows from \cite[Theorem II.1A.6]{BridsonHaefliger} and Alexandrov's Patchwork  \cite[Proposition II.4.9]{BridsonHaefliger}.  The metric completion of a $CAT(0)$ space is $CAT(0)$, by \cite[Corollary II.3.11]{BridsonHaefliger}.
\end{proof}

\begin{proposition}[The flow is $\text{a.e.}$ defined for all time]\label{p=flowdefined} If I.--III. hold, then for almost every $v\in T^1M$, there exists an infinite geodesic (necessarily unique) tangent to $v$.   
\end{proposition}

Before proving this we state and prove another lemma that will be useful later as well.
Let $\pi:T^1N\to N$ be the natural projection. Let $$U_\rho =\{v\in T^1N \,:\, d(\pi(v),\partial N) < \rho\},$$
and let $S^+(\rho)$ be the set of all tangent vectors that flow into $U_\rho$ 
in some forward time $0\leq t \leq 1$.

\begin{lemma}\label{l=specialnbdvol} If I.--III. hold, then for $\rho<1$
$$m(S^+(\rho)) = O(\rho^{1+\nu}).$$
\end{lemma}

\begin{proof} Consider the ``shell" $S^+_k(\rho)$ of vectors $v$ that flow into  $U_\rho$ at times between $k\rho$ and $(k+1)\rho$. Any 
vector in this shell is in $U_{2\rho}$ at time  $(k+1)\rho$. Volume-preservation of the flow implies that the
the volume of $S_k^+(\rho)$ is at most the volume of $U_{2\rho}$, which is $O(\rho^{2+\nu})$, by assumption III.  The set $S^+(\rho)$ is contained in a union of the shells $S^+_0(\rho),\ldots, S^+_m(\rho)$, where $m$ 
is $O(\rho^{-1})$.
It follows that the volume of $S^+(\rho)$ is $O(\rho^{-1} \rho^{2+\nu}) = O(\rho^{1+\nu})$.
\end{proof}

\begin{proof}[Proof of Proposition~\ref{p=flowdefined}] 
The set of vectors such that the flow is not defined for some $0\leq t\leq 1$ is contained in $S^+(\rho)$ for all $\rho>0$.  By Lemma~\ref{l=specialnbdvol} this set has measure $0$.  It follows that the set of vectors for which the flow is defined for all time has full measure.  
\end{proof}

%The next proposition will be used later to prove ergodicity of the geodesic flow on $T^1N$,
%under further assumptions on the metric.

Suppose that $v\in TM$ determines an infinite geodesic ray $\gamma_v\colon [0,\infty) \to M$
tangent to $v$ at $0$. Since $M$ is a $CAT(0)$ space, the functions $b_{v,t}^s\colon M \to \RR$ defined by 
$$
b_{v,t}^s(y)  = d(y,\gamma_v(t)) - t
$$
converge uniformly on compact sets as $t\to\infty$ to a function $b_v^s\colon  M\to \RR$,
called a {\em (stable) Busemann function} \cite[Lemma II.8.18]{BridsonHaefliger}. For a fixed $v$, the Busemann function $b_v^s$ is clearly Lipschitz continuous, with Lipschitz norm $1$. 
If we assume that I. holds, then we can say more. 

 \begin{proposition}\label{p=nicehoro2}
Assume that  I. holds.  For any $v$ that determines an infinite geodesic ray $\gamma_v$, the function $b_v^s$ is convex and $C^1$,  and  $\|\grad b_v^s\| \equiv 1$.
 
For every $y\in M$, the unit vector
$$w_v^s(y) := -\grad  b_v^s(y)$$
defines an infinite geodesic ray $\gamma_{w_v^s(y)}\colon [0,\infty)\to \bar M$ tangent to
$w_v^s(y)$ at $0$ with the property that
$$
d(\gamma_v(t) , \gamma_{w_v^s(y)}(t)) \leq d(\gamma_v(0), y),
$$
for all $t\geq 0$.
\end{proposition}

\begin{proof}
Since $\gamma_v$ is an infinite ray, and $M$ is a geodesically convex Riemannian manifold, the functions $b_{v,t}^s$ are convex, $C^1$ and have the property that $\|\grad b_{v,t}^s(y)\| = 1$, for every $y\in M$.   Since $M$ is nonpositively curved, and $b_{v,t}^s$ converges uniformly on compact sets in $M$ to $b_v^s$, the desired properties of $C^1$ smoothness of $b_v^s$, convexity
and $\|\grad b_v^s\| \equiv 1$ follow from \cite[Lemma 3.4, and the following Remark]{BallmannGromovSchroeder}.  
 The final conclusion follows from \cite[Proposition II.8.2]{BridsonHaefliger}\label{p=nicehoro}.
\end{proof}

Suppose that $v\in T^1M$ determines an infinite  geodesic ray.
Proposition~\ref{p=nicehoro2} implies that  for each $t\in \RR$, the set $\cH^s_v(t) := (b_v^s)^{-1}(t)$ is a connected, codimension-$1$,  complete $C^1$ submanifold of $M$, called a {\em stable horosphere at level $t$}.  
For  such a $v$, we define:
$$\cW^s(v) = \{w_v^s(y)  \,:\, y\in \cH^s_v(0)\}.
$$
The set of basepoints  $\pi(\cW^s(v))$ in $M$ is the horosphere $\cH^s_v(0)$, and 
$\cW^s(v)$ is a continuous, codimension-$1$ submanifold of $T^1M$.
Similarly,  if $\gamma_v$ projects to a backward recurrent geodesic ray in $N$, we define the {\em unstable Busemann function} and {\em unstable manifold}:
$$
b_{v}^u(y)  =  \lim_{t\to\infty} d(y,\gamma_v(-t)) - t, \qquad\hbox{and}\quad \cW^u(v) = \{w^u_v(y)  \,:\, y\in \cH^u_v(0)\},
$$
 where $w^u_v(y) = -\grad b^u_v(y)$, and  $\cH^u_v(t) := (b_v^u)^{-1}(t)$ is the
{\em unstable horosphere} at level $t$ determined by $v$.

Our next proposition justifies the terminology ``stable and unstable manifolds" for $\cW^s(v)$
and $\cW^u(v)$.  The results stated up to this point all hold true when $M$ is nonpositively curved, but the proposition uses the negative curvature assumption on $M$ in an essential way.

We say that a geodesic ray $\gamma\colon [0,\infty)\to N$ is {\em (forward) recurrent} if the tangent vector $\dot\gamma(0)$ is an accumulation point for the tangent vectors
$\{\dot\gamma(t): t>0\}$.  We similarly define backward recurrence for a geodesic ray
 $\gamma\colon (-\infty,0]\to N$. An infinite geodesic is {\em recurrent} if it is both forward and backward recurrent.  Under assumptions I.-III., Proposition~\ref{p=flowdefined} and Poincar\'e recurrence imply that almost every $v\in T^1N$ determines
an infinite recurrent geodesic $\gamma_v\colon \RR\to N$ with $\dot\gamma_v(0) = v$.

\begin{proposition}[Contraction of horospheres]\label{p=largehoro} 
Assume I.--III. Let $v\in T_xM$ be tangent to an infinite geodesic ray $\gamma_v$ whose projection to $N$  is  forward recurrent. 
Let $y\in M$ be any other point, and let $w=w_v^s(y)\in T_yM$. Then $w$ is tangent to an
infinite geodesic ray $\gamma_w\colon [0,\infty) \to M$, and 
$$\lim_{t\to\infty}d(\gamma_v(t),\gamma_w(t+b_v^s(y)))=0;$$
moreover,
$$\lim_{t\to\infty}d_{Sas} (\varphi_t(v),\varphi_{t+b_v^s(y)}(w))=0.$$
In particular, if $\gamma_v$ projects to a forward recurrent geodesic ray in $N$, then for every 
$t>0$, 
$\varphi_t(\cW^s(v)) = \cW^s(\varphi_t(v))$, and
for every $w\in \cW^s(v)$, we have
$\lim_{t\to\infty}d_{Sas} (\varphi_t(v), \varphi_t(w))=0$.

Similarly, if $v$ is tangent to a backward ray $\gamma_v\colon (-\infty,0]\to M$ whose projection is recurrent, then 
$w=w^u_v(w)$ is tangent to a backward ray $\gamma_w\colon (-\infty,0]\to M$, and
$$\lim_{t\to -\infty}d_{Sas} (\varphi_t(v),\varphi_{t+b_v^s(y)}(w))=0.$$
 In particular, for every $w\in \cW^u(v)$, we have $\lim_{t\to-\infty} d(\varphi_t(v),\varphi_t(w))=0$.
\end{proposition}
Before  beginning the proof we remark that  in  \cite{DaskWent} a property called {\em nonrefraction} was proved for the WP metric. Using that result, a short proof of the above proposition was given in the WP case in \cite{BMM}. 

\begin{proof} Let $\gamma_v\colon [0,\infty)\to M$ be an infinite geodesic ray whose projection to $N$ is recurrent, and let $x=\gamma_v(0)$ be the footpoint of $v$. Suppose that
$x'\in M$ is another point, and let $v' = w_v^s(x')$. 
Since $\bar M$ is $CAT(0)$, the distance $d(\gamma_{v}(t), \gamma_{v'}(t))$  is a convex function of $t$; since it is bounded, it must be nonincreasing, and hence bounded above for all $t$ by $d(x,x')$.  We claim that if $d(x,x') < d(x,\partial M)$, then the image of $\gamma_{v'}$
must lie entirely in $M$. Since the projection of $\gamma_v$ to $N$  is recurrent, there exist
sequences $g_n\in\Gamma$ and $t_n\to\infty$ such that 
%\begin{eqnarray*}
%d(x, g_n \gamma_v(t_n)) & <& d(x,\partial M) - d(x,x')\\
%&\leq & d(x,\partial M) - d(\gamma_{v}(t_n), \gamma_{v'}(t_n)),\\
%&= & d(x,\partial M) - d(g_n\gamma_{v}(t_n), g_n\gamma_{v'}(t_n)),
%\end{eqnarray*}
$$d(x, g_n \gamma_v(t_n))  < d(x,\partial M) - d(x,x').$$
Then
$$d(x, g_n \gamma_v(t_n)) <  d(x,\partial M) - d(\gamma_{v}(t_n), \gamma_{v'}(t_n)) =  d(x,\partial M) - d(g_n\gamma_{v}(t_n), g_n\gamma_{v'}(t_n)),$$
which implies, by the triangle inequality,  that $d(x, g_n\gamma_{v'}(t_n)) < d(x,\partial M)$.
Hence $g_n\gamma_{v'}(t_n)\in M$, and so  $\gamma_{v'}(t_n)\in M$; geodesic convexity
of $M$ implies that $\gamma_{v'}[0,t_n]\subset M$, for all $n$, which proves the claim.

Now a standard ruled surface argument using geodesic convexity and the negative curvature of $M$  (see e.g. \cite[Theorem 4.1]{BMM}, where it is proved in the WP context) 
shows that for every $\gamma_v$ that projects to a recurrent geodesic ray in $N$,
and any $y\in M$ with the property that $\gamma_{w_v^s(y)}[0,\infty)\subset M$, 
the distance
$d(\gamma_{w_v^s(y)}(t), \gamma_v[0,\infty))$ is strictly decreasing in $t$ and tends to $0$ as $t\to\infty$. (Alternately, one can show this using Jacobi fields).  What is more,
this convergence takes place in the tangent bundle:
$$
\lim_{t\to\infty} d_{Sas}(\dot\gamma_{w_v^s(y)}(t) , \dot\gamma_v[0,\infty) ) =0.
$$

Now suppose that $y\in M$ is an arbitrary point. Connect $y$ to $x=\gamma_v(0)$ by a geodesic arc $\sigma$ in $M$. Fix $\epsilon_0>0$ such that $d(x,\partial M) <\epsilon_0$. We claim that if $x'$ is any point on $\sigma$ that satisfies
$$\lim_{t\to\infty} d(\gamma_{w_v^s(x')}(t), \gamma_v[0,\infty)) = 0,$$
then for any point $y'$ on $\sigma$ such that $d(x',y')<\epsilon_0/3$:
$$\lim_{t\to\infty} d(\gamma_{w_v^s(y')}(t), \gamma_v[0,\infty)) = 0.$$
From the claim it follows that $\lim_{t\to\infty} d_{Sas}(\dot\gamma_{w_v^s(y)}(t), \dot\gamma_v[0,\infty)) = 0$.

To prove the claim, suppose that $x'$ and $y'$ are given. Since the distance
$d((\gamma_{w_v^s(x')}(t), \gamma_{w_v^s(y')}(t))$ is bounded for all $t>0$ and
convex, it is nonincreasing, and hence bounded above by $\epsilon_0/3$, for
all $t>0$.  If $T>0$ is sufficiently large,
then the distance from $\gamma_{w_v(x')}(t)$ to $\gamma_v$ is less than $\epsilon_0/3$ for
all $t>T$. Since $\gamma_v$ projects to a recurrent ray in $N$, there exist $g_n\in \Gamma$ and $t_n\to\infty$ such that $d(\gamma_v(t_n), g_nx)<\epsilon_0/3$.   It follows that
$\gamma_{w_v^s(y')}(t_n)\in M$ when $t_n > T$, which implies that $\gamma_{w_v(y')}[0,\infty)\subset M$.  The claim follows.  

A simple application of the triangle inequality shows that the property 
$\lim_{t\to\infty} d (\gamma_{w_v^s(y)}(t), \gamma_v[0,\infty))= 0$ implies
that $$\lim_{t\to\infty}d(\gamma_v(t),\gamma_{w_v^s(y)}(t+b_v^s(y)))=0.$$
Since $\lim_{t\to\infty} d_{Sas}(\dot\gamma_{w_v^s(y)}(t), \dot\gamma_v[0,\infty))= 0$
for every $y\in M$, we conclude that 
$$\lim_{t\to\infty}d_{Sas}(\varphi_t(v),\varphi_{t+ b_v^s(y)}(w_v^s(y))=0.
$$
\end{proof}

The proof of Theorem~\ref{t=generalergodicity} now proceeds in several steps. The first is to
establish nonuniform hyperbolicity.  This is a classical result for closed manifolds with negative curvature; see, e.g.,  \cite[Section 17.6]{KatokHasselblatt}.

We need the following lemma.
\begin{lemma}\label{l=integrability} Assume that hypotheses I.--III. hold.
Let $\varphi_1$ be the time-$1$ map of the geodesic flow. Then 
$$
\int_{T^1N} \log^+\|D\varphi_1\|\, dm <\infty 
\quad\text{and}\quad 
\int_{T^1N} \log^-\|D\varphi_1\|\, dm <\infty.$$
\end{lemma}

\begin{proof}  
Lemma~\ref{l=specialnbdvol} implies that for $n\geq 1$, $m(S^+(1/n))=O((1/n)^{1+\nu})$.
On $S^+(1/n)$ we have $\log^+\|D\varphi_1\| =O(\log n)$, and hence 
 $$\int_{S^+(1/n)} \log^+\|D\varphi_1\|\,dm  =O(\log n/n^{1+\nu}).$$ 
Summing over $n$ gives the first half of the conclusion. The second half follows from the first and equivariance of the geodesic flow under the $m$-preserving involution $u \mapsto -u$: if $w = \varphi_1(v)$, then $-v = \varphi_1(-w)$.
\end{proof}

It follows from the lemma that $\log\|D\varphi_1\|$ is integrable. Consequently Oseledec's theorem can be applied to the cocycle $D\varphi_1$. It implies that  
for $m$-almost every $v \in T^1N$ there exist $k(v)\leq 2n-1$ real numbers $$\lambda_1(v)< \lambda_2(v) < \cdots < \lambda_{k(v)}(v)$$ and a $D\varphi_t$-invariant  splitting $T_vT^1N = \bigoplus_{i=1}^{k(v)} E_i(v)$ such that for every nonzero vector
$\xi\in E_i(v)$:
$$
\lim_{t\to\pm\infty} \frac{1}{t} \log\|D_v\varphi_t(\xi)\| = \lambda_i(v).
$$
The functions $k(v), \lambda_i(v)$, and $E_i(v)$ depend measurably on $v$.
The numbers $\lambda_i(v)$ are called the {\em Lyapunov exponents of $\varphi_t$ at $v$},
and $E_i(v)$ the {\em Lyapunov subspaces}. Since the orthocomplement $\dot\varphi^\perp$ is $D\varphi_t$-invariant, and the restriction of $D\varphi_t$ preserves a natural symplectic form, the Lyapunov exponents of 
$\varphi_t$ are paired: if $\lambda$ is a Lyapunov exponent, then so is $-\lambda$. Moreover,
 since the generating vector field $\dot\varphi$ is preserved by $D\varphi_t$, it follows that 
$$\lim_{t\to\pm\infty} \frac{1}{t} \log\|D_v\varphi_t(\xi)\| = 0,$$
for any $\xi$ tangent to the orbits.

For $v\in T^1N$ such that the geodesic $\gamma_v(t)$ is defined for all $t$, let  $E^u(v)$ be the subspace of $T_vT^1N$ spanned by the unstable perpendicular Jacobi fields at $v$, and $E^s(v)$ the subspace spanned by the stable perpendicular Jacobi fields at $v$. These spaces each have dimension $n-1$ and 
 $$
 T_vT^1N = E^s(v)\oplus E^0(v)\oplus E^u(v),
 $$
 where $E^0(v)$ is the one dimensional subbundle tangent to the orbits of the flow $\varphi_t$. The splitting at $v$ is mapped to the splitting at $\varphi_t(v)$ by $D\varphi_t$. 

%Moreover the $E^u(v)$ and $E^s(v)$
%depend continuously on $v$ by Proposition~\ref{prop:unstableiscts}.

\begin{lemma} There is a $\varphi_t$-invariant set $\Lambda_0 \subset T^1N$ of full measure with respect to $m$ such that for every $v \in \Lambda_0$ we have 
$$
 E^s(v) = \bigoplus_{\lambda_i(v) < 0} E_i(v) 
\quad\text{and}\quad
 E^u(v) = \bigoplus_{\lambda_i(v) > 0} E_i(v).
$$
\end{lemma}

\begin{proof} We choose $\Lambda_0$ to be the set of vectors $v \in T^1N$ such that
\begin{enumerate}
\item \label{orbdef}
$\varphi_t(v)$ is defined for all $t$;
\item \label{expdef}
the exponents $\lambda_i(v)$ are defined for $i =1,\dots,k(v)$; and
\item \label{urec}
 $v$ is uniformly forward and backward recurrent under the flow $\varphi_t$.
\end{enumerate} 
The last property means the following:
\begin{enumerate} 
\item[(\ref{urec}$'$)] \label{urec'} 
for any neighborhood $U$ of $v$, there is $\delta > 0$ such that for all large enough $T$ the sets $R_+(T) = \{t \in [0,T]: \varphi_t(v) \in U\}$ and $R_-(T) = \{t \in [0,T]: \varphi_{-t}(v) \in U\}$ both have Lebesgue measure at least $\delta T$. This ensures that both sets contain finite subsets of cardinality at least $\delta T - 1$ in which distinct elements differ by at least $1$.
\end{enumerate}

Properties (\ref{orbdef})--(\ref{urec}) hold for $m$-almost all vectors in $v \in T^1N$.  For (\ref{orbdef}) this is Proposition~\ref{p=flowdefined}, for (\ref{expdef}) it is a part of Oseledec's theorem, and for (\ref{urec}) it follows from a standard argument using the Birkhoff ergodic theorem.

Since the set $\Lambda_0$  is invariant under the involution $u \mapsto -u$ and the derivative of this involution maps $E^s(u)$ to $E^u(-u)$, it will suffice to prove the second statement. To this end, recall that if $J$ is 
a nonzero unstable Jacobi field along a geodesic $\gamma$, then $\|J(t)\|$ is a strictly increasing convex function. Given $v \in \Lambda_0$, we can choose a neighborhood $U$ of $v$ and $\eta > 0$ such that if $J(t)$ is an unstable Jacobi field along a geodesic $\gamma$ with $\dot\gamma(0) \in U$, then $\|J(1)\| \geq (1+ \eta)\|J(0)\|$. With $\delta$ chosen as in (\ref{urec}$'$), we obtain 
$$
\|J(T)\| \geq (1+ \eta)^{\delta T - 1}\|J(0)\|,
$$
for any unstable Jacobi field $J(t)$ along the geodesic $\gamma_v(t)$.
\end{proof}

We summarize the consequences of the discussion since Lemma~\ref{l=integrability} in the following:
\begin{proposition}[Nonuniform hyperbolicity]\label{p=nuh} Under assumptions I.-VI., the geodesic flow is nonuniformly hyperbolic.  On the   full measure, $\varphi_t$-invariant subset $\Lambda_0\subset T^1N$ defined above there is a measurable $D\varphi_t$-invariant  splitting of the tangent bundle:
$$
T_{\Lambda_0}(T^1N) = E^s\oplus E^0 \oplus E^u
$$
such that, for every $v\in \Lambda_0$:
\begin{enumerate}
\item $E^0(v)$ is tangent to the orbits of the flow:  $E^0(v) = \RR\dot \varphi(v)$;
\item $E^u(v)$ is spanned by the unstable perpendicular Jacobi fields at $v$, and 
$E^s(v)$ is spanned by the stable perpendicular Jacobi fields at $v$; and
\item for every nonzero $\xi^u \in E^u(v)$, $\xi^s \in E^u(v)$:
$$
\lim_{t\to \infty}  \frac{1}{t} \log \|D_v\varphi_t(\xi^u)\|  >0,\quad\text{ and }\,
\lim_{t\to \infty}  \frac{1}{t} \log \|D_v\varphi_t(\xi^s)\|  < 0,
$$
and the limits are finite.
\end{enumerate}
\end{proposition}

This completes the first step. The next is to introduce the local stable and unstable manifolds.

\begin{proposition}[Existence and absolute continuity of families of local stable manifolds]\label{p=aclam} 
Assume I.-VI.  Let $n = \dim(N)$, and let $\Lambda_0\subset T^1N$ be given by Proposition~\ref{p=nuh}. 
There exist a full volume, $\varphi_t$-invariant subset $\Lambda_1\subset \Lambda_0$, a measurable function $r\colon \Lambda_1\to \RR_{>0}$
and measurable families of $C^\infty$, $(n-1)$-dimensional embedded disks $\cW_\loc^s = \{\cW_\loc^s(v): v\in \Lambda_1\}$ and $\cW_\loc^u = \{\cW_\loc^u(v) : v\in \Lambda_1\}$ with the following properties.  For each $v\in \Lambda_1$:
\begin{enumerate}
\item[(1)] $\cW^s_\loc(v)$ is tangent to $E^s(v)$ and $\cW_\loc^u(v)$ is tangent to 
$E^u(v)$ at $v$;
\item[(2)] for all $t>0$, 
$$\varphi_t(\cW^s_\loc(v))\subset \cW^s_\loc(\varphi_t(v)),\quad \text
{ and } \,\varphi_{-t}(\cW^u_\loc(v))\subset \cW^u_\loc(\varphi_{-t}(v));
$$
\item[(3)] $w\in \cW_\loc^s(v)$ if and only if $d(v,w)< r(v)$ and
$$\lim_{t\to\infty} d_{Sas}(\varphi_t(v),\varphi_t(w)) = 0;$$
\item[(4)]$w\in \cW^u_\loc(v)$ if and only if $d(v,w)< r(v)$ and
$$\lim_{t\to-\infty} d_{Sas}(\varphi_t(v),\varphi_t(w)) = 0.$$
\end{enumerate}
Moreover, for $\ast\in\{s,u\}$, the family $\cW_\loc^\ast$ is absolutely continuous.  In particular:
\begin{enumerate}
\item[(5)] if $Z\subset T^1N$ has volume $m(Z)=0$, then for $m$-almost every $v\in \Lambda_1$, 
the set  $Z\cap \cW_\loc^\ast(v)$ is a zero set  in $\cW_\loc^\ast(v)$ (with respect to the induced $(n-1)$-dimensional Riemannian volume); and
\item[(6)] if $D\subset T^1N$ is any $C^1$-embedded, $n$-dimensional open disk, and $B\subset D$ has induced Riemannian volume zero in $D$, then
 $m(\Sat_\loc^\ast(B)) = 0$, where
$$
\Sat_\loc^\ast(B) := \bigcup_{\{v\in \Lambda_1 \,: \, \cW_\loc^\ast(v)\cap B\neq \emptyset\}} \cW_\loc^\ast(v).
$$
\end{enumerate}
\end{proposition}

The conclusions of Proposition~\ref{p=aclam} will follow from the main results in \cite{KS}.  To apply these results, it is necessary to verify a list of hypotheses, some of a technical nature, concerning the $C^3$ properties of the Sasaki metric and the geodesic flow.
We defer the verification of these properties, assuming I.-VI.,  to Appendix B and now show how 
Proposition~\ref{p=aclam} can be used to prove ergodicity of $\varphi_t$.  Properties (5) and (6) in Proposition~\ref{p=aclam} are the heart of the matter in proving ergodicity.  Property (5) is a form of ``leafwise absolute continuity" and (6) is a form of ``transverse absolute continuity."  

Properties (5) and (6) are obvious if $\cW_\loc^s(v)$ and $\cW_\loc^u(v)$ depend smoothly on $v$, as they do for the geodesic flow of a manifold of constant negative curvature. But this is rarely the case. Examples of compact manifolds for which  the bundles $E^s$ and $E^u$ are only H\"older continuous have been given by  Anosov \cite{An67} and Hasselblatt \cite{Hassel}, and their techniques extend to the present context. However these examples do not appear to rule out the curious and extremely unlikely possibility that the bundles are smooth for the special case of the WP metric.

Let $\Omega_1$ be the full measure set of $v\in T^1M$ such that $\gamma_v$ projects to a (forward and backward) recurrent 
geodesic in $TN$. Each $v\in \Omega_1$ has a stable manifold $\cW^s(v)$ and an unstable manifold $\cW^u(v)$. 
For  $\delta < \inj(\pi(v'))$, where    $v'=Dp(v)\in Dp(\Omega_1)$,   denote by 
$\cW^\ast(v,\delta)$ the connected component of $\cW^\ast(v)\cap B_{T^1M}(v,\delta)$ containing $v$, where $B_{T^1M}(v,\delta)$ is the Sasaki ball  of radius $\delta$  in $T^1M$
centered at $v$. We denote  by $\cW^s(v',\delta)$ the projection $Dp(\cW^s(v,\delta))$; it is an $(n-1)$-dimensional embedded disk.

Notice that, for every $v\in \Omega_1$, if  $v' = Dp(v)$ belongs to the full measure
set $\Lambda_1$ of Proposition~\ref{p=aclam}, then  
the local stable manifold $\cW^s_{loc}(v')$ through $v'$  must coincide with
$\cW^s(v', r(v'))$, where $r\colon \Lambda_1\to \RR_{>0}$ is the function
given by Proposition~\ref{p=aclam}.  

  At this point, we have established the almost everywhere existence of the global, complete submanifolds $\cW^s(v)$ and $\cW^u(v)$  in $T^1M$, invariant under the flow, but we have not shown them to have any absolute continuity properties.  On the other hand, the
local Pesin stable and unstable manifolds $\cW^s_{loc}(v)$ and   $\cW^u_{loc}(v)$ 
have good absolute continuity properties, but they are not complete submanifolds -- they are open disks with measurably varying radii.  To prove ergodicity, we would like a collection of complete subbmanifolds forming an absolutely continuous (almost everywhere) foliation with controlled geometry.  The key step in showing this is to use this  almost everywhere coincidence of the global submanifolds with the local Pesin disks to obtain absolute continuity of the global foliation.  This is the content of the next proposition.

\begin{proposition}[Smoothness and absolute continuity of horospherical laminations]\label{p=achoro}
Assume I.-VI. There is a full volume subset $\Omega_2\subset \Omega_1$ 
such that for $\ast\in\{s,u\}$ and
for $v\in \Omega_2$,  the Busemann function $b_v^\ast\colon M\to \RR$ is $C^\infty$. 
The leaves of the lamination $\cW^\ast = \{\cW^\ast(v) : v\in \Omega_2\}$ 
are $C^\infty$ submanifolds of $T^1M$ diffeomorphic to $\RR^{n-1}$.

Let $\Lambda_2 = Dp(\Omega_2)$.  The family of manifolds 
$$\{\cW^\ast(v, \delta) : v\in \Lambda_2, \,\delta<\inj(\pi(v))\}
$$
has the following absolute continuity properties.
\begin{enumerate}
\item if $Z\subset T^1N$ has volume $m(Z)=0$, then for $m$-almost every $v\in \Lambda_2$, 
and every $\delta < \inj(\pi(v))$,
the set  $Z\cap \cW^\ast(v,\delta)$ is a zero set  in $\cW^\ast(v,\delta)$ (with respect to the induced $(n-1)$-dimensional Riemannian volume); and
\item if $D\subset T^1N$ is any smoothly embedded, $n$-dimensional open disk, and $B\subset D$ has induced Riemannian volume zero in $D$, then  for any $\delta< \frac12\inf_{v\in D}\inj(\pi(v))$, we have $m(\Sat^\ast(B, \delta)) = 0$, where
$$
\Sat^\ast(B,\delta) := \bigcup_{\{v\in \Lambda_2 \,: \, \cW^\ast(v,\delta)\cap B\neq \emptyset\}} \cW^\ast(v,\delta).
$$
\end{enumerate}
\end{proposition}

\begin{proof} 
We first show that $\cW^s(v)$ is a $C^\infty$ submanifold of $T^1M$, for
almost every $v\in T^1M$. For any $\epsilon>0$ there exists a compact set $\Delta_\epsilon \subset \Lambda_1$ of measure $m(\Delta_\epsilon) > 1-\epsilon$ such that the restriction of the function $r$ from Proposition~\ref{p=aclam} to $\Delta_\epsilon$ is continuous and bounded from below by a constant $r_\epsilon>0$.  Fix $\epsilon>0$, and let $\Delta_\epsilon^s \subset \Delta_\epsilon$ be the set of vectors $v'\in \Delta_\epsilon$ such that $\varphi_{k_n}(v')\in \Delta_\epsilon$ for a sequence of integers $k_n\to\infty$.  Poincar{\'e} recurrence implies that $m(\Delta_\epsilon\setminus \Delta_\epsilon^s) = 0$.

 Fix $v'\in \Delta_\epsilon^s\cap Dp(\Omega_1)$. Let  $v\in Dp^{-1}(v')$ be an  arbitrary lift of $v'$ to $T^1M$, and let $w\in \cW^s(v)$.  We show that $\cW^s(v)$ is $C^\infty$ in a neighborhood of $w$; as $w$
is arbitrary, this implies that $\cW^s(v)$ is $C^\infty$.  
Since $v' = Dp(v) \in \Delta_\epsilon^s$, there exists a sequence $k_n\to \infty$ such that  $\varphi_{k_n}(v')\in \Delta_\epsilon$.  At the same time, Proposition~\ref{p=largehoro} implies
that
$$
\lim_{t\to\infty} d_{Sas}(\varphi_t(v), \varphi_t(w)) = 0, 
$$
and so for $n$ sufficiently large, $d_{Sas}(\varphi_{k_n}(v),\varphi_{k_n}(w)) < r_\epsilon/2$,
where $r_\epsilon>0$ is the lower bound on the restriction of $r$ to $\Delta_\epsilon$.
But this implies that  $Dp( \varphi_{k_n}(w)) \in \cW^s_{\loc}(\varphi_{k_n}(v'))$.  Since $\varphi_{k_n}$
is a diffeomorphism,  we conclude that there is a neighborhood of $w$ in $\cW^s(v)$ that is
diffeomorphic to the $C^\infty$ submanifold $\cW^s_{\loc}(\varphi_{k_n}(v'))$.  Since $w$ was arbitrary,
this implies that $\cW^s(v)$ is a $C^\infty$ submanifold of $T^1M$.  The intersection 
$\Lambda_2^s := \bigcap_{\epsilon>0} \Delta_\epsilon^s \cap Dp(\Omega_1)$ is a full volume 
subset of $T^1N$, and we have shown
that for every $v \in \Omega_2^s :=  Dp^{-1}(\Lambda_2^s)$, the
submanifold $\cW^s(v)$ is $C^\infty$.

For each $v\in \Omega_2^s$, consider 
the map $\psi$ from $\cH^s_v\times \RR$ to $M$ that sends
$(y,t)$ to $\pi(\varphi_t(w^s_v(y)))$, where  $w^s_v(y) = -\grad b^s_v(y)$. 
Since $\cW^s(v)$ is $C^\infty$, the function $w^s_v(y)$ is $C^\infty$ along $\cH^s_v$;
it follows that $\psi$ is a diffeomorphism.  In the coordinates on $M$ given by 
$\psi$, the Busemann function $b^s_v$  assigns the value $-t$ to the point $(x,t)$.  It follows that $b^s_v$ is $C^\infty$, for every $v\in\Omega_2^s$.
 Similarly, there is a set $\Omega_2^u$ of full measure
such that $b^u_v$ is $C^\infty$ for every $v\in \Omega_2^u$.  Setting $\Omega_2=\Omega_2^u \cap \Omega_2^s$, we obtain
the full measure set where the conclusions of the proposition will hold. 

We establish the absolute continuity properties of $\cW^s$; analogous arguments show the properties for $\cW^u$.   The preceding arguments show that for every $v\in \Lambda_2$ there exists an integer $k\geq 0$ such that
\begin{eqnarray}\label{e=kflowinclude}
\varphi_k(\cW^s(v,\delta)) \subset \cW^s_\loc(\varphi_k(v)),\,\text{ for every }
\delta<\inj(\pi(v))
\end{eqnarray}
For a fixed $k\geq 0$, denote by $X_k$ the set of $v\in \Lambda_2$ for which (\ref{e=kflowinclude}) holds. Then $\Lambda_2 = \bigcup_{k\geq 0} X_k$.

 Suppose that  $m(Z)=0$, for some $Z\subset T^1N$.  Then the set $\hat Z = \bigcup_{k\geq 0} \varphi_k(Z)$ also has measure $0$.  It follows from Proposition~\ref{p=aclam} that for almost every $w\in \Lambda_1$,  the induced Riemannian measure of $\hat Z$ in $\cW^s_\loc(w)$ is zero.   But this implies in particular that for every $k\geq 0$ and for almost every $v\in X_k$,  the induced Riemannian measure of $\varphi_k(Z)\subset \hat Z$ in $\varphi_k(\cW^s(v,\delta)) \subset \cW^s_\loc(\varphi_k(v))$  is zero; hence the induced volume of $Z$ in $\cW^s(v,\delta)$ is $0$, for all $\delta<\inj(\pi(v))$. This establishes (1).

Suppose that $D$ is a $C^1$-embedded, $n$-dimensional disk in $T^1N$.  Fix $\delta< \frac12\inf_{v\in D} \inj(\pi(v))$.   Suppose that $B\subset D$ has induced Riemannian volume $0$.  Let 
$$B_k = B \cap \bigcup_{w\in X_k} \cW^s(w,\delta)$$
and note that
$$
\Sat^s(B,\delta) = \bigcup_{k\geq 0} \Sat^s(B_k,\delta);
$$
hence it suffices to show that $m(\Sat^s(B_k,\delta)) = 0$, for all $k\geq 0$.

Fix $k\geq 0$. For each $w\in B_k$, there an $n$-dimensional open ball $D_w\subset D$ centered at $w$ in the induced Riemannian metric in $D$,  such that 
$\bigcup_{j=0}^k\varphi_{j}(D_w)\subset T^1N$.  
Since $\varphi_{k}$ is a diffeomorphism, the set $\varphi_{k}(B_k\cap D_w)$ has induced
Riemannian volume zero in the $n$-dimensional disk $\varphi_{k}(D_w)$.  It follows from 
Proposition~\ref{p=aclam} that $m(\Sat_\loc^s(\varphi_{k}(B_k\cap D_w)))=0$, and so
$$
m\left(\varphi_{-k}\left(\Sat_\loc^s(\varphi_{k}(B_k\cap D_w))\right)\right)=0.
$$
But  (\ref{e=kflowinclude}) implies that 
$$
\Sat^s(B_k\cap D_w,\delta )   \subset 
\varphi_{-k}\left(\Sat_\loc^s(\varphi_{k}(B_k\cap D_w))\right),
$$
and so $m(\Sat^s(B_k\cap D_w,\delta)) = 0$.  Now fix a countable cover $\{D_{w_i}\,:\, w_i\in B_k\}$ of $B_k$ in $D$ by such balls (this is possible by the Besicovitch covering theorem, since $D$ is an embedded $C^1$ submanifold). Then 
$$
\Sat^s(B_k,\delta) \subset \bigcup_{i} \Sat^s(B_k\cap D_{w_i},\delta ),
$$
and so $m(\Sat^s(B_k,\delta) )= 0$.  Conclusion (2) follows.
\end{proof}

We remark that  Proposition~\ref{p=largehoro} and Proposition~\ref{p=achoro} show that the horospheres  $\cH^*_v(0)$ 
are the level sets of  regular  values of $C^\infty$ functions. Consequently they are complete submanifolds of $T^1M$. As remarked above, the smooth manifolds given by Propositon~\ref{p=achoro} may be open and hence have boundary. 

\begin{proof}[Proof of ergodicity.]   Assume I.-VI. The proof that $\varphi_t$ is ergodic is an adaptation of the standard ``Hopf Argument," along the lines of the proof of local ergodicity in \cite{KatokBurns}.  To prove ergodicity, it suffices to show that for every continuous function
$f\colon T^1N \to \RR$ with compact support:
\begin{eqnarray}\label{e=ergav}
\lim_{T\to \infty} \frac{1}{T}\int_{0}^T f(\varphi_t(v))\, dt = \int_{T^1N} f\, dm,\quad\text{for }\, m-a.e.\,\, v\in T^1N
\end{eqnarray}
Indeed, if (\ref{e=ergav}) holds for a dense set of functions $f$ in $L^2$, then by continuity of the  projection $f\mapsto B(f) = \lim_{T\to \infty} \frac{1}{T}\int_{0}^T f\circ \varphi_t \, dt$,
(\ref{e=ergav}) will hold for every $f$ in $L^2$.

Fix then a continuous function $f$ with compact support and define measurable functions
$f^s$ and $f^u$ by:
$$
f^s(v) = \limsup_{T\to \infty} \frac{1}{T}\int_{0}^T f(\varphi_t(v))\, dt, \,\text{ and }\, 
f^u(v) = \limsup_{T\to \infty} \frac{1}{T}\int_{-T}^0 f(\varphi_t(v))\, dt. 
$$
The Birkhoff Ergodic Theorem implies that there is a set $G\subset T^1N$ of full measure
such that for every $v\in G$, we have $f^s(v) = f^u(v) = B(f)(v)$.    
Since $f$ is continuous with compact support, and the leaves of $\cW^s$ are contracted by $\varphi_t$, it follows that $f^s$ is constant along leaves of $\cW^s$.  Similarly, $f^u$ is constant along leaves of $\cW^u$. Finally, all three functions $f^s, f^u, B(f)$ are invariant under the flow $\varphi_t$.

Now fix a arbitrary element $v\in T^1N$.  
We will show that there is a neighborhood $U_v$ of $v$
on which $B(f)$ is almost everywhere constant.  Since $T^1N$ is connected, this will imply that
$B(f)$ is almost everywhere constant on $T^1N$.  Since $\int_{T^1N} B(f)\, dm = \int_{T^1N} f\, dm $, it will then follow that (\ref{e=ergav}) holds, and so $\varphi$ is ergodic.

Let $\delta = \delta(v) = \frac14 \min\{\inj(\pi(v)), d(v,\partial N)\}$, and
let $V$ be the $\delta$-neighborhood of $v$ in $T^1N$.
For $w\in \Lambda_2\cap V$,  consider the set 
$$
N_\delta(w) =\Sat^u\left( \varphi_{(-\delta,\delta)}\left(\cW^s(w,\delta)\right),\delta \right);
$$
We claim:
\begin{itemize}
\item[(a)] for almost every $w\in \Lambda_2\cap V$, $B(f)$ is almost everywhere constant on $N_\delta(w)$;
\item[(b)] there is a neighborhood $U_v\subset V$ of $v$ 
such that for almost every $w\in U_v$,  the set
$N_\delta(w)\cap U_v$ has full measure in $U_v$. 
\end{itemize}
Together, these statements imply that there is a neighborhood $U_v$ of $v$ on which
$B(f)$ is $\text{a.e.}$ constant, completing the proof of ergodicity.

We first establish part (a) of this claim.   Let $G$ 
be the full measure subset of vectors in $\Lambda_2$ where the limit (\ref{e=ergav}) exists and $f^u=f^s=B(f)$.  The absolute continuity property (1) of $\cW^s$ in Proposition~\ref{p=achoro} implies that for almost every $w\in V\cap \Lambda_2$,
the intersection $G\cap \cW^s(w,\delta)$ has full volume in $\cW^s(w,\delta)$ (that is, its complement has induced volume $0$). Fix such a $w$.  On $\cW^s(w,\delta)$, $f^s$ takes a constant value $f^s \equiv a$.  On the full volume subset $G\cap \cW^s(w,\delta)$,
$f^u$ coincides with $f^s$ and therefore also takes the constant value $a$. Since $f^u$ is $\varphi_t$-invariant, and $\varphi_t$ is a $C^\infty$ flow, $f^u$ takes the constant value $a$ on the full measure subset  $G':= \varphi_{(-\delta,\delta)}\left(G\cap \cW^s(w,\delta)\right)$ of the $n$-dimensional $C^\infty$ submanifold $ D = \varphi_{(-\delta,\delta)}\left(\cW^s(w,\delta)\right)$. 

But $f^u$ is constant along $\cW^u$ manifolds and so takes the constant value $a$ on $\Sat^u(G', \delta)$.  Since $\cW^u$ satisfies the absolute continuity property (2) in Proposition~\ref{p=achoro}, and $G'$ has full measure in $D$, it follows that
$\Sat^u(G',\delta)$ has full measure in $\Sat^u(D,\delta) = N_\delta(w)$.  Hence $f^u$ is constant on a full measure subset of $N_\delta(w)$.  Since $f^u = B(f), \text{ a.e.}$, it follows that   $B(f)$ is almost everywhere constant on $N_\delta(w)$, proving part (a).

We next establish part (b) of the claim.   Let $\cC^u$ and $\cC^s$ be the closed, continuous conefields spanning $\dot\varphi^\perp$ over $T^1N$  that are given by  Corollary~\ref{c=cones}. 
For $\ast\in\{u,s\}$, the absolute  continuity property (1) of $\cW^\ast$ implies that
for almost every $w\in \Lambda_2\cap V$, the disk $\cW^\ast(w,\delta)$ is almost everywhere tangent to $E^\ast$, which by  Corollary~\ref{c=cones} is contained in the continuous conefield $\cC^\ast$. Hence for almost every $w$,  the tangent bundle $T(\cW^\ast(w,\delta))$ is everywhere contained in $\cC^\ast$.
The invariance of $\cW^s$ under $\varphi_t$ implies that for almost every $w\in  \Lambda_2\cap V$, the tangent bundle to the disk  $ D(w) = \varphi_{(-\delta,\delta)}\left(\cW^s(w,\delta)\right)$ is everywhere contained in $\cC^s\oplus E^0$.  The line field $E^0 = \RR\dot\varphi$ is smooth on the whole of $T^1N$, while $E^u \oplus E^s$ at any $v$ is the orthogonal complement of $E^0$ everywhere that the subspaces $E^u$ and $E^s$ are defined.  By Corollary~\ref{c=cones}, the conefields $\cC^u$ and $\cC^s$ intersect only at $0$.
It follows that there exists a neighborhood $U_v\subset V$ of $v$ such that for any $w, w'\in \Lambda_2\cap U_v$:
$$\cW^u(w',\delta)\cap D(w)\neq \emptyset;
$$
in other words, for every $w\in \Lambda_2\cap U_v$, the set
$N_\delta(w) = \Sat^u(D(w),\delta)$ intersects $\Lambda_2\cap U_v$ in a full measure subset.
This completes the proof of part (b) of the claim, and the proof of ergodicity.
\end{proof}

\begin{proof}[Proof of the Bernoulli property] Recall that a {\em contact form} on a $2n+1$ dimensional manifold is a differential one-form $\beta$ with the property that $\beta\wedge (d\beta)^n$ is nondegenerate.  A {\em contact flow} is a flow that preserves a contact form.
It is a well-known fact that every geodesic flow $\varphi_t$, when restricted to the unit tangent bundle,  is a contact flow; the one-form that assigns the value $1$ to $\dot\varphi$ and vanishes on $\dot\varphi^\perp$ is contact and is $D\varphi_t$-invariant.  This follows from the fact that $\varphi_t$ preserves the symplectic form on the full tangent bundle and that $\dot\varphi^\perp$ is $D\varphi_t$-invariant.

Theorem 3.6 of \cite{KatokBurns} states that any ergodic, nonuniformly hyperbolic contact flow defined on an invariant, positive volume subset of a compact contact manifold is Bernoulli on that subset.  Compactness is a simplifying assumption in the proof, and the same proof works for a nonuniformly hyperbolic contact flow that satisfies the conclusions of Proposition~\ref{p=aclam}. Returning to the context of Theorem~\ref{t=generalergodicity}, we have just proven that the geodesic flow is nonuniformly hyperbolic and ergodic.  Since it is contact, it is therefore Bernoulli.
\end{proof}

This completes the proof of the ergodicity/Bernoulli conclusion in  Theorem~\ref{t=generalergodicity}. In Appendix B, we complete the verification of the hypotheses of \cite{KS} and prove the conclusion that $\varphi_t$ has finite, positive entropy.\end{proof}

\section{Bounds on the derivative of $\varphi_1$ in the WP metric}\label{s=firstderivWP}

%{The WP metric is negatively curved, so the analysis in the previous section applies in this setting. In this section, we fix a Riemann surface $S$ and denote by $\cT$ the Teichm\"uller space $\Teich(S)$, by $\bar\cT$ the augmented Teichm\"uller space $\overline\Teich(S)$, and by$\partial\cT$ the boundary $\bar\cT\setminus \cT$.  We will omit the dependence on $S$ in the notation for the mapping class group $\MCG = \MCG(S)$ and curve complex $\cC = \cC(S)$. For $\sigma\in\cC $, the boundary stratum corresponding to the noded surface $X_\sigma$ is denoted by $\cT_\sigma$. We denote by $\pi\colon T\cT \to \cT$ the natural projection and by $\varphi_t$ the geodesic flow on $T^1\cT$. }

In this section we use the notation of Section~\ref{ss=dmcompact}, omitting the dependence on $S$.
For each unit WP tangent vector  $v\in T^1\cT$  and $t\geq 0$,  we denote by $\rho_t(v)$ the minimum WP distance from the geodesic segment $\pi(\varphi_{[-t,t]}(v))$ in $\cT$ to the singular locus $\partial\cT$.  If  $\varphi_{[-t,t]}(v)$ is not defined on the interval because the geodesic hits the singular locus in this time interval,  then we set $\rho_t(v)=0$. The main result of this section is:

\begin{theorem}\label{t=firstderivest} There are constants $\beta>0$, $0< \delta \leq 1$, $\rho_0 > 0$  and
$C \geq 1$ such that, if   $\tau\in [0,\delta]$ and   $v\in T^1\cT$ satisfies $\rho_t(v) \in (0,\rho_0)$, then
$$
\| D_v \varphi_\tau\|_{WP} \leq C (\rho_\tau(v))^{-\beta}.
$$
\end{theorem}

Since it will not cause confusion, we omit the subscript ``WP" from the notation for inner product, norm and distance functions in this section.  These subscripts will return in Section 5, where we need comparisons between the WP and Teichm\"uller metric.

\subsection{Bounding the derivative of the geodesic flow}

Theorem~\ref{t=firstderivest} is based on an estimate on the derivative of the geodesic flow that holds in 
any manifold with negative curvature. The estimate is not optimal, but will suffice for our purposes. There are simpler bounds on the derivative of the geodesic flow in \cite{Ma} and the appendix of \cite{BallmannBrinBurns1}, but they are not adequate for us.

\begin{theorem}\label{thm:bound} 
Let $M$ be a negatively curved manifold, and for  $\tau\leq 1$, let
$\gamma: [-\tau,\tau]\to M$ be  a geodesic. Let $\kappa:[-\tau,\tau]\to \R_{>0}$
be a Lipschitz function such that, for $-\tau \leq t \leq \tau$, the sectional curvature of any plane containing $\dot\gamma(t)$ is greater than $-\kappa(t)^2$ and let $u: [-\tau,\tau]\to [0,\infty)$ be the solution of the Riccati equation 
 $$
 u' + u^2 = \kappa^2
 $$
 such that $u(-\tau) = 0$. 
Then
 $$
 \|D_{\dot\gamma(0)}\phi_\tau\| \leq 1+ 2(1+ u(0)^2)\left(1+\sqrt{1+u(\tau)^2}\right)\exp\left(\int_0^\tau u(s)ds\right).
 $$
\end{theorem}

This theorem is proved at the end of this subsection. 
To prove Theorem~\ref{t=firstderivest} we will apply 
Theorem~\ref{thm:bound} to the WP geodesic segment 
$\gamma_v: [-\tau,\tau] \to \cT$ with a suitable choice of the function $\kappa$.
At the end of this   section in Proposition~\ref{p=controlfork} we show, using results of Wolpert, that  there are universal  constants $Q,L \geq 1$ such that if $v$ and $\tau$ satisfy  the hypotheses of 
Theorem~\ref{t=firstderivest}, then we can chose the positive Lipschitz function  
$\kappa$ of Theorem~\ref{thm:bound} to have the following properties:

\begin{enumerate}
\item [($\kappa 1$)] \label{kappa1} $\kappa$ is $Q$-controlled on $[-\tau,\tau]$, by which we mean that $\kappa$ is differentiable from the right and there is a constant $Q \geq 1$ such that
 $$
 D_R \kappa \geq  \frac{1-Q^2}Q \kappa^2.
 $$

\item [($\kappa 2$)]\label{kappa2}  There is a constant $L > 0$ such that 
 $$
 \int_{-\tau}^\tau \kappa(t)\,dt  \leq L|\ln(\rho_\tau(v))|.
 $$
\item [($\kappa 3$)] \label{kappa3} There is a constant $P>0$ such that 
$$\kappa(\tau) \leq P (\rho_{\tau}( \dot\gamma(0) ) )^{-1},$$

\end{enumerate}

Assuming  these estimates we have 
\begin{proof} [Proof of Theorem \ref{t=firstderivest}]
We first observe that if $u$ is the solution of $u'+u^2=\kappa^2$ with $\kappa$ Lipschitz and $Q$-controlled and $u(-\delta)=0$ then 
%$\kappa: [-\delta,\delta] \to \R_{>0}$ is Lipschitz and $Q$-controlled as above then  
 %for  $u$  the solution of  $u' +u^2 = \kappa^2$ such that $u(-\delta) = 0$, 
$u \leq Q\kappa$ on  $[-\delta,\delta]$.  For if 
$u(t) = Q\kappa(t)$ for some $t$, then $u'(t) \leq (1 - Q^2)\kappa^2(t) \leq D_RQ\kappa(t)$.

Now 
Theorem~\ref{t=firstderivest} follows immediately from Theorem~\ref{thm:bound}, and  
 the  estimates ($\kappa2$) and ($\kappa3$). 
 \end{proof}

\begin{proof} [Proof of Theorem \ref{thm:bound}]
Let $\cX$ and $\cY$ be the fundamental solutions of the matrix Jacobi equation 
$$
\cJ''(t) + \cR(t) \cJ(t) = 0
$$
such that $\cX'(-\tau) = 0$, $\cY(\tau) = 0$ and $\cX(0) = Id = \cY(0)$. The matrices $U(t) = \cX'(t)\cX^{-1}(t)$ and
$V(t) = \cY(t)'\cY^{-1}(t)$ are symmetric, since it is obvious that the condition given  in Section~\ref{SSMJRE} is satisfied by $\cX$ at $-\tau$ and by $\cY$ at $\tau$.
Moreover $U(-\tau) = 0$ and it  follows from  \cite[Section 1.10]{Eb}  that $U(t)$ is positive definite for each $t \in (-\tau,\tau]$.
\begin{lemma}\label{lem:uU} 
$\|U(t)\| \leq u(t)$ for each $t \in [-\tau,\tau]$.
\end{lemma}

\begin{proof}
For each unit vector $e \in \R^{\dim(M)-1}$, let $u_e(t) = \langle U(t)e,e\rangle$. Then $u_e(-\tau) = 0$ and 
$u_e > 0$ on $(-\tau,\tau]$ for each $e$. 
Since $U$ is symmetric, $\|U\| = \sup_e u_e$. 
The matrix Riccati equation, the symmetry of $U$, the assumption that $-\kappa^2$ is a lower bound for the sectional curvatures,  and Cauchy-Schwarz give
$$ 
 u_e' =\langle U' e,e\rangle = \langle\cR e,e\rangle  - \langle U^2 e, e\rangle \leq 
 \kappa^2 - \langle Ue,Ue\rangle \leq \kappa^2 -\langle Ue,e\rangle ^2= \kappa^2-u_e^2.
$$
It follows that $u_e \leq u$ on $[-\tau,\tau]$ for each $e$. Hence $\|U\| \leq u$.
\end{proof}

\begin{lemma}
\label{lem:grow}
 For any non trivial orthogonal Jacobi field $X$ such that $X'(-\tau) = 0$ we have
$$
\frac{\|(X(\tau),X'(\tau))\|_{Sas}}{\|(X(0),X'(0))\|_{Sas}}  \leq
\sqrt{1 + \|U(\tau)\|^2} \exp\left(\int_0^\tau \|U(t)\|\, dt\right).
$$
\end{lemma}

\begin{proof}
We have $\|X'\| \leq  \|U\|\|X\|$ by the definition of $U$. Hence
$$
\|(X(\tau),X'(\tau))\|_{Sas} = \sqrt{X(\tau)^2 + X'(\tau)^2} =   \|X(\tau)\|\sqrt{1 +
\|U(\tau)\|^2}.
$$
Since $\|X\|'(t) = \langle X',X/\|X\|\rangle \leq \|X'(t)\|$, we have
\begin{align*}
\frac{\|X(\tau)\|}{\|X(0)\|} &\leq \exp\left( \int_0^\tau
\dfrac{\|X'(t)\|}{\|X(t)\|}\,dt \right)
     \leq \exp\left(\int_0^\tau\|U(t)\|\,dt\right).
\end{align*}
Putting these last two inequalities together gives the desired estimate.
\end{proof}

\begin{lemma} For any orthogonal Jacobi field $Y$ such that $Y(\tau) = 0$ we have
\begin{equation} \label{Lderiv}
 \|Y'(0)\| \geq \|Y(0)\|/\tau \geq \|Y'(\tau)\|.
 \end{equation}
\end{lemma}
\begin{proof}
$\|Y\|$ is convex,
by Lemma~\ref{lem:negcon}, and decreases from $\|Y(0)\|$ to $0$ across
the interval $[0,\tau]$. Hence
 $$
 -\|Y\|'(0) \geq \|Y(0)\|/\tau \geq - \lim_{t \to \tau^-} \|Y\|'(t).
 $$
Since $\|Y\|' = \langle Y',Y/\|Y\|\rangle$, the Cauchy-Schwarz inequality  gives 
$\|Y'(0)\| \geq -\|Y\|'(0)$.
Since $Y(\tau) = 0$, we have $Y(t) = (t-\tau)Y'(\tau) + o(|t-\tau|)$ for $t$ near $\tau$,
whence
 $$
 \lim_{t \to \tau^-} \|Y\|'(t) = -\lim_{t \to \tau^-} \|Y'(t)\| =  -\|Y'(\tau)\|.
 $$ 
\end{proof}

Two immediate consequences of this lemma are:
\begin{enumerate}
\item \label{Y1} All eigenvalues of $V(0)$ are less than or equal to $-1$, and hence all eigenvalues of $U(0) - V(0)$ are greater than or equal to $1$.
\item \label{Y2} If $Y$ is as in the Lemma, then 
 $\|(Y(\tau),Y'(\tau))\|_{Sas} \leq \|(Y(0),Y'(0))\|_{Sas}$.
\end{enumerate}

We now consider an arbitrary orthogonal Jacobi field $(J,J')$ and in the next lemma decompose it as
 $$
  (J,J') = (X,X') + (Y,Y'),
  $$
  where $X'(\tau) = 0$ and $Y(\tau) = 0$. 
  
    \begin{lemma}\label{lem:Xsize}
The decomposition of the Jacobi field $(J,J')$ as $(X,X')+(Y,Y')$ as above  satisfies     $\|(X(0),X'(0))\|_{Sas} \leq  2(1+ \|U(0)\|^2)\|(J(0),J'(0))\|_{Sas}$.
  \end{lemma}

  \begin{proof}
  Let $v = J(0)$, $v' = J'(0)$ and $w =  [U(0) - V(0)]^{-1}[v' - U(0)v]$.  Then
 \begin{align*}
 (v,v') &= (v,U(0)v) + (0,v' - U(0)v)\\  
        &= (v,U(0)v) + ( w-w, [U(0) - V(0)] w) \\
        &= (v + w, U(0)(v +w))  - (w, V(0)w).
 \end{align*}
 This is the desired decomposition. Since $\|[U(0) - V(0)]^{-1}\| \leq 1$ by (\ref{Y1}) above, we obtain
  \begin{align*}
 \|(X(0),X'(0))\|_{Sas} &\leq \|v+w\|(1+ \|U(0)\|)^{1/2}\\  
                                 &\leq (\|v\| + \|v'\| +\|U(0)\|\|v\|)(1+ \|U(0)\|^2)^{1/2} \\
                                 &\leq \sqrt{2}(\|v\| + \|v'\|)(1+ \|U(0)\|^2)\\ 
                                & \leq 2\|(J(0),J'(0))\|_{Sas}(1+ \|U(0)\|^2) ,
 \end{align*}
 as desired. \end{proof}
 
  Using (\ref{Y2}) above we see that
  \begin{align*}
  \|(J(\tau),J'(\tau))\|_{Sas} &\leq \|(X(\tau),X'(\tau))\|_{Sas} +  \|(Y(\tau),Y'(\tau))\|_{Sas} \\
                                          &\leq \|(X(\tau),X'(\tau))\|_{Sas} +  \|(Y(0),Y'(0))\|_{Sas} \\
                                          &\leq \|(X(\tau),X'(\tau))\|_{Sas} +  \|(J(0),J'(0))\|_{Sas}  
                                                                 +  \|(X(0),X'(0))\|_{Sas}.
  \end{align*}
  The theorem now follows from 
   Lemmas~\ref{lem:uU}, \ref{lem:Xsize} and \ref{lem:grow}.
\end{proof}

The remainder of this section is devoted to work leading up to the proof of Proposition~\ref{p=controlfork}, whose proof will conclude that of Theorem~\ref{t=firstderivest}.
We begin with the next two subsections which summarize  work of Wolpert in an important constellation of papers \cite{WolpertExtension,WolpertGeometry, WolpertUnderstanding, WolpertBehavior}.

\subsection{Combined length bases}

Wolpert's precise estimates for the WP metric are stated in terms of a local system of vector fields on $\cT$ that are especially adapted to the pinched curves in nearby strata.  To define this system
of vector fields in a neighborhood  in $\cT$ of a stratum $\cT_\sigma$, where $\sigma\in \cC$, one first chooses carefully a complementary collection of curves $\chi$  (disjoint from the curves in $\sigma$, but not necessarily from each other) so that
the length functions $\ell_\beta$ for $\beta\in \chi$ give local coordinates on $\cT_\sigma$.  The pair $(\sigma, \chi)$ is called a combined length basis. Having found a combined length basis $(\sigma,\chi)$, the vector fields in the neighborhood of $\cT_\sigma$ are defined using the almost complex structure $J$ and the length functions $\ell_\alpha$ and $\ell_\beta$, for $\alpha\in\sigma$ and $\beta\in \chi$.   For the purposes of our arguments, it is important that these choices be made uniformly.  Here we describe Wolpert's construction of combined length basis and explain how they can be chosen in a uniform manner by using the compactness of $\overline\cM$.

If $\chi$ is an arbitrary finite collection of vertices in $\cC $ and $X\in \cT$, we define:
$$\uell_\chi(X)  = \min_{\beta\in\chi} \ell_\beta(X),\, \text{ and }\,
\oell_\chi(X)  = \max_{\beta\in\chi} \ell_\beta(X).$$
 For $X\in \cT$, we continue to denote by $\uell(X)$ the systole of $X$, which is the length of the shortest closed hyperbolic geodesic in $X$.
%By a ``simple closed curve in $S$" we will always mean a vertex of $\cC $; that is, the homotopy class of a nonperipheral, %homotopically nontrivial,  simple closed curve in $S$.   
%Two curves are {\em disjoint} if they
%are connected by an edge in $\cC $. 
 Let $\cB$ be the set of pairs $(\sigma,\chi)$, where $\sigma\in\cC $ and  $\chi$  is a collection of  simple closed curves in $S$ such that each $\beta\in\chi$ is disjoint from every $\alpha\in\sigma$ (we allow for the possibility that $\chi=\emptyset$). 

For each simple closed curve $\alpha$ in $S$, the {\em root length  function} $$\ell_\alpha^{1/2} \colon \cT \to \RR_{>0}$$ 
plays an important role in various asymptotic expansions of the WP metric.  Wolpert proved that the functions $\ell_\alpha$ and
 $\ell_\alpha^{1/2}$ are convex along WP geodesics in $\cT$ (see Corollary 3.4 and Example 3.5 of \cite{WolpertBehavior} and Corollary 8.2 of \cite{Wolf09}).  In \cite{Wolf} Wolf  gave another proof of the convexity of $\ell_\alpha$.
The WP gradient of $\ell_\alpha^{1/2}$ defines a vector field
$$\lambda_{\alpha} =\grad \ell_{\alpha}^{1/2}.$$ 
Following Wolpert, we say that  $(\sigma,\chi)\in \cB$ is a {\em combined (short and relative) length basis at $X\in \cT$} if  the collection $$\{\lambda_\alpha(X), J\lambda_\alpha(X), \grad\ell_\beta(X)\}_{\alpha\in\sigma,\beta\in\chi}$$ 
is a basis for $T_X\cT$.  

For each $\eta >  0$, let
$$U(\eta)  = \{X\in \cT \,\mid\, \uell (X)<\eta\},$$
which is a deleted open neighborhood of $\partial\cT$ in $\bar\cT$. 

\begin{proposition}\label{p=thickthin}  There exist constants $c>1$, $\eta, \delta>0$ and a countable collection  $\cU$ of open sets
in $\cT$ with the following properties. 
\begin{enumerate}
\item For each $U\in \cU$, there exists a combined length basis $(\sigma, \chi)\in \cB$
such that, for every $X\in U$:
$$ 1/c <  \underline \ell_{\chi}(X) \leq \overline \ell_{\chi}(X) < c.
$$
\item For each $X\in U(\eta)$, there exists  $U\in \cU$  such that for any $Y\in \cT$, 
$$d(X,Y) < \delta\,\implies\,  Y\in U;$$
in particular, the sets in $\cU$ cover $U(\eta)$.
\end{enumerate}
\end{proposition}

Before proving this proposition, we discuss further the properties of the WP metric in a neighborhood of
the boundary strata of $\cT$.
Let $\sigma\in \cC $ be a simplex,  and consider a marked noded Riemann surface 
$f\colon S\to X_\sigma$ representing an element of the boundary stratum $\cT_\sigma$. 
Recall that the hyperbolic surface $\hat X_\sigma$ is obtained from $X_\sigma$ by deleting its nodes.  If  $\beta$ is a  simple closed curve in $S$ that is disjoint from the curves in $\sigma$, then $f_\ast[\beta]$ is uniquely  represented as a closed geodesic on $\hat X_\sigma$.   In this way, the definition of $\ell_\beta$ extends continuously to the boundary stratum $\cT_\sigma$; for such $\beta$, we define $\ell_\beta([f\colon S\to X_\sigma])$ to be the hyperbolic length of the geodesic representative 
of $f_\ast[\beta]$ on $\hat X_\sigma$.   For $X_\sigma\in \cT_\sigma$, we can also define a relative systole $\uell(X_\sigma)$ 
to be the infimum of $\ell_\beta(\hat X_\sigma)$,  taken over all curves $\beta$ disjoint from the curves in $\sigma$.

Recall that the boundary stratum $\cT_\sigma$ is isomorphic to a  product of 
Teichm\"uller spaces.  In particular,  
$\cT_\sigma$ itself carries a WP metric,
which is the product of the WP metrics on the Teichm\"uller spaces of the pieces of $X_\sigma$, for any
$X_\sigma\in \cT_\sigma$.  
 We say that $\chi$ is a {\em relative length basis at $X_\sigma$} if $(\sigma,\chi)\in\cB$ and the functions $\{\ell_\beta\}_{\beta\in\chi}$ 
give local coordinates for $\cT_\sigma$ at $X_\sigma$.  
Equivalently, $\chi$ is a relative length basis at $X_\sigma$ if the vectors
$\{\grad\ell_\beta(X_\sigma)\}_{\beta\in\chi}
$ in the induced $WP$ metric on $\cT_\sigma$
span the tangent space $T_{X_\sigma}\cT_\sigma$.  
The following proposition is well-known; see, for example, Section 4 of \cite{WolpertBehavior}.
\begin{proposition}[Existence of relative length bases]\label{l=wolpertexistence} 
For each $\sigma\in \cC $,  and each marked noded Riemann surface 
$X_\sigma \in \cT_\sigma$, there exists $(\sigma,\chi)\in \cB$ such that
$\chi$ is a relative length basis at $X_\sigma$.
\end{proposition}

We remark that, unlike Fenchel-Nielsen coordinates, the local coordinates $\{\ell_\beta\}_{\beta\in\chi}$ never extend to a global coordinate system on $\cT_\sigma$; the reason is that there are points in $\cT_\sigma$ where the   geodesic representatives of the curves in $\chi$ cross each other orthogonally.  At these points, the coordinate system hits a singularity.  Proposition~\ref{l=wolpertexistence}  ensures, however, that if one works locally these issues can be ignored. Wolpert proves:
\begin{theorem}[\cite{WolpertBehavior}, Corollary 4.5.] \label{t=relativeWPexpansion}
The WP metric is comparable to a sum of differentials of
geodesic-length functions for a simplex $\sigma$ of short geodesics and corresponding
relative length basis $\chi$ as follows
$$\langle \quad,\quad\rangle \asymp
\sum_{\alpha\in\sigma}
(d\ell_\alpha^{1/2})^2
+(d\ell_\alpha^{1/2}\circ J)^2 +
\sum_{\beta\in\chi} (d\ell_\beta)^2,$$
where, given $X_\sigma \in \cT_\sigma$  and $\chi$ there is a neighborhood $U$ of $X_\sigma$ in 
$\bar\cT$ in which the comparison holds uniformly.
\end{theorem}

This has the immediate corollary:
\begin{corollary}\label{c=extend} If $\chi$ is a relative length basis at $X_\sigma \in \cT_\sigma$, then there 
is a neighborhood $V$ of $X_\sigma$ in $\bar\cT$  such that for every $X\in V \cap \cT$,
$(\sigma,\chi)$ is a combined length basis at $X$.
\end{corollary}

%We introduce more notation that will be used in the proof of  Proposition~\ref{p=thickthin}. For $\sigma\in\cC $  and %$\eta>0$, we set 
%$$U(\sigma, \eta)  = \{X\in \cT \,\mid\, \oell_\sigma (X)<\eta\},$$
%so that
%$$U(\eta) = \bigcup_{\sigma\in\cC } U(\sigma, \eta) .$$

\begin{proof}[Proof of Proposition~\ref{p=thickthin}]

Let $P\colon \bar\cT \to \bar\cM$ be the quotient map from $\bar\cT$ to the Deligne-Mumford compactification $\bar\cM$ under the action of the  mapping class group $\MCG$.  
Note that $P(U(\eta))$ is a deleted open neighborhood of $\partial\cM$ in $\bar\cM$.
Since the action of the mapping class group on  $\cC$ has finitely many orbits, we can choose a finite number of simplices $\sigma_1,\dots,\sigma_k \in \cC$ such that $\partial\cT$ is the union of the translates by the mapping class group of the sets $\T_{\sigma_1},\dots, \T_{\sigma_k}$.

For each  $X \in\T_{\sigma_i}$,  we can choose a simplex $\chi$ such that  $(\tau,\chi) \in \cB$ and $(\tau,\chi)$ gives a combined length basis at each point of some neighborhood $U(X)$ of $X$; this is Corollary 3.5. The neighborhood can be chosen small enough so that there is a constant $c(X) > 1$ such that 
 $$
 1/c(X) < \underline\ell_\chi(Y) \leq \overline\ell_\chi(Y) < c(X)
 $$
 for all $Y \in U(X)$. 
Since $\bar \cM$ is compact, we can choose a finite number of points $X_1,\dots, X_N$  such that the sets $PU(X_i)$ cover  $\bar \cM$. 
The set $\U$ in the statement of Proposition 3.2 can be chosen to be the collection of all translates by elements of the mapping class group of the sets $U(X_i)$.  The desired constant $c$ is the maximum of the constants $c(X_i)$.
Part (2) is obvious from the way in which  $\U$ was chosen. We started with a finite cover of a compact set and then translated them by the mapping class group.
\end{proof}

\subsection{First and second order properties of the WP metric}\label{s=wprops}

For each $c> 1$,  and $(\sigma,\chi)\in \cB$,  let
$$\Omega(\sigma,\chi, c )=\{X\in \cT  \mid  \oell_{\sigma\cup\chi} (X) <  c,
\text{ and }  1/c < \uell_\chi(X) \, \}.$$
Wolpert proved key estimates on the WP metric in $\Omega(\sigma,\chi, c )$, which we summarize
in the following three propositions.

The first set of estimates expands upon and refins the statement in Theorem~\ref{t=relativeWPexpansion}.
\begin{proposition}[First order estimates]\cite{WolpertExtension}\label{p=wolpertest1} Fix $c > 1$. 
For any $(\sigma,\chi) \in \cB$, the following estimates hold uniformly on $\Omega(\sigma, \chi, c )$:
\begin{enumerate}
\item  if $\alpha,\alpha' \in \sigma$, then
$$\langle J \lambda_\alpha, J\lambda_{\alpha'} \rangle =\langle\lambda_{\alpha},\lambda_{\alpha'} \rangle  = \frac{1}{2\pi}\delta_{\alpha,\alpha'}+O((\ell_{\alpha}\ell_{\alpha'})^{3/2});
$$
\item if $\alpha, \alpha' \in\sigma$ and $\beta\in\chi$, then
$$\langle \lambda_{\alpha}, J\lambda_{\alpha'} \rangle =\langle J\lambda_\alpha,\grad\ell_\beta\rangle =0;
$$  
\item if $\beta,\beta' \in\chi$, then
$$
\langle \grad \ell_{\beta},\grad\ell_{\beta'}\rangle \asymp 1;
$$ 
moreover, $\langle \grad \ell_{\beta},\grad\ell_{\beta'}\rangle$ extends continuously
to $\cT_\sigma$; 
\item if $\alpha\in\sigma$ and $\beta\in\chi$, then
$$\langle \lambda_\alpha,\grad\ell_{\beta} \rangle = O(\ell_\alpha^{3/2});
$$
\item if  $X\in \Omega(\sigma,\chi,c)$, then
$$
d(X, \cT_\sigma)  =  \left(2\pi  \sum_{\alpha\in \sigma} \ell_{\alpha}(X)\right)^{1/2} + O  \left(\sum_{\alpha\in \sigma} \ell_{\alpha}^{5/2}(X)\right) .
$$
\end{enumerate}
\end{proposition}

The second set of Wolpert's estimates are formulae for covariant derivatives, which are described in the next proposition. In each formula in the next proposition, the  error term is a vector, and the expression $v = O(a)$ means that the WP length of $v$ is $O(a)$.
\begin{proposition}[Second order estimates]\cite{WolpertExtension}\label{p=wolpertest2}  Fix $c>1$. For any $(\sigma,\chi) \in \cB$, the following estimates hold uniformly on $\Omega(\sigma, \chi, c)$:
\begin{enumerate}
\item for any vector $v \in T\Omega(\sigma,\chi, c)$, and  $\alpha\in\sigma$, we have
$$
\nabla_{v}\lambda_{\alpha} = \frac{3}{2\pi  \ell_\alpha^{1/2}} \langle v, J\lambda_{\alpha} \rangle J\lambda_{\alpha}+O( \ell_\alpha^{3/2} \|v\|_{WP});
$$
\item for $\beta\in\chi$ and  $\alpha\in\sigma$, we have
$$
\nabla_{\lambda_{\alpha}} \grad\ell_\beta = O(\ell_\alpha^{1/2}), \quad \nabla_{J\lambda_{\alpha}} \grad \ell_\beta=O(\ell_\alpha^{1/2});
$$ 
\item for $\beta, \beta'\in \chi$, $\nabla_{\grad \ell_{\beta}}\grad\ell_{\beta'}$ extends continuously to $\cT_\sigma$. 
\end{enumerate}
\end{proposition}

The final set of Wolpert's estimates we use involve the WP  curvature tensor.
\begin{proposition}[Bounds on curvature]\cite{WolpertExtension}\label{p=wolpertest3}     Fix $c>1$. For any $(\sigma,\chi) \in \cB$, the following estimates hold uniformly on $\Omega(\sigma, \chi, c)$.
For all  $\alpha\in\sigma$ we have
 \begin{equation}
\label{eq:curve}\langle R(\lambda_\alpha, J\lambda_\alpha)J\lambda_\alpha, \, \lambda_\alpha
 \rangle=\frac{3}{16\pi^2 \ell_\alpha}+O(\ell_\alpha).
\end{equation}
Moreover for any quadruple  $(v_1,v_2,v_3,v_4) \in \{\lambda_\alpha, J\lambda_\alpha,\grad\ell_\beta\}_{\alpha\in\sigma, \beta\in\chi}^4$ that is not a curvature-preserving permutation of $(\lambda_\alpha,J\lambda_\alpha, J\lambda_\alpha, \lambda_\alpha)$ for some $\alpha\in\sigma$, 
we have: \begin{equation}
\label{eq:others}
\langle R(v_1,v_2)v_3, v_4 \rangle =O(1).\footnote{Wolpert actually proves more:
each vector $v_i$ appearing in this expression that is of the form $\lambda_\alpha$ or $J\lambda_\alpha$ introduces a multiplicative  bound $o(\ell_\alpha)$ in the curvature tensor.   This means that there are sectional curvatures that are arbitrarily close to  $0$.}
\end{equation}
\end{proposition}

\subsection{Curvature estimates along a geodesic}

Fix a  unit speed WP geodesic $\gamma \colon I\to \cT$ in Teichm\"uller space. For each simple closed curve $\alpha$ we define
functions $f_\alpha = f_{\alpha,\gamma} \colon  I\to \RR_{>0}$ and  $r_\alpha = r_{\alpha,\gamma} \colon  I\to \RR_{>0}$ 
by 
$$f_\alpha(t) = \ell_\alpha^{1/2}(\gamma(t)),\, \text{ and } \,
  r^2_\alpha(t) =  \langle \lambda_\alpha,\dot\gamma(t)\rangle^2 + 
                           \langle J\lambda_\alpha,\dot\gamma(t)\rangle^2.
 $$
Roughly, $r_\alpha$ measures  the speed of the geodesic $\gamma$ in the complex line field  spanned by
$\{\lambda_\alpha,  J\lambda_\alpha\}$.
Wolpert used the function $r_\alpha$ to study the behavior of geodesics terminating in the boundary strata of $\cT$.
We will use $r_\alpha$ and $f_\alpha$ to bound sectional curvatures along $\gamma$.  We summarize in the next few lemmas
the key properties of $r_\alpha$ and $f_\alpha$ that will be used in the sequel.

The first property is an immediate consequence of part (5) of Proposition~\ref{p=wolpertest1} and explains
the significance of the quantity $f_\alpha$.
\begin{lemma}\label{l=cfomptodist} Fix $c>1$.  For every $(\sigma,\chi)\in\cB$ and any $\gamma$, 
if $\gamma(t)\in\Omega(\sigma,\chi,c)$,
then 
$$
d(\gamma(t), \cT_\sigma) =  \left(2\pi \sum_{\alpha\in\sigma} f_\alpha^2(t)\right)^{1/2} + O(\sum_{\alpha\in\sigma}f_\alpha^5 (t) ).
$$
\end{lemma}

The next two lemmas will allow us to bound the variations of $r_\alpha$ and $f_\alpha$ along a geodesic.  As was pointed out to us by  Scott Wolpert, the next lemma can be seen 
as the WP analogue of the first Clairaut equation for the model surface of revolution for $y=x^3$  discussed in the Introduction (see \cite{WolpertUnderstanding}).
 
\begin{lemma}\label{l=rprimeest} Fix $c>1$.  For every $(\sigma,\chi)\in\cB$ and any $\gamma$, 
if $\gamma(t)\in\Omega(\sigma,\chi,c)$, then
 $$r_{\alpha}'(t) = O(f_\alpha^3(t)),$$
for every $\alpha\in\sigma$.
\end{lemma}

\begin{proof}
Since the WP metric is  K\"ahler, the almost complex structure $J$ is parallel, and so we have 
\begin{eqnarray*}
2r_{\alpha}
(t)r_{\alpha}'(t)&=& 2 \langle\lambda_{\alpha}, \dot\gamma(t)\rangle \langle \frac{D}{\partial t}\lambda_{\alpha}, \dot\gamma(t)\rangle +
2 \langle J\lambda_{\alpha}, \dot\gamma(t)\rangle \langle J \frac{D}{\partial t}\lambda_{\alpha}, \dot\gamma(t)\rangle.
\end{eqnarray*}

By part (1) of Proposition~\ref{p=wolpertest2},  we have
 $$
 \frac{D}{\partial t}\lambda_{\alpha} = \langle\dot\gamma,J\lambda_{\alpha}\rangle\frac{3}{2\pi f_\alpha} J\lambda_{\alpha} + O(f_\alpha^3)
\quad
\text{and}
\quad
J \frac{D}{\partial t}\lambda_{\alpha} =  - \langle\dot\gamma,J\lambda_{\alpha}\rangle\frac{3}{2\pi f_\alpha} \lambda_{\alpha} + O(f_\alpha^3).
$$
Plugging this into the formula for $2r_{\alpha}(t)r_{\alpha}'(t)$, and noting that
$$
\max\{|\langle\lambda_{\alpha}, \dot\gamma\rangle|, |\langle J\lambda_{\alpha}, \dot\gamma\rangle|\} < r_{\alpha},
$$
 we get:
$$
2r_{\alpha}(t)r_{\alpha}'(t) =  \frac{3}{\pi f_\alpha}\langle\lambda_{\alpha},\dot\gamma\rangle\langle\dot\gamma, J\lambda_{\alpha}\rangle^2  - 
  \frac{3}{\pi f_\alpha} 
 \langle\lambda_{\alpha},\dot\gamma\rangle\langle\dot\gamma, J\lambda_{\alpha}\rangle^2+  O(r_{\alpha}f_\alpha^3) =O(r_{\alpha} f_\alpha^3).
$$
 \end{proof}

\begin{lemma}\label{l=festimate} 
 Fix $c>1$.  For every $(\sigma,\chi)\in\cB$ and any $\gamma$, 
if $\gamma(t)\in\Omega(\sigma,\chi,c)$, then
 $$
 r_{\alpha}^2(t) = (f_\alpha'(t))^2+\frac{2\pi}{3}f_\alpha(t) f_\alpha''(t) 
                           + O(f_\alpha^4(t)),
 $$
for every $\alpha\in\sigma$.
\end{lemma}

\begin{proof}  Since  $\lambda_{\alpha} = \grad \ell_\alpha^{1/2}$, it follows
that 
$$f_\alpha' = \langle \lambda_{\alpha}, \dot\gamma \rangle.$$ 
Differentiating this expression, we obtain using part (1) of Proposition~\ref{p=wolpertest2}:
$$
f_\alpha'' = \frac{d}{dt} \langle \lambda_{\alpha}, \dot\gamma \rangle = \langle \nabla_{\dot\gamma}\lambda_{\alpha}, \dot\gamma \rangle =
 \frac{3}{2\pi f_\alpha(t)} \langle \dot\gamma, J\lambda_{\alpha}\rangle^2 + O(f_\alpha^3).
$$
Now multiply this last expression by $\frac{2\pi}{3}f_\alpha$ and add it to the above expression for $f_\alpha'^2$.  The result then follows from the definition of $r_{\alpha}^2$. 
\end{proof}

Let $$\ok^2(t) = \sup_{v\in T^1_{\gamma(t)}\cT} -\langle R(v,\dot\gamma(t))\dot\gamma(t), v\rangle.$$ 

We next bound $\ok^2$ in terms of $r_\alpha$ and $f_\alpha$.
%Specifically, 
%Roughly speaking, the curvature  is bounded by the sum of the squares of projections of $\dot\gamma$ divided by length. 
%the preceding estimates on $r_\alpha$  and $f_\alpha$ will then allow us to bound
%solutions to the differential inequality $u' \leq \ok^2 - u^2$.

\begin{lemma}\label{l=ksquaredbound} Fix $c>1$. 
 For any $(\sigma,\chi)\in\cB$ and any unit speed geodesic $\gamma$, 
if $(\sigma,\chi)$ is a combined length basis in $U\subset \Omega(\sigma,\chi,c)$,
and $\gamma(t)\in U$, then
$$\ok^2(t)=\sum_{\alpha\in\sigma}  O\left(\frac{r_\alpha^2(t)}{f_\alpha^2(t)}\right).$$
\end{lemma}

\begin{proof} 
Since $(\sigma,\chi)$ is a combined length basis, we can write $v\in T^1 \Omega(\sigma,\chi,c)$ and $\dot \gamma$ as 
$$v=\sum_{\alpha_\in \sigma} (a_\alpha \lambda_\alpha+ b_\alpha J\lambda_\alpha)+\sum_{\beta\in\chi} c_\beta \grad\ell_{\beta}$$ and $$\dot\gamma=\sum_{\alpha_\in \sigma}  (A_\alpha \lambda_\alpha+ B_\alpha J\lambda_\alpha)+\sum_{\beta\in\chi} C_\beta \grad\ell_{\beta}.$$ 
Now $v$ and $\dot\gamma$ are  unit vectors, the above estimates on the metric say that all coefficients $a_\alpha,b_\alpha,c_\beta,A_\alpha,B_\alpha,C_\beta$ are $O(1)$. 
Moreover by these same estimates and the definition of $r_\alpha$, we have 
$$r_\alpha^2=\frac{1}{4\pi^2}(A_\alpha^2+B_\alpha^2)+O(f_\alpha^3).$$ 

It now follows from Proposition~\ref{p=wolpertest3} that
\begin{eqnarray*}
-\langle R (v,\dot\gamma) \dot\gamma, v\rangle  & = & -\sum _{\alpha \in\sigma} (a_\alpha^2B_\alpha{^2}+A_\alpha{^2}b_\alpha^2)  \langle R(\lambda_\alpha, J\lambda_\alpha)  J\lambda_\alpha , \lambda_\alpha\rangle +
O(1)\\ &=&\sum_{\alpha\in \sigma} O\left(\frac{r_\alpha^2}{f_\alpha^2}\right)+O(1).
\end{eqnarray*}
\end{proof}

\subsection{Estimates on $r_\alpha/f_\alpha$}

%In the next two  subsections, we establish how bounds on $\ok^2$ 
%can be used to bound solutions to the Riccati inequality  $u' \leq \ok^2 - u^2$.
%The main result in this subsection, the next proposition,  will give us control over $\ok^2$.

We now estimate $r_\alpha/f_\alpha$; in vew of the previous lemma, this 
will give us control over $\ok^2$.

\begin{proposition}\label{p=roverfbound}  Fix  $c>1$.
There is a constant $A = A(c) > 0$ such that for any $(\sigma,\chi) \in \cB$,
for any unit speed WP segment $\gamma\colon[-\delta, \delta] \to \Omega(\sigma,\chi,c)$, with $0\leq\delta\leq 1$, and any 
$\alpha\in\sigma$, we have 
  $$
  \frac {r_\alpha(t)}{f_\alpha(t)} \leq A \max\left(1, \frac {r_\alpha(t_0)}
                          { r_\alpha(t_0)|t - t_0| + f_\alpha(t_0)} \right) 
          \qquad\text{for $0 \leq t \leq \delta$,}
  $$
where $t_0$ is the unique time in $[0,\delta]$ such that $f_\alpha(t) \geq f_\alpha(t_0)$ for $0 \leq t \leq \delta$.
\end{proposition}

\begin{proof}[Proof of Proposition~\ref{p=roverfbound}]   
The time $t_0$ is uniquely defined since $f_\alpha(t)$ is a convex function of $t$. It will suffice to prove the  proposition under the additional assumption that $f_\alpha(t)$ is increasing for $ t \geq 0$.  If    $t_0 = 0$, we apply this restricted form of the  proposition directly to the geodesic $\gamma$; if $t_0 = \delta$, we apply it to the geodesic $t \mapsto \gamma(\delta - t)$;  and if $0 < t_0 < \delta$, we consider both of the geodesics $t \mapsto \gamma(t -t_0)$ and $t \mapsto \gamma(t_0 -t)$.

We choose $C \geq 1$ large enough so that:
\begin{enumerate}

\item[$(C_1)$] the $O(f^4_\alpha)$ term in the equation $r^2_\alpha = (f'_\alpha)^2 + \dfrac{2\pi}3 
f_\alpha f''_\alpha + O(f^4_\alpha)$  given by Lemma~\ref{l=festimate} is at most $Cf^4_\alpha$;

\item[$(C_2)$] $\dfrac{r_\alpha}{C} \leq \dfrac12$;

\item[$(C_3)$] $|r'_\alpha| \leq Cf_\alpha^3$  (which is possible by Lemma~\ref{l=rprimeest}).
\end{enumerate}

Conditions $(C_1)$ and $(C_2)$ give a lower bound on $f''_\alpha$ when ${r_\alpha}/{f_\alpha} \geq C$ 
and $|f'_\alpha|$ is small.

\begin{lemma} \label{2ndderiv}
If $\dfrac{r_\alpha}{f_\alpha} \geq C$ and $|f'_\alpha| \leq \dfrac{r_\alpha}2$, then
$f''_\alpha \geq \dfrac3{4\pi}\dfrac{r^2_\alpha}{f_\alpha}$.
\end{lemma}

\begin{proof} By $(C_1)$ and $(C_2)$,
 $$
 r^2_\alpha = (f'_\alpha)^2 + 
                \dfrac{2\pi}3 f_\alpha f''_\alpha + O(f^4_\alpha) 
    \leq \frac{r^2_\alpha}{4}  +  \dfrac{2\pi}3 f_\alpha f''_\alpha +       
                \frac{r^2_\alpha}{4C} 
   \leq \frac{r^2_\alpha}{2}  +  \dfrac{2\pi}3 f_\alpha f''_\alpha.
 $$ 
\end{proof}

We continue with the proof of Proposition~\ref{p=roverfbound}. Recall we are assuming $t_0=0$. 
We have that $f_\alpha(t)$ is  increasing for $t\geq 0$.  We shall show that  
$$ 
\dfrac {r_\alpha(t)}{f_\alpha(t)} \leq \max\left(4C, \dfrac{32\pi r_\alpha(0)}
                                                   {f_\alpha(0) + tr_\alpha(0)}\right) 
                    \qquad\text{for $0 \leq t \leq \delta$.}
$$ 

 If $\dfrac {r_\alpha(t)}{f_\alpha(t)} \leq 4C$ for $0 \leq t \leq \delta$ we are done. Otherwise, let 
 $$
 b = \sup\{t \in [0,\delta] : \dfrac {r_\alpha(t)}{f_\alpha(t)} \geq 4C\}.
 $$
Since $\dfrac {r_\alpha(t)}{f_\alpha(t)} \leq 4C$ for $ b \leq t \leq \delta$, it will suffice to show that
\begin{equation}\label{wanted0}
 \dfrac {r_\alpha(t)}{f_\alpha(t)} \leq \dfrac{32\pi r_\alpha(0)}{f_\alpha(0) + tr_\alpha(0)} \qquad\text{for $0 \leq t \leq b$.}
 \end{equation}

 The following lemma is based on the existence of the value $b$ defined above. We show 
that the function $r_\alpha$ is approximately constant and $r_\alpha/f_\alpha$ is large on the interval $[0,b]$.

\begin{lemma} \label{ratio}
For $0 \leq t \leq b$ we have: 
\begin{enumerate}
\item[(i)] $\dfrac{r_\alpha(0)}2 \leq r_\alpha(t) \leq 2r_\alpha(0)$; 

\item[(ii)] $\dfrac {r_\alpha(t)}{f_\alpha(t)} \geq C$.
\end{enumerate}
\end{lemma}

\begin{proof} By $(C_3)$, $|r'_\alpha| \leq Cf^3_\alpha$ on the interval $[0,b]$. Since $b \leq \delta\leq 1$ and $f_\alpha$ is increasing on $[0,b]$, we have
$|r_\alpha(b) - r_\alpha(t)| \leq  Cf^3_\alpha(b)$ for $0 \leq t \leq b$. The definition of $b$ and $(C_2)$ ensure that $f_\alpha(b) \leq \dfrac{r_\alpha(b)}{2C} \leq \dfrac14$. 
Hence
 \begin{align*}
 \frac{|r_\alpha(b) - r_\alpha(t)|}{r_\alpha(b)} &\leq 
           \frac{Cf^3_\alpha(b)}{2C f_\alpha(b)} \leq \frac1{2}
                  f^2_\alpha(b) \leq  \frac1{32}.
\end{align*}
 Thus 
$\dfrac{31}{32} \leq \dfrac{r_\alpha(t)}{r_\alpha(b)} \leq \dfrac{33}{32}$ for $0 \leq t \leq b$, and (i) follows easily.
Claim (ii) follows from (i) since $r_\alpha(b)/f_\alpha(b) \geq 2C$ and $f_\alpha$  is increasing on $[0,b]$.
\end{proof}

Using this lemma we see that inequality (\ref{wanted0}) will follow if we prove
 \begin{equation} \label{wanted}
16\pi f_\alpha(t) \geq f_\alpha(0) + tr_\alpha(0)\qquad\text{for $0 \leq t \leq b$}.
 \end{equation}

Lemma~\ref{ratio}(i) ensures that $r_\alpha(0) > 0$, so we can set $a = \dfrac{f_\alpha(0)}{r_\alpha(0)}.$  Now for  $0 \leq t \leq \min(a,b)$, we have
 $$
 f_\alpha(0) + tr_\alpha(0)
    \leq f_\alpha(0) + a r_\alpha(0) = 2f_\alpha(0) \leq 2f_\alpha(t), 
 $$ 
since $f_\alpha$ is  increasing on $[0,\delta]$. This gives (\ref{wanted}) for $0 \leq t \leq \min(a,b)$.

We are done if $a \geq b$. It remains to show that
if $a \leq b$, then inequality~(\ref{wanted}) also holds for $a \leq t \leq b$. Since $f_\alpha$ is convex and (\ref{wanted}) already holds for $t = a$, it will suffice to show that if $a\leq b$ that 
\begin{equation}
\label{newwanted}
16\pi f'_\alpha(a) \geq r_\alpha(0)
\end{equation} 

 We may assume that $4 f'_\alpha(a) \leq r_\alpha(0)$, since otherwise there is nothing to prove.  
 Then $f'_\alpha(t) \leq f'_\alpha(a) \leq r_\alpha(0)/4$ for $0 \leq t \leq a$, because $f_\alpha$ is convex and increasing on $[0,a]$.  Since $a \leq b$, we can now apply Lemma~\ref{ratio}(ii)  to see that on $[0,a]$ we have
  $$
  f'_\alpha (t) \leq r_\alpha(0)/4 \leq  r_\alpha(t)/2 \quad\text{and}\quad 
                                          \dfrac{r_\alpha(t)}{f_\alpha(t)}  \geq C.
  $$
  Thus both hypotheses of Lemma~\ref{2ndderiv} are satisfied on $[0,a]$. Lemmas~\ref{2ndderiv} and \ref{ratio} give us 
 $$
 f''_\alpha(t) \geq  \frac3{4\pi}\frac{r^2_\alpha(t)}{f_\alpha(t)}    
               \geq  \frac3{16\pi}\frac{r^2_\alpha(0)}{f_\alpha(t)}
                >    \frac1{8\pi}\frac{r^2_\alpha(0)}{f_\alpha(t)}
 $$
 for $0 \leq t \leq a$. Since $f'_\alpha \leq r_\alpha(0)/4$ on $[0,a]$, we have 
 $$
 f_\alpha(a) \leq f_\alpha(0) + ar_\alpha(0)/4 = f_\alpha(0) +  f_\alpha(0)/4 < 2f_\alpha(0),
 $$
 and hence
 $$
 f''_\alpha(t) \geq \frac1{16\pi}\frac{r^2_\alpha(0)}{f_\alpha(0)}
, $$
for $0 \leq t \leq a$. 
Finally, since $f_\alpha$ is increasing on $[0,a]$, we have $f'_\alpha(0) \geq 0$ and
 $$
 f'_\alpha(a) \geq \frac{a}{16\pi}\frac{r^2_\alpha(0)}{f_\alpha(0)} 
     =  \frac{r_\alpha(0)}{16\pi},
 $$ which is the desired inequality  ~(\ref{newwanted}). 

\end{proof}

Combining Lemma~\ref{l=ksquaredbound} and Proposition~\ref{p=roverfbound} we obtain the immediate corollary:

\begin{corollary}\label{c=producecontrol}  Fix $c>1$.
There is a constant $B = B(c) > 0$ such that for any $(\sigma,\chi) \in \cB$,
if $(\sigma,\chi)$ is a combined length basis in an open set $U\subset  \Omega(\sigma,\chi,c) $
and $\gamma\colon[-\delta, \delta] \to U$ is a unit-speed WP geodesic segment,  then
 $$
\ok(t) \leq B  \max_{\alpha\in\sigma} \left(1, \frac {r_\alpha(t_\alpha)}
                          { r_\alpha(t_\alpha)|t - t_\alpha| + f_\alpha(t_\alpha)} \right) 
          \qquad\text{for $0 \leq t \leq \delta$,}
  $$
where $t_\alpha$ is the unique time in $[-\delta,\delta]$ such that $f_\alpha(t) \geq f_\alpha(t_\alpha)$ for $-\delta \leq t \leq \delta$.
\end{corollary}

\subsection{Controlled bounds on the curvature}
In this subsection we show that it is possible to choose an upper bound $\kappa$ for $\ok$ with the 
properties ($\kappa1$),  ($\kappa2$) and ($\kappa3$) used in the proof of Theorem~\ref{t=firstderivest}.
This will complete the proof of Theorem~\ref{t=firstderivest}.

We begin with some simple properties of controlled functions. 
 If $\kappa$ is $Q$-controlled, then it is $Q'$-controlled, for all $Q' > Q$. 
If  $\kappa$ is $Q$-controlled, then so is $t \mapsto \kappa(t -t_0)$ for any $t_0$, and for any $A > 0$, the function $A\kappa$ is $\frac{Q-1}{A} + 1$-controlled. 
The maximum of two $Q$-controlled functions is $Q$-controlled. Moreover $\kappa$ is $1$-controlled if $\kappa \equiv 1$ and $2$-controlled if 
$\kappa(t) = \dfrac{1}{|t| + a}$ where $a > 0$. 

\begin{proposition}\label{p=controlfork} There exist constants $P, Q, L\geq 2$ and  $\delta \in(0,1)$  such that for any positive $\delta' < \delta$ and
any geodesic segment $\gamma\colon (-\delta',\delta') \to \cT$, there exists a $Q$-controlled function $\kappa\colon  (-\delta',\delta') \to  \RR_{>0}$  such that
for every $t\in (-\delta',\delta')$:
\begin{enumerate}
\item $\ok^2(t) \leq \kappa^2(t)$, where 
 $$\ok^2(t) = \sup_{v\in T^1_{\gamma(t)}\cT}  - \langle R(v,\dot\gamma(t))\dot\gamma(t), v\rangle;$$ 
\item 
$$ \int_{-\delta'}^{\delta'}  \kappa(s)\, ds   \leq L | \ln (\rho_{\delta'}  ( \dot\gamma(0) ) )  |, $$
and
\item $$\kappa(\delta') \leq P (\rho_{\delta'}( \dot\gamma(0) ) )^{-1},$$
where $\rho_{\delta'}(\dot\gamma(0))$ is the distance from the geodesic segment $\gamma[-\delta',\delta']$ to $\partial\cT$.
\end{enumerate}
\end{proposition}

\begin{proof} 
Let $c$, $\eta$ and $\delta$ be the constants and let
$\cU$ be the collection of open sets in $\cT$  given by Proposition~\ref{p=thickthin}.
We write
$$
\cT = U(\eta)  \sqcup   \Theta;
$$
the set $\Theta = \cT\setminus   U(\eta)$ lies in the thick part of Teichm\"uller space in which the WP sectional curvatures are negative and bounded below by a constant $-b^2$.
By shrinking the value of $\delta$ if necessary, we may assume that for every $X\in \Theta$,
and $Y\in \cT$, if $d(X,Y) < \delta$, then:
 $$\sup_{v,w\in T^1_{Y}\cT}  - \langle R(v,w)w, v\rangle < b^2.$$ 
 Let $B = B(c) > 0$ be the constant given by  Corollary~\ref{c=producecontrol}.

Fix $\delta' <\delta$.
It follows from Proposition~\ref{p=thickthin} that if $\gamma\colon (-\delta', \delta') \to \cT$ is a unit-speed WP geodesic, then either $\gamma(0) \in \Theta$, or
$$\gamma (-\delta', \delta')  \subset U,
$$
for some $U\in\cU$.

If  $\gamma(0) \in \Theta$, then 
we define $\kappa$  to be the constant function  $b$.  Then by construction we have $\ok^2 \leq \kappa^2$.
Since the $WP$ distance from any point in $\cT$ to $\partial\cT$ is bounded above by a uniform constant,
it also follows that in this case:
$$ \int_{-\delta'}^{\delta'}  \kappa(s)\, ds   =  2b\delta' = O (| \ln (\rho_{\delta'}  ( \dot\gamma(0) ) )  |), $$
and
$$\kappa(\delta') = O (\rho_{\delta'}( \dot\gamma(0) ) )^{-1}.$$

Suppose on the other hand that
$\gamma(\delta',\delta' ) \subset U$, for some $U\in\cU$. 
Let $(\sigma,\chi)$ be the combined length basis in $U$ given by Proposition~\ref{p=thickthin}
satisfying:
$$ 1/c <  \uell_{\chi}(X) \leq \oell_{\chi}(X) < c,
$$
for every $X\in U$.  For $\alpha\in\sigma$, 
define $\kappa_{\alpha}\colon (-\delta',\delta')\to \RR_{>0}$ by:
$$\kappa_{\alpha}(t) = \frac {r_{\alpha}(t_{\alpha})}
                          { r_{\alpha}(t_{\alpha})|t - t_{\alpha}| + f_{\alpha}(t_{\alpha})},$$
where $t_{\alpha}$ is the unique time in $[-\delta',\delta']$ such that $f_{\alpha}(t) \geq f_{\alpha}(t_{\alpha})$ for $t\in [-\delta',\delta']$. 
Observe that $\kappa_{\alpha}$ is a $2$-controlled function and attains its maximum value of $ \frac {r_{\alpha}(t_{\alpha})} { f_{\alpha}(t_{\alpha})}$
at $t=t_{\alpha}$.

Applying Corollary~\ref{c=producecontrol}, 
we obtain  that for all $t\in (-\delta',\delta')$: 
 $$
\ok(t) \leq B  \max_{\alpha\in\sigma} \{1, \kappa_{\alpha} \}.
  $$
We define $\kappa\colon (-\delta',\delta')\to \RR_{>0}$ by:
$$\kappa =  B \max_{\alpha\in\sigma} \{1,\kappa_{\alpha}\}.$$  
Since $\kappa_\alpha$ is $2$-controlled, for each $\alpha$, it follows that $\kappa$ is $\frac{1}{B}+1$-controlled.
By its construction $\kappa$ satisfies the inequality $\ok^2 < \kappa^2$ on $(-\delta',\delta')$.  

It remains to estimate the integral of $\kappa$ over the interval $(-\delta',\delta')$. 
 Simple integration shows that
\begin{eqnarray*}
\int_{-\delta'}^{\delta'} \kappa_{\alpha}(s)\,ds  &=& O(\max\{\delta', \left |\ln (  f_{\alpha}(t_{\alpha})    ) \right| )\}),
\end{eqnarray*}
since $r_\alpha(t_{\alpha}) = O(1)$.

Note that
$f_{\alpha}(t_{\alpha})$ is the minimum value  of the function  $\ell_\alpha^{1/2}$ along the geodesic
segment $\gamma[-\delta,\delta]$.  Lemma~\ref{l=cfomptodist} implies that there exists a constant
$r>0$ such that $f_{\alpha}(t_{\alpha}) \geq r\rho_{\delta'}(\dot\gamma(0))$.  
This implies that
\begin{eqnarray*}
\int_{-\delta'}^{\delta'}  \kappa(s)\, ds &\leq&
 B \max_{\alpha\in\sigma}\{2\delta', \int_{-\delta'}^{\delta'} \kappa_{\alpha} (s)\,ds\} 
=O(|\ln ( \rho_{\delta'}(\dot\gamma(0) ) ) | ).
\end{eqnarray*}
Similarly,
\begin{eqnarray*}
\kappa(\delta') &\leq& B \max_{\alpha\in\sigma}\{1, \frac{r_\alpha(t_\alpha)}{f_\alpha(t_\alpha)}\} 
=O( \rho_{\delta'}(\dot\gamma(0) ) ) ^{-1} .
\end{eqnarray*}
\end{proof}

\section{Higher order control of the WP metric}\label{s=higherorder}

In this section, we show how to control higher order derivatives of the WP metric. This will verify Assumption IV. in Theorem~\ref{t=generalergodicity}.
  The main result in this section is
\begin{proposition}\label{p=derivcurvbounds} 
There exist $C, \beta_1>0$ such that for any 
$X_0\in \cT$,
the $WP$ curvature tensor  $R_{WP}$  satisfies:
$$ \max\{  \| (\nabla R_{WP})_{X_0} \| , \| (\nabla^2 R_{WP})_{X_0} \| \}   \leq C \rho_0^{-\beta_1}
$$
where $\rho_0 = \rho_0(X_0)$ is the distance from $X_0$ to the singular locus
$\partial \cT$.
\end{proposition}

We remark that similar bounds on higher derivatives of the $WP$ curvature tensor can also be obtained using the methods in this section.

\subsection{Estimates on the WP metric in special coordinates}

Following \cite{CTMKahler}, we introduce coordinates on $\Teich(S)$ in which 
we can bound the derivatives of the WP metric.  In this subsection we denote by $N$ the complex dimension of $\Teich(S)$.
Let $\Delta^N$ denote the Euclidean unit polydisk in $\CC^N$.  We will denote by $z = (z_1,\ldots, z_N)$ an
element of $\Delta^N$, where $z_k$ is a complex coordinate, and by $x_k = \Re(z_k)$, $y_k = \Im(z_k)$ the real coordinates.
Let $e_i$ be the vector field $\partial/\partial x_i$, for $1\leq i\leq N$, and $\partial/\partial y_{i-N}$, for $N+1 < i\leq 2N$. The main content of this subsection is the proof of the following proposition.

\begin{proposition}\label{p=metricest} There exists $C\geq 1$ such that for any $X_0\in \Teich(S)$, there is a
holomorphic embedding $\psi=\psi_{X_0} \colon \Delta^N \to \Teich(S)$ with the following properties:
\begin{enumerate}
\item $\psi(0) = X_0$;
\item  setting $G_{ij}(z) = (\psi^\ast g_{WP})_z(e_i,e_j)$ for $z\in \Delta^N$, we have
\begin{enumerate}
\item $\| G^{-1}(z)\| \leq C \uell(X_0)^{-2}$, and
\item  for any $i,j\in\{1,\ldots ,2N\}$ and any $k\geq 0$,
$$
\sup_{(\xi_1, \ldots, \xi_k) \in \{x_1,\ldots, x_N, y_1, \ldots, y_N\}^k} \left|\frac{\partial^k G_{i,j}} {\partial \xi_1 \cdots \partial \xi_k }(z) \right| \leq Ck! \,.
$$
\end{enumerate}
\end{enumerate}
\end{proposition}

We will use Proposition~\ref{p=metricest} to bound the covariant derivatives of the $WP$ curvature in terms of the the distance to the singular strata.

\begin{proof} 
The {\em Teichm\"uller cometric} on the cotangent bundle $T^*\Teich(S)$ is the
Finsler metric which is given on each cotangent space  $T_X^\ast\Teich(S)$ by the 
$L^1$ norm on $Q(X)$:
$$\|\phi\|_T = \|\phi\|_1 = \int_X|\phi|. 
$$
The Teichm\"uller norm on $T\Teich(S)$ is then induced by the standard pairing (\ref{eq:pair})
between quadratic and Beltrami differentials.  

The following lemma is proved in \cite{CTMKahler} and follows from Nehari's bound and the fact
that the Teichm\"uller and Kobayashi metrics agree on the image of the Bers embedding.
\begin{lemma} \cite[cf. Theorem 2.2 and Proof of Theorem 8.2]{CTMKahler}\label{l=psibounds} There exists $C_0\geq 1$ such that for any $X_0\in \Teich(S)$, there is a
holomorphic embedding $\psi=\psi_{X_0} \colon \Delta^N \to \Teich(S)$,
sending $0\in \Delta^N$ to $X_0$ and
such that for every $v\in T\Delta^N$, we have:
$$
\frac{1}{C_0} \|v\| \leq \| D \psi(v) \|_T \leq C_0\|v\|,
$$
where $\|\cdot\|$ is the Euclidean norm on $\Delta^N$, and $\|\cdot\|_T$ is the Teichm\"uller Finsler norm on $\Teich(S)$.
\end{lemma}

Fix a point $X_0\in \Teich(S)$, and let $\psi = \psi_{X_0}$ be the holomorphic embedding given by this lemma.  This is a holomorphic embedding satisfying part (1) of Proposition~\ref{p=metricest}.
Since the metric $g_{WP}$ on $\Teich(S)$ is K\"ahler with respect to the $2$-form $\omega_{WP}$, and $\psi$ is holomorphic,  it follows that the  pullback metric $\psi^\ast g_{WP}$ on $\Delta^N$ is K\"ahler with respect to the pullback form $\psi^\ast\omega_{WP}$ and the standard almost complex structure on $\Delta^N$. 

To establish Part (2)  of Proposition~\ref{p=metricest} we need  a comparison between the Teichm\"uller and WP metrics.  For a given Riemann surface $X$, recall that $\uell(X)$ denotes the length of the shortest simple closed curve in the hyperbolic metric. 
\begin{lemma} 
\label{prop:WPcomparenew} There exists $C>0$ such that
for any $X\in \Teich(S)$ and any tangent vector $[\mu]\in T_X\Teich(S)$, we have 
$$
\|\mu\|_{WP}\geq C\uell(X)\|\mu\|_{T}.
$$
\end{lemma}
A more refined analysis can improve the exponent of $\uell(X)$ in Lemma~\ref{prop:WPcomparenew}  to $1/2$,  but that will not be needed.  We are grateful to Scott Wolpert for suggesting the proof given here.

\begin{proof}[Proof of Lemma~\ref{prop:WPcomparenew}]  We establish the dual statement of Lemma~\ref{prop:WPcomparenew} 
in the Teichm\"uller and WP cometrics: there exists $C>0$ such that for any $\phi\in Q(X)$:
 $$
 \|\phi\|_{WP}\leq C\uell(X)^{-1}\|\phi\|_{T}.
$$
To this end, write $X=\HH^2/\Gamma$, normalized so that the covering transformation corresponding to the  shortest curve is the transformation $T(z)=\lambda z$.  Then $\log \lambda=\uell(X)$.
Fix a Dirichlet fundamental domain $D$ for the action of $\Gamma$ centered at  the point $i$.  For $\uell$ sufficiently small, by the collar lemma, the union of $\uell(X)^{-1}$ copies of $D$ contains a ball $B$ of fixed radius centered at any point $z$ of $D$.  
Then for any $\phi\in Q(X)$  the  Cauchy integral formula gives that
 $$|\phi(z)| =O\left(\int_B |\phi| \right)=O({\uell(X)}^{-1}\|\phi\|_T),$$ with the last estimate following from the fact that $B$ is covered by at most ${\uell(X)}^{-1}$ copies of $D$.

On the other hand, we can bound the $L^2$ norm by the $L^\infty$ norm as follows.
Since the hyperbolic metric $\rho$ is bounded away from $0$, the above bound for $|\phi(z)|$  
 on $D$ gives  $$\|\phi\|_{WP}^2=\int_X\frac{|\phi|^2}{\rho^2}=O(\uell(X)^{-2}\|\phi\|_T^2).$$
\end{proof}

Part (2a) of the Proposition now follows immediately from Lemma~\ref{prop:WPcomparenew}.
The proof of part (2b) uses in a crucial way results of McMullen in \cite{CTMKahler}.
Using the embedding $\psi$, we define
an embedding $\Psi \colon \Delta^N\times\Delta^N \to QF(S)$
by 
$$\Psi(z,w) =  (\psi(z) , \overline{\psi(\bar w)}).$$
Since $\psi$ is holomorphic and $X\mapsto\overline X$ is antiholomorphic,
the map $\Psi$ is holomorphic. Note that the image of the
antidiagonal $\{(z,\overline z) : z\in\Delta^N \}$ under $\Psi$ lies in the 
Fuchsian locus $F(S)\subset QF(S)$. 
Denote by ${\alpha}\colon \Delta^N \to \Delta^N\times \Delta^N$ the antidiagonal
embedding ${\alpha}(z) = (z,\bar z)$, and by $\hat{\alpha}\colon \Teich(S) \to QF(S)$ the antidiagonal
embedding $\hat{\alpha}(X) = (X,\overline{X})$.  Then we have the following commutative diagram:

\[
\begin{xy}
\xymatrixcolsep{4pc}
\xymatrixrowsep{3pc}
\xymatrix{
 & \Delta^N\times\Delta^N\ar[dr]^{\Psi} & \\
 \Delta^N\ar[ur]^{{\alpha}}\ar[dr]^{\psi} & & QF(S)\\
 & \Teich(S)\ar[ur]^{\hat{\alpha}} & \\
}
\end{xy}
\]

Note that the maps ${\alpha}$ and $\hat{\alpha}$ are not holomorphic, although their derivatives are bounded in the Euclidean and Teichm\"uller metrics, respectively.

Since $\Teich (S)$ and $QF(S)$ are complex manifolds, so are their cotangent bundles $T^*\Teich(S)$
and $T^*QF(S)$, and $T^*QF(S) = T^* \Teich (S) \oplus T^*\Teich (\overline S)$.
Fixing $Z\in \Teich(\overline{S})$ we define a map $\tau\colon QF(S) \to T^*\Teich(S)$
by:
$$
\tau(X,Y) = \sigma_{QF}(X,Y) - \sigma_{QF}(X,Z).
$$
Since $T^* \Teich(S)$ embeds as the first factor in $T^* QF(S)$, we may regard $\tau$ as a $1$-form on $QF(S)$, which by Theorem~\ref{thm:quasi} in the introduction is holomorphic and bounded in the Teichm\"uller Finsler norm on $\Teich(S)$.
%:
%$$
%%\|\tau(X,Y)\|_T = \sup \{ |\tau(X,Y)(\mu)|  :  [\mu] \in T_X \Teich(S), \|\mu\|_T = 1 \}\leq C,
%$$
%where $C$ is independent of $(X,Y)$. 
Furthermore the $1$-form  $\theta= -\hat{\alpha}^*\tau$ on $\Teich(S)$ is a primitive for the WP K\"ahler form:
$$
d(i\theta) = \omega_{WP}.
$$
Pulling the holomorphic $1$-form $\tau$ back to $\Delta^N\times \Delta^N$, 
we thus obtain a holomorphic $1$-form $\kappa = \Psi^*\tau$.
Then  $\kappa$ is  bounded in the Euclidean metric on $\Delta^N\times \Delta^N$, 
since $\tau$ is bounded in the Teichm\"uller metric, and the Euclidean metric
is comparable to the $\Psi$-pullback of the Teichm\"uller metric, by Lemma~\ref{l=psibounds}.
This bound is independent of $X_0$. Moreover, the commutativity of the diagram above implies:
\begin{lemma}  The holomorphic $2$-form $\Omega = d(i\, \kappa)$ on $\Delta^N\times \Delta^N$ satisfies
${\alpha}^\ast\Omega = \psi^\ast \omega_{WP}$, which is the K\"ahler $2$-form for the pullback metric $\psi^\ast g_{WP}$. The holomorphic $2$-form $\Omega = d(i\, \kappa)$ on $\Delta^N\times \Delta^N$ satisfies
${\alpha}^\ast\Omega = \psi^\ast \omega_{WP}$, which is the K\"ahler $2$-form for the pullback metric $\psi^\ast g_{WP}$.
\end{lemma}
We now finish the proof of  Proposition~\ref{p=metricest}.
In complex coordinates $(z_1,\ldots,z_N,w_1,\ldots,w_N)$ on $\Delta^N\times \Delta^N$ one can write
$$
\kappa =\sum_{i=1}^N a_i dz_i,
$$
 where $a_i\colon \Delta^N\times \Delta^N \to \CC$ 
are bounded holomorphic functions.  
Now 
$$
\Omega = d(i\,\kappa) =
\sum_{j,k=1}^N i\frac{\partial{a}_j}{\partial z_k}dz_k\wedge dz_j + i\frac{\partial{a}_j}{\partial w_k}dw_k\wedge dz_j,
$$ 
and so
$$
{\alpha}^\ast\Omega  
= \sum_{j,k=1}^N
 i\frac{\partial{a}_j}{\partial z_k}dz_k\wedge dz_j + i\frac{\partial{a}_j}{\partial \bar z_k}d\bar z_k \wedge dz_j. 
$$ 
The Euclidean coefficients of the K\"ahler metric $\psi^\ast  g_{WP}$ are 
hence linear combinations, with bounded coefficients, of $\partial{a}_j/ \partial z_k$ and  
$\partial{a}_j / \partial \bar z_k$, which in turn are pullbacks of the complex partial derivatives  $\partial{a}_j/\partial z_k$ and  $\partial{a}_j/\partial w_k$.  Since the ${a}_j$ are bounded holomorphic functions, Cauchy's Theorem implies that 
the derivatives  $\partial{a}_j/\partial z_k$ and  $\partial{a}_j/\partial w_k$ are bounded 
for $\|(z,w)\| < 1/2$;  it follows that the (real) partial derivatives of ${a}_i$ 
are bounded for $\|z\|<1/2$. The same applies to all higher order partial
derivatives (where the bound for the $k$th order derivatives incorporates a factor of $k!$).  By rescaling the map 
$\psi$ by a dilation, we may assume that these estimates hold for all $z\in \Delta^N$.
This completes the proof of (2).
\end{proof}

\subsection{Proof of Proposition~\ref{p=derivcurvbounds}}

Fix $X_0\in\Teich(S)$ and local coordinates $\psi = \psi_{X_0}$ as in Proposition~\ref{p=metricest}. 
For $z\in \Delta^N$, let $G(z) = G_{X_0}(z) = (G_{ij}(z))$ be the matrix
for the pullback metric $\psi^\ast g_{WP}$, and let $G^{ij}(z) = \left(G(z)^{-1}\right)_{ij}$.

The curvature tensor for $G$ can be calculated in these Euclidean coordinates using the Christoffel symbols
%:
%$$
%\Gamma^m_{ij} = \frac12 G^{km}\left( \frac{\partial}{\partial x^i} G_{kj} 
%+ \frac{\partial}{\partial x^j} G_{ik} 
%- \frac{\partial}{\partial x^k} G_{ij}
%\right)
%$$
and the Riemannian curvature tensor coefficients,
%:
%$$
%R^\ell_{ijk} = \frac{\partial}{\partial x^j}\Gamma^\ell_{ik}
%-\frac{\partial}{\partial x^k}\Gamma^\ell_{ij}\
%+\Gamma^\ell_{js} \Gamma^s_{ik}
%-\Gamma^\ell_{ks} \Gamma^s_{ij}.
%$$
all of which can be expressed as sums of products of the coefficients $G^{ij}$ and first and second order partial derivatives of the coefficients $G_{ij}$.
%At each $z\in\Delta^N$ the curvature tensor $R_z$ for the WP pullback  metric $G$ is defined by setting $R_z(e_i, e_j)e_k = \sum_\ell R^\ell_{ijk} e_\ell$.
%It is a function from $\Delta^N\times \CC^{3N}$ to $\CC^N$.
%Finally, the covariant tensor defined by
%$$
%R_{i,j,k,\ell} = \sum_{\nu} G_{\ell \nu} R^{\nu}_{ijk},
%$$
%satisfies $\la R_z(e_i, e_j)e_k, e_\ell \ra = R_{ijk\ell}$.
%The formula for the covariant derivative of a covariant tensor $T_{b_1\ldots b_s}$ is given by:
 %$$(\nabla_c T)_{b_1 \ldots b_s} = \frac{\partial}{\partial x^c}T_{b_1 \ldots b_s}
%-\,\Gamma ^d {}_{b_1 c} T_{d \ldots b_s} - \cdots - \Gamma ^d {}_{b_s c} T_{b_1 \ldots b_{s-1} d}.
%$$
Since $\|D\psi\|$ and $\|D\psi^{-1}\|$ are bounded by Lemma~\ref{l=psibounds}, the quantities  $\|(\nabla R_{WP})_{\psi(z)}\|$ and $\|(\nabla^2 R_{WP})_{\psi(z)}\|$ can therefore be bounded by a 
(universal) polynomial function of the quantities $|G^{ij}(z)|$, $|G_{ij}(z)|$ and
$$ \left|\frac{\partial^k G_{i,j}} {\partial \xi_1 \cdots \partial \xi_k }(z) \right|,$$
for $k=1,\ldots, 4$.
But Proposition~\ref{p=metricest} implies that  the entries $G^{ij}(z)$ are $O(\uell(X_0)^{-2})$ and
the entries $G_{ij}(z)$ and their first $k$ derivatives are $O(1)$; the conclusion of
 Proposition~\ref{p=derivcurvbounds} then follows. $\diamondsuit$

\section{Ergodicity and finite entropy of the WP geodesic flow}

Fix a Riemann surface $S$, and let $\cT = \Teich(S)$, $\MCG = \MCG(S)$ and $\cM = \cM(S)$.
We describe here first how the results of Section 6 can be applied to obtain ergodicity and finite entropy of the geodesic flow on the quotient $\cM^1 = T^1\cT/\MCG$.   Note that the results in Section~\ref{s=generalsetup} cannot be applied directly with $M=\cT$ and $\Gamma = \MCG$, since $\MCG$  does not act freely on $\cT$.  Our strategy is to prove ergodicity first for a finite branched cover $T^1\cT/\MCG[3]$. Here $\MCG[k]$ is the {\em level $k$ congruence subgroup}:
$$\MCG[k] = \{\phi\in \MCG  : \phi_*=0 \,\text{ acting on }\, H_1(S,\ZZ/k\ZZ)\},
$$
which is clearly a finite index subgroup of $\MCG$.  It is a well-known fact that for $k\geq 3$, $\MCG[k]$ is torsion-free and so acts freely and properly discontinuously by isometries on $\cT$
\cite{Serre}.  The quotient $T^1\cT/\MCG[k]$ has finite volume. We obtain ergodicity for the flow on $T^1\cT/\MCG[k]$ for any $k\geq 3$ using the setup of the previous section. 
\subsection{Ergodicity of the flow on $T^1(\cT/\MCG[k])$}

Fix $k\geq 3$.
To establish ergodicity and finite metric entropy of the $WP$ geodesic flow on $T^1(\cT/\MCG[k])$, we show that the assumptions I.-VI. of Theorem~\ref{t=generalergodicity} in Section~\ref{s=generalsetup} are satisfied for
$M=\cT$,  $\Gamma=\MCG[k]$ and the WP metric. We recall that the distance from $X\in \cT$ to the singular locus $\partial\cT$ is comparable to $\uell(X)^{1/2}$ (Proposition~\ref{p=wolpertest1}, part(5)).

The fact that in the Weil-Petersson metric $\cT$ is geodesically convex was proved by Wolpert \cite{WolpertBehavior}.  Since the completion $\bar\cM$ of $\cM$ is compact, and
$\cT/\MCG[k]$ is a finite branched cover of $\cM$, it follows that the completion of
$\cT/\MCG[k]$ is compact as well.  Hence assumptions I. and II. hold true.

The curvature bound in assumption IV. is due to Wolpert and was stated as Proposition~\ref{p=wolpertest3}.  The bounds on $\|\nabla R_{WP}\|$ and $\|\nabla^2 R_{WP}\|$ in assumption IV.  are the content of Proposition~\ref{p=derivcurvbounds}.  Assumption VI. was proved in Theorem~\ref{t=firstderivest}. 
It remains to prove Assumptions III. and V. 
For $X\in\cT$, we continue to denote by $\rho_0(X)$ the WP distance from $X$ to $\partial \cT$. 

\medskip

\noindent{\bf Verifying assumption III.:  $\partial \left(\cT/\MCG[k]\right)$  is volumetrically cusplike.} 

\medskip

Given $\rho>0$, let $$E_\rho=\{X \in \cT/\MCG[k]: \rho_0(X)\leq \rho\}.$$
\begin{lemma}
\label{lem:volume}
We have $\Vol(E_\rho)=O(\rho^4)$
\end{lemma}

\begin{proof}
%Cover $E_\rho$ with a finite number of open sets in which we can use the plumbing coordinates $(s,t)$ introduced in %Section~\ref{ss=plumbing}. Fix one such neighborhood $U$ with  a collection $\{\alpha\}$ of short curves.
Fix a pants decomposition $\sigma$ that includes the short curves.  For each  curve $\alpha\in\sigma$, denote by $x_\alpha$ the function satisfying
$2\pi^2x_\alpha^2=\ell_\alpha$, where  
$\ell_\alpha$ is the length function. The theorem on p. 284  in 
\cite{WolpertBehavior} gives the asymptotic expansion 
$$g(\cdot,\cdot)\asymp \sum_\sigma 4dx_\alpha^2+x_\alpha^6d\theta_\alpha^2,$$ where $\theta_\alpha$ is the twist function. This gives that the volume element, which is the square root of the determinant of the metric $|g|^{1/2}$, is of the order $\prod_\alpha x_\alpha^3$. For the short curves,  $x_\alpha$ is comparable to  the distance to the boundary stratum in which $\alpha$ is pinched.  Thus we have  
%
%By the estimates in Section~\ref{s=comparison} the determinant of the matrix
%$(\la \phi_i,\phi_j \ra)$
%is $$\asymp \prod_{i=1}^p -|t_i|^2(\log|t_i|)^3.$$
%The volume element is then the determinant of the inverse and the volume is then found by  integration;  
$\Vol(E_\rho)=O(\rho^4)$.
%\prod_{i=1}^p\int_0^{|t_i|}\frac{|dt_i|}{-|t_i|^2(\log|t_i|)^3}=
%\prod_{i=1}^p\frac{1}{(-\log |t_i|)^2}.$$

\end{proof}

\noindent{\bf  Verifying assumption IV.:  $\cT/\MCG[k]$ has controlled injectivity radius.} 

\medskip

For $\alpha\in\cC$, denote by  $\tau_\alpha\in \MCG$  the Dehn twist about the curve $\alpha$. 
Given a simplex $\sigma=\{\alpha_1,\ldots, \alpha_p\}\in\cC(S)$,  let 
 $\Gamma(\sigma) = \la \tau_1,\ldots, \tau_p \ra$ be the abelian group generated by the Dehn twists about the $\alpha_i$.  Given $\epsilon>0$ let $\Omega(\sigma,\epsilon)=\{X:\forall \alpha\in\sigma, \ell_\alpha(X)<\epsilon\}$.  
\begin{lemma}
\label{lem:injective}
There exists $j_0\geq 1$ with the following property. For each $\epsilon>0$ there exists $c_0>0$ such that if $\phi\in \MCG[k]$ and $d_{WP}(X,\phi(X))<c_0$, then there exists $\sigma\in \cC(S)$ such that 
\begin{enumerate}
\item $X\in\Omega(\sigma,\epsilon)$, and 
\item for some $j\leq j_0$, $\phi^{j}\in \Gamma(\sigma)$.
\end{enumerate}
\end{lemma}

\begin{proof} Let $\epsilon>0$ be given.
Let $j_0$ be the product of $(3g-3+n)!$ and the product of the maximum orders of finite order elements on  surfaces of lower complexity.   The first conclusion (1) holds since  $\MCG[k]$ acts properly discontinuously without fixed points. 
  Now suppose the second statement (2) is not true; i.e., there exists $\epsilon$,  a sequence $X_m\in \Omega(\sigma,\epsilon)$,  and a sequence $\phi_m$ such that $d_{WP}(X_m,\phi(X_m))\to 0$ and yet for all $j\leq j_0$,  $\phi_m^{j}\notin \Gamma(\sigma)$.

Passing to a subsequence and applying an element $\psi_m\in\Gamma(\sigma)$  we can assume there is $\sigma$ such that $X_m$ converges to a noded surface $X_\sigma$.  For $\beta\in\sigma$ we have  $\ell_{\phi_m(\beta)}(\phi_m(X_m))=\ell_\beta(X_m)\to 0$. This implies that for $m$ sufficiently large, $\phi_m(\beta)\in\sigma$ as well.  Then for some $j\leq j_0$  the mapping class $\phi_m^{j}$ preserves  
the individual curves of $\sigma$.

The classification of elements of $\MCG$ implies that the restriction of  $\phi_m^{j}$  to each piece of  $X_\sigma$  is the composition of Dehn twists about boundary curves with an element that is either  pseudo-Anosov or finite order.  If  it is finite order in each piece then raising $\phi_m^j$ to a higher power we can assume $\phi_m^j$  is the product of Dehn twists, hence  in $\Gamma(\sigma)$, contrary to assumption. Thus  $\phi_m^{j}$ must be pseudo-Anosov on some piece.    
But then there is a uniform lower bound \cite[Theorem 7.6]{DaskWent}  for $d_{WP}(X_\sigma,\phi_m^{j}(X_\sigma))$ and thus a lower bound for 
$d_{WP}(X_m,\phi_m(X_m))$ for $m$ sufficiently large, a contradiction.  
\end{proof}

\begin{lemma}
There is a constant $c>0$ such that for any $X\in \cT/\MCG[k]$:
$$\inj(X)\geq c\rho_0(X)^3.$$
\end{lemma}
\begin{proof}
By Proposition 15 of \cite{WolpertGeometry} there is a positive constant $c>0$ such that for $X\in\Omega(\sigma,\epsilon)$,  $d_{WP}(X,\Gamma(\sigma)(X))\geq c\rho_0(X)^3$.  This bounds the injectivity radius from below. 
\end{proof}

Applying Theorem~\ref{t=generalergodicity}, we have now proved 
\begin{theorem}
\label{t=cover}
The Weil-Petersson flow on $T^1\cT/\MCG[k]$ is ergodic and has finite entropy.
\end{theorem}

\subsection{Ergodicity of the flow on $\cM^1(S)$: Proof of Theorem~\ref{t=main}}

The manifold $\cT/\MCG[k]$ is a finite branched cover over $\M$.
Let $h:X\to X$ be a conformal automorphism of finite order, and let
 $F(h)$ be the fixed point set of the induced action on $\cT$.
 It is known \cite{Ra} that if $S$ is compact and $h$ is not the hyperelliptic involution in genus $2$,  then $F(h)$
 has complex dimension at most $3g-5$. In fact $F(h)$ is the  Teichm\"uller space of the quotient orbifold  $X/h$. In genus $2$ the action induced by the hyperelliptic involution fixes every point of $\cT$. 
In the noncompact case where $S$ has punctures,  the complex dimension of $F(h)$ is at most $3g-4$. Let $F$ denote the union of the fixed point sets of the actions of all finite order elements of $\MCG(S)$, excluding the genus $2$ hyperellipic case. This is a countable union of lower dimensional Teichm\"uller spaces. 

\begin{lemma} 
$F$ is a closed subset of $\cT$, of  codimension at least $2$.
\end{lemma}
\begin{proof}
We have already seen that each fixed point set has real codimension at least $2$ so we need only check that 
the union is locally finite. Fix a compact set $K\subset\cT$. By the proper discontinuity of the action of 
$\MCG$ on $\cT$, there cannot be an infinite set of finite order elements each with a fixed point in $K$. Thus $K$ is intersected by only finitely many of the fixed point sets $F(h)$, and so the union of these sets is closed.
\end{proof}
 
We now finish the proof of ergodicity. 
Since the fixed point set of $\MCG$ has codimension at least $2$, the geodesic flow is defined almost everywhere on the quotient $\M^1$. If one has a positive measure invariant set in $E\subset \M^1$, then the lift of $E$ is a  positive measure invariant measure set in $T^1\cT/\MCG[k]$, which by ergodicity must have $0$ or full measure. The same is then true for $E$.  Hence the geodesic flow on $\cM^1$ is ergodic.  Moreover, any nontrivial factor of a Bernoulli flow is Bernoulli, and so the the geodesic flow on $\cM^1$ is Bernoulli as well.

Since the geodesic flow on $T^1\cT/\MCG[k]$ covers the geodesic flow on a full measure subset of $\cM^1$, it follows that the entropy of the flow on $\cM^1$ is also finite. This completes the proof of Theorem~\ref{t=main}. $\diamondsuit$

%\begin{remark} Among the class of Riemannian manifolds satisfying the general hypotheses of Theorem~\ref{t=generalergodicity}, it is possible to find examples where the bundles $E^u$ and $E^s$ fail to extend to continuous subbundles and to find examples where $E^u$ and $E^s$ are defined everywhere but fail to be $C^1$ (such examples were first constructed by Anosov in \cite{An67}).  None of the methods for constructing such systems apply to the Weil-Petersson metric, which leads to the natural question: in the case of  the WP geodesic flow,  do $E^u$ and $E^s$ extend to smooth subbundles? Equivalently, do the Busemann functions $b^u_{v}$ and $b^s_v$ depend smoothly on $v$?
%\end{remark}

\section{Appendix A: Bounding the second derivative of the geodesic flow}\label{s=secondderiv}

In this appendix we give precise estimates relating the norm of the first derivative of the geodesic flow, local bounds on the derivative of curvature, and the norm of the second derivative of the geodesic flow.  The results here will be used in Appendix B.

\subsection{More on the Sasaki metric and statement of the general result}

Let $M$ be a Riemannian manifold,
and let $\pi\colon TM\to M$ be the canonical projection. The Sasaki metric on $T^1M$ induces a Sasaki metric on  $T T^1M$,
which for brevity we will also call the Sasaki metric (although strictly speaking it is some sort of Sasaki Sasaki metric).  In general we will denote the Sasaki distance on $T^1M$ by 
$d_{Sas}$ and on $TT^1M$  by $\bd_{Sas}$.

Recall that for $v\in T^1 M$, each vector $\xi\in T_v T^1 M$ can be naturally
identified with a pair $(u,w)\in T^1_{\pi(v)} M\times  T^1_{\pi(v)} M$. 
The distance $\bd_{Sas}$  on $T T^1M$ induced by this metric can be estimated as follows.
Let  $\xi_0 = (u_0, w_0) \in  (T^1_{\pi(v_0)} M)^2$ and $\xi_1 = (u_1,w_1) \in (T^1_{\pi(v_1)} M)^2$ be tangent vectors in $T T^1 M$ based at $v_0$ and $v_1$ respectively.   Let $\sigma$ be a Sasaki geodesic
in $T^1 M$ from $v_0$ to $v_1$.  Let $P_\sigma \colon T^1_{\pi(v_0)} M \to  T^1_{\pi(v_1)} M$ be parallel translation along the curve of basepoints $\pi\circ \sigma$ in $M$.
The following lemma follows from the discussion in Section 2.
\begin{lemma}\label{l=sasasasa} For each $v_0$ there exists an $\epsilon>0$
such that if  $d_{Sas}(v_0,v_1)<\epsilon $, then
\begin{eqnarray*} 
\bd_{Sas}(\xi_0, \xi_1 )& \leq&  d_{Sas}(v_0,v_1)  + \|u_1 -  P_\sigma(u_0)\| +   \|w_1 -  P_\sigma(w_0)\|\\
&\leq& 2 \bd_{Sas}(\xi_0, \xi_1 ).
\end{eqnarray*}
\end{lemma}

The main result in this section is:
\begin{proposition}\label{p=generalsecondderiv} Let $M$ be an $m$-dimensional Riemannian manifold, possibly incomplete, and let $t_0\leq 1$ be a positive number. Let $\gamma\colon [-t_0, t_0]\to M$ be  a unit-speed  geodesic segment. 

Suppose that there exist constants $C_1, C_2, C_3 > 1 $ and $\epsilon_0>0$ such that 
for all $t\in (-t_0,t_0)$: 
\begin{enumerate}
\item  if $v\in T^1M$ satisfies $d_{Sas}(v,\dot\gamma(0)) <\epsilon_0$,
then
$$\max\{ \|D_{v} \varphi_{t}\| ,    \|D_{\varphi_{t}(v)} \varphi_{-t}\|\} \leq C_1;$$
\item  if $p\in M$ satisfies  $d( p, \gamma(t))<\epsilon_0$, then
$$\|R_p\|\leq C_2\quad\text{ and }\quad \|\nabla R_p\| \leq C_3.
$$
\end{enumerate}
Then there exists  $\epsilon_1>0$ such that for every $t\in (-t_0,t_0)$, for every pair $v_0,v_1 \in T^1 M$, with $d_{Sas}(v_i,\dot\gamma(0)) <\epsilon_1$,
and for all $\xi_i\in T^1_{v_i}T^1 M$, ($i=0,1$), we have:
$$
\bd_{Sas}( D\varphi_t(\xi_0),  D\varphi_t(\xi_1) ) \leq 
 (8m   C_1^4 C_2^2 C_3 )  \bd_{Sas}(\xi_0,\xi_1).
$$
\end{proposition}

\subsection{Variations of solutions to linear ODEs}
To prove Proposition~\ref{p=generalsecondderiv}, we first treat the linearized version of the problem. We begin with a basic fact about solutions to linear ODEs.
Consider a second-order linear ODE 
\begin{eqnarray}\label{e=jacobiODE}
x''(t) =  - \cR(t) x(t) 
\end{eqnarray}
where $\cR\colon [0,T] \to L(\RR^m)$ is continuous; in our application
$\cR(t)$ will be a matrix representing the
sectional curvature operator along a geodesic $\gamma$
and (\ref{e=jacobiODE}) will be the Jacobi equation in a suitably chosen 
coordinate system along $\gamma$.

Then  (\ref{e=jacobiODE})  can be 
transformed into a first order system in the standard way
by introducing the variable $z(t) = \left(\begin{array}{c}x(t)\\  y(t)\end{array}\right) \in \RR^{2m} $ and the additional
constraint $x'(t) = y(t)$.  Then $z$ satisfies the first order ODE:
 \begin{eqnarray}\label{e=1storder}
z'(t) = \left(\begin{array}{cc}
0 & I \\
-\cR(t) & 0
\end{array}\right)z(t).
\end{eqnarray}
The fundamental solution  $F(t)$  to this equation has the property 
that if $x(t)$ is a solution
to  (\ref{e=jacobiODE}) with initial values $x(0) = x_0$, $x'(0) = y_0$, then 
$ \left(\begin{array}{c}x(t)\\  x'(t)\end{array}\right) =  F(t)  \left(\begin{array}{c}x_0\\  y_0 \end{array}\right)$.

The following is a basic fact from the theory of ODEs.
\begin{proposition}\label{c=fundamentalbound}  Let $F_i\colon  [0,T] \to L(\RR^{2m})$ be the  fundamental solution to the differential equation
$x''(t) = -\cR_i(t) x(t)$, for $i=0,1$.  Then
$$
\|F_0(T) - F_1(T)\| \leq  T  \|F_0 \|_0  \|F_0^{-1}\|_0 \|\cR_0 - \cR_1\|_0 \|F_1\|_0.
$$
\end{proposition}

\subsection{Proof of Proposition~\ref{p=generalsecondderiv}}

We now return to the setting of differential geometry and finish the proof of Proposition~\ref{p=generalsecondderiv}.
Let $\gamma\colon [-t_0 ,t_0] \to M$ be given. We start with a lemma.
\begin{lemma}Under the assumptions of Proposition~\ref{p=generalsecondderiv} suppose 
   $p_0, p_1\in M$ satisfy  $d( p_i, \gamma(t))<\epsilon_0$, then for all
$v_i, w_i \in T_{p_i}^1M$, $i=0,1$, the curvature tensor  $R$  satisfies:
$$\|R(v_0, w_0 )w_0\| \leq C_2,
$$
and
$$d_{Sas}(R(v_0, w_0 )w_0 ,  R(v_1, w_1)w_1)  \leq C_3 (d_{Sas}(v_0, v_1) +
d_{Sas}(w_0,w_1)  ).
$$ 
\end{lemma}
\begin{proof} This follows in a straightforward way from the Mean Value Theorem and the 
hypotheses that $\|R_p\|\leq C_2$ and $\|\nabla R_p\| \leq C_3$, for all $p\in M$ with  $d( p, \gamma(t))<\epsilon_0$.
\end{proof}

Let $v_0,v_1 \in T^1M$ be unit tangent vectors in a neighborhood of $\dot\gamma(0)$, and let $\sigma\colon (-2,2)\to T^1M$ be a Sasaki geodesic with $\sigma(0)= v_0$ and $\sigma(1) = v_1$. Each $\sigma(s)$ determines a unit speed geodesic $\gamma_s\colon (-t_0,t_0)\to M$ with $\dot\gamma_s(0) = \sigma(s)$.
In this way $\sigma$ determines a variation of geodesics $\alpha \colon (-2,2)\times (-t_0,t_0) \to  M$  with the property that $\alpha(s,t) = \gamma_s(t)$.

We may assume that the norms of the derivatives of $\alpha$ are uniformly 
bounded from above by a constant, say $1$.  For $s\in (-2,2)$, let $L_s(t)  = \partial \alpha /\partial s (s,t)$ be the induced Jacobi field along $\gamma_s$.  Choose $\epsilon_1$
such that  if $d_{Sas}(v_i, \dot\gamma(0)) < \epsilon_1 $ for $i=0,1$, then
$d_{Sas}(\dot\gamma(t), \dot\gamma_s(t) ) < \epsilon_0$ for all $(s,t)\in (-2,2)\times (-t_0,t_0) $, where $\epsilon_0$ is given by the hypotheses of the proposition.
 %The Sasaki distance  between $\dot\gamma_0(t)$ and %\dot\gamma_s(t)$ is  bounded by the Sasaki length of a path joining them. The path $\alpha(s,t)$ joins the footpoints of these two vectors.
If   $d_{Sas}(v_i, \dot\gamma(0)) < \epsilon_1 $, then 
for any $(s,t)\in (-2,2)\times (-t_0, t_0)$ we  have
\begin{eqnarray*}
d_{Sas}(\dot \gamma_s(t) , \dot\gamma_0(t) )& \leq &\int_0^s \|(L_u(t), L_u'(t)))\|_{Sas}\, du.
%&\leq & 2 d_{Sas}(\dot \gamma_s(t) , \dot\gamma_0(t) ).
\end{eqnarray*}
Since $\sigma$ is a Sasaki geodesic the above inequality is an equality in the case of $t=0$; that is, $$\int_0^s\| (L_u(0),L_u'(0)) \|_{Sas}\, du=d_{Sas}(\dot\gamma_s(0),\dot\gamma_0(0)).$$ By the assumed bound (1) on the first derivative of the geodesic flow (which bounds the growth of Jacobi fields), we also have that
$\|(L_u(t), L_u'(t))\|_{Sas} \leq C_1 \|(L_u(0), L_u'(0))\|_{Sas}$,  for any $u,t$, and so 
$$\int_0^s \|(L_u(t), L_u'(t)))\|_{Sas}\, du\leq C_1\int_0^s\| (L_u(0),L_u'(0)) \|_{Sas}\, du.$$
Putting these inequalities together, we obtain:
\begin{eqnarray}\label{e=sasgrow}
d_{Sas}(\dot\gamma_s(t), \dot\gamma_0(t)) \leq  C_1  d_{Sas}(\dot\gamma_s(0), \dot\gamma_0(0)).
\end{eqnarray}

Our goal is to bound the Lipschitz norm of the derivative of the time-$t$ map of the geodesic flow  $\varphi_t$ at  $\dot\gamma_0(0)$.   The conclusion of
Proposition~\ref{p=generalsecondderiv} will follow if we show that
for any $(s,t)\in (-2,2)\times (-t_0,t_0)$, and any $\xi_0 \in T_{\dot \gamma_0(0)}^1 T^1 M$ and  $\xi_s \in T_{\dot \gamma_s(0)}^1 T^1 M$, we have:
\begin{equation}
\label{eq:goal}
\bd_{Sas} (D\varphi_t(\xi_0) ,   D \varphi_t(\xi_s) ) \leq  (4m   C_1^4 C_2^2 C_3 )  \bd_{Sas}( \xi_0, \xi_s).
\end{equation}

Recall that under the standard identification of $\xi_s  \in T_{\dot \gamma_s(0)} TM$ with a pair $(u_s, w_s) \in (T_{\gamma_s(0)} M)^2$, the vector $D_{\dot\gamma_s(0)} \varphi_t(\xi_s)$ is identified with $(J_s(t), J_s'(t))$, where $J_s$ is the solution to the 
 (second-order) Jacobi equation
\begin{eqnarray}\label{e=jacobigammas}
J'' + R(J,\dot\gamma_s)\dot\gamma_s = 0
\end{eqnarray}
with initial condition $(J_s(0), J_s'(0)) = (u_s,w_s) $.

To analyze the variation of solutions to this ODE, we  
fix convenient coordinates for the tangent bundle to the geodesic $\gamma_s$  in order to express  (\ref{e=jacobigammas}) as a matrix equation of the form (\ref{e=jacobiODE}). 
To this end, let $\{e_j(0,0) : j=1,\ldots, m\} $ be an orthonormal frame at $\gamma_0(0) = \alpha(0,0)$ spanning the tangent space $T_{\gamma_0(0)} M $. 
We first parallel translate this frame along
$\alpha(s,0)$ to obtain an orthonormal frame $\{ e_j(s,0)\}$ at $\gamma_s(0)$,
for $s\in (-2,2)$. We next parallel translate the frame  $\{ e_j(s,0)\}$ along
$\gamma_s(t)$, for $t\in (-t_0,t_0)$ to obtain a  frame $\{ e_j(s,t)\}$ at each point $\alpha(s,t)$.

\begin{lemma}\label{l=frametrans} For $j\in\{1,\ldots,m\}$, we have:
% $$d_{Sas} (e_j(0, t), e_j(s, t)) \leq 2 C_1  C_2 \,
%d_{Sas} (\dot\gamma_0(0), \dot\gamma_s(0)),$$ 
 $$d_{Sas} (e_j(0, t), e_j(s, t)) \leq  d(\gamma_0(0),\gamma_s(0))  + 2 C_1C_2  \, d_{Sas} (\dot\gamma_0(0), \dot\gamma_s(0)),$$ 
for all $(s,t)\in (-2,2)\times (-t_0,t_0)$.
\end{lemma}

\begin{proof}  Fix $j$.  Our construction of $e_j$ (using parallel translation)
gives that for all $s,t$:
\begin{eqnarray}\label{e=parallelej} 
\frac{D}{\partial s} e_j (s,0) = 0, \quad \hbox{ and }  \quad  \frac{D}{\partial t} e_j (s,t) = 0;
\end{eqnarray}
we would like to estimate $\frac{D}{\partial s} e_j (s,t)$ for general $s, t$.   To do this, we first estimate $\frac{D}{\partial t}  \frac{D}{\partial s} e_j (s,t)$.

It follows directly from the definition of the Riemannian curvature tensor and the joint integrability of the pair $\{L_s, \dot\gamma_s\}$  that
%\begin{eqnarray*}
$$ R(L_s(t), \dot\gamma_s(t) ) e_j (s,t) =
\frac{D}{\partial t} \frac{D}{\partial s} e_j (s,t) - \frac{D}{\partial s} \frac{D}{\partial t} e_j (s,t) 
 =   \frac{D}{\partial t} \frac{D}{\partial s} e_j (s,t),$$
%\end{eqnarray*}
where we have used the second part of (\ref{e=parallelej}) in the last step.
Applying the bound $\| R(L_s(t), \dot\gamma_s(t) ) e_j\| \leq C_2 \|L_s(t)\|$, we obtain that 
$ \left\|\frac{D}{\partial t} \frac{D}{\partial s} e_j (s,t)\right \| \leq	C_2 \|L_s(t)\|$.
Integrating this expression with respect to $t$, we then have the bound:
$$
\|\frac{D}{\partial s} e_j (s,t)\| \leq \|\frac{D}{\partial s}e_j(s,0)\|+ C_2 \int_0^t  \|L_s(u)\|\,du
= C_2 \int_0^t  \|L_s(u)\| \,du. $$
%$$\leq  C_2  \int_0^t C_1\|(L_s(0), L'_s(0))\|_{Sas}\,du 
%\leq  C_1 C_2 \|(L_s(0), L'_s(0))\|_{Sas},$$
%where we have used the first part of  (\ref{e=parallelej}).
Integrating again, this time with respect to $s$, and using Lemma~\ref{l=sasasasa} and (\ref{e=sasgrow}), we obtain:
\begin{eqnarray*}
d_{Sas}(e_j (0,t), e_j (s,t)) &\leq & d(\gamma_0(t),\gamma_s(t)) + \int_0^s \| \frac{D}{\partial s} e_j (u,t)\|\, du\\
&\leq& d(\gamma_0(0),\gamma_s(0)) +   \int_0^s\int_0^t \|L_w'(u)\|\,du\,dw  +  C_2 \int_0^s \int_0^t  \|L_w(u)\| \,du\,dw\\
&\leq& d(\gamma_0(0),\gamma_s(0)) +   2 C_2 \int_0^s\int_0^t   \|(L_w(u), L'_w(u))\|_{Sas}\,du\,dw \\
&\leq&  d(\gamma_0(0),\gamma_s(0)) + 2 C_1 C_2  \int_0^s  \|(L_w(0), L'_w(0))\|_{Sas}\,dw\\
&= & d(\gamma_0(0),\gamma_s(0))  + 2  C_1C_2  \, d_{Sas} (\dot\gamma_0(0), \dot\gamma_s(0)),
\end{eqnarray*}
%&\leq & (3 C_1C_2 )\, d_{Sas} (\dot\gamma_0(0), \dot\gamma_s(0)),
which is the desired bound.
\end{proof}

For $(s,t)\in (-2,2)\times (-t_0,t_0)$, the frame $\{e_j(s,t)\}$ gives an isometric linear isomorphism between
$\RR^m$ and $T_{\alpha(s,t)} M $:
$$
(x_1,\ldots, x_m) \mapsto  \sum_{j=1}^m x_j e_j(s,t).
$$
This in turn induces for each $(s,t)$ an isometric linear isomorphism $$I_{s,t} \colon \RR^{2m} \to T_{\dot\gamma_s(t)} T M \cong T_{\alpha(s,t)} M\times T_{\alpha(s,t)}M .$$
Lemma~\ref{l=frametrans} has the following immediate corollary. 
\begin{corollary}\label{c=saseuclid} For each $(s,t)\in (2,2)\times (-t_0,t_0)$ and each (Euclidean) unit vector $z \in \RR^{2m}$, we have
$$ d_{Sas}(I_{s,t}(z), I_{0,t}(z)) \leq  d(\gamma_0(0),\gamma_s(0))  + 2  C_1C_2  \, d_{Sas} (\dot\gamma_0(0), \dot\gamma_s(0)).
$$
\end{corollary}

Expressing the Jacobi equation 
 (\ref {e=jacobigammas}) along $\gamma_s$ in the coordinates on $T_{\gamma_s}M$ given by $I_{s,t}$,
we obtain the ODE:
\begin{eqnarray}\label{e=ODEincoords}
x''(t) = - \cR_s(t) x(t),
\end{eqnarray}
where
$
\left(\cR_s(t)\right)_{i,j} =  \langle \, R\left(e_i(s,t), \, \dot\gamma_s(t)\right) \dot\gamma_s(t), \, e_j(s,t)\,\rangle
$ .
%It is here one can see the convenience in choosing the frame $\{e_{j}\}$ to be parallel along the geodesic $\gamma_s$; otherwise, the Jacobi equation would include terms involving $x'$,
%which would complicate the linear analysis.
Denote by $F_s\colon (-t_0,t_0)\to L(\RR^{2m})$ the fundamental solution to  (\ref{e=ODEincoords}).  Proposition~\ref{c=fundamentalbound} implies that for any $(s,t) \in (-2,2)\times (-t_0,t_0)$,
we have
$$
\|F_0(t) - F_s(t)\| \leq  \|F_0 \|_0  \|F_s^{-1}\|_0 \|\cR_0 - \cR_s\|_0 \|F_s\|_0.
$$

Now the main hypotheses of the proposition, when combined with ~\ref{l=frametrans} and (\ref{e=sasgrow}), and the triangle inequality can be seen to give the upper bound:
\begin{eqnarray*}
\|\cR_0(t) - \cR_s(t)\| &\leq & (m C_1C_2^2 C_3 )  \, d_{Sas}(\dot\gamma_0(0), \dot\gamma_s(0))  ,
\end{eqnarray*}
(where we omit the  details)
and so
$$
\|F_0(t) - F_s(t)\| \leq ( m C_1C_2^2 C_3 )  \|F_0 \|_0    \|F_s\|_0  \|F_s^{-1}\|_0    d_{Sas} (\dot\gamma_0(0), \dot\gamma_s(0)),
$$
for all $(s,t)\in (-2,2)\times (-t_0,t_0)$.
The bounds on the first derivative of $\varphi_t$ imply that for all $(s,t)\in (-2,2)\times (-t_0,t_0)$, we have $\max\{\|F_s(t)\|, \|F_s^{-1}(t)\|\}  \leq   C_1$, which implies that
$$
\|F_0(t) - F_s(t)\| \leq  ( m C_1^4 C_2^2 C_3 )  \, d_{Sas}(\dot\gamma_0(0), \dot\gamma_s(0) ).
$$

Finally,  suppose that $\xi_0  = I_{0,0}(z_0)$ and
 $\xi_s  = I_{s,0}(z_s)$ are arbitrary unit tangent vectors
to $T^1M$ based at $\dot \gamma_0(0)$ and $\dot \gamma_s(0)$, respectively
(where $z_0,z_s$ are Euclidean unit vectors in $\RR^{2m}$).
Since $\frac{D}{\partial s} I_{s,0} = 0$, Lemma~\ref{l=sasasasa} implies that
the Sasaki distance $\bd_{Sas}(\xi_0, \xi_s )$ between $\xi_0$ and $\xi_s$
is uniformly comparable to $\|z_0 - z_1\| + d_{Sas}(\dot \gamma_0(0), \dot \gamma_s(0))$;
in particular:
\begin{eqnarray}\label{e=finalsascompare}
\|z_0 - z_1\| + d_{Sas}(\dot \gamma_0(0), \dot \gamma_s(0)) \leq 2 \bd_{Sas}(\xi_0, \xi_s ). 
\end{eqnarray}

We may then conclude using Corollary~\ref{c=saseuclid} and our previous estimates that:
\begin{eqnarray*}
&&\bd_{Sas}(D\varphi_t(\xi_0), D\varphi_t(\xi_s) )\, = \,\bd_{Sas}(I_{0,t} \left(F_0(t) z_0\right) ,
I_{s,t}\left( F_s(t) z_s \right) )\\
&\leq &   \bd_{Sas}(I_{0,t} \left(F_0(t) z_0\right) ,
I_{0,t}\left( F_s(t) z_s \right) ) +  \bd_{Sas}(I_{0,t} \left(F_s(t) z_s\right) ,
I_{s,t}\left( F_s(t) z_s \right) )   \\
&\leq&
 \|F_0(t) z_0 - F_s(t) z_s \|  + \|F_s(t) z_s\| \,\left( d(\gamma_0(0),\gamma_s(0))  + 2 C_1C_2  \, d_{Sas} (\dot\gamma_0(0), \dot\gamma_s(0)) \right)\\
&\leq&    \|(F_0(t) - F_s(t) ) z_0 \|  
+   \|F_s(t) (z_0 - z_s)\| + C_1 \left( d(\gamma_0(0),\gamma_s(0))  + 2 C_1C_2  \, d_{Sas} (\dot\gamma_0(0), \dot\gamma_s(0)) \right)\\
&\leq& (m C_1^4 C_2^2 C_3 )\, d_{Sas}(\dot\gamma_0(0), \dot\gamma_s(0) )
+     C_1 \|z_0 - z_s \|  + C_1 \left( d(\gamma_0(0),\gamma_s(0))  + 2  C_1C_2  \, d_{Sas} (\dot\gamma_0(0), \dot\gamma_s(0)) \right)\\
&\leq& ( m C_1^4 C_2^2 C_3 )\, d_{Sas}(\dot\gamma_0(0), \dot\gamma_s(0) )
+     C_1 \|z_0 - z_s \|   + 3 C_1^2 C_2  \, d_{Sas} (\dot\gamma_0(0), \dot\gamma_s(0))\\
 &\leq & ( 4m C_1^4 C_2^2 C_3 )\, \left( \|z_0 - z_s \| +  d_{Sas}(\dot\gamma_0(0), \dot\gamma_s(0) )\right )  \\
&\leq &  (8 m   C_1^4 C_2^2 C_3 )\, \bd_{Sas}(\xi_0, \xi_s ),
\end{eqnarray*}
where we used (\ref{e=finalsascompare}) in the last step.
This proves the desired inequality (\ref{eq:goal}) and  completes the proof of Proposition~\ref{p=generalsecondderiv}.

\section{Appendix B: Proof of Proposition~\ref{p=aclam}:  verifying the conditions of  Katok-Strelcyn-Ledrappier}

In this appendix we prove  Proposition~\ref{p=aclam}.
We assume the conditions I.-VI. in Theorem~\ref{t=generalergodicity}.
Let $V\subset T^1 N$ be the  set of $v\in T^1N$ such that $\varphi_t(v)\in T^1 N$, for
all $t\in (-1,1)$.  Fix $t_0\in (0,1)$ and consider the restriction of the time-$t_0$ map
$\varphi_{t_0}$ to $V$.   To prove Proposition~\ref{p=aclam}, we verify that the main  hypotheses in \cite{KS} hold for the map  $\varphi_{t_0}\colon V \to T^1N$.  The main results in \cite{KS} then imply the conclusions of Proposition~\ref{p=aclam}.  
To paraphrase \cite{KS}, the  conditions we will verify ensure that the set of singularities of the map $\varphi_{t_0}$  is ``thin" and that the first and second derivatives of $\varphi_{t_0}$ grow moderately near this set. 

In the setup of \cite{KS}, the background hypotheses are: $X$ is a compact metric space, and $V$ is an open and dense subset of $X$ carrying a Riemann structure with controlled singularities near $X\setminus V$.  In our application, $V$ is the set described above, endowed with the Sasaki Riemann structure, and $X = \bar{T^1 N}$ is the  completion of $T^1 N$ in the Sasaki 
distance metric $d_{Sas}$.  We first verify that $X$ is compact, which establishes condition (A) of \cite{KS}.

\begin{lemma} $(\bar{T^1N}, d_{Sas})$ is compact.
\end{lemma}

\begin{proof} Let $\la v_{n,m} \ra_m$ be a sequence of elements of $\bar{T^1N}$, 
where for each
$m\geq 1$, $\la v_{n,m}\ra$ is a Sasaki Cauchy sequence in $T^1N$.  Since
$d_{Sas}(v,w) \geq d(\pi(v),\pi(w))$, it follows that for each $m$, the sequence $\la \pi(v_{n,m} ) \ra $  is  Cauchy in $N$; since $\bar N$ is compact, by passing to a subsequence in the $m$'s, we may assume that
$\la \pi(v_{n,m}) \ra_m$  converges to a 
Cauchy sequence $\la x_n \ra$ in $N$.  What this means is that for every $\epsilon>0$ there exists an $m_0>0$ such that for $m\geq m_0$, we have
$$
\lim_{n\to \infty} d(\pi(v_{n,m}), x_n) <\epsilon.
$$

Now for each $n$,
consider the collection $\{\hat v_{n,m} \,\mid\, m\geq 1\} \subset T_{x_n}^1 N$ obtained by parallel 
translating each $v_{n,m}$ along a geodesic from $T_{\pi(v_{n,m})}N$ to $T_{x_n}N$.
Using compactness of $T_{x_n}^1 N$ and a diagonal argument, we obtain a subsequence $m_k$ such that for each $n$, $\hat v_{n,m_k}$  converges as $k\to \infty$ to an element $\hat v_n\in T_{x_n}^1N$, uniformly in $n$; that is, for every $\epsilon>0$, there exists $k_0 > 0$ such that 
for all $k> k_0$ we have 
$$
\lim_{n\to \infty} \|\hat v_{n,m_k} - \hat v_n\| <\epsilon.
$$

Since the Sasaki distance $d_{Sas}(v_{n,m_k},\hat v_n)$ is bounded by
$d(\pi(v_{n,m_k}), x_n) + \|\hat v_{n,m_k} - \hat v_n\|$, it follows that for every $\epsilon>0$ there exists a $k_1>0$ such that for all $k\geq k_1$,
\begin{eqnarray*}
\lim_{n\to \infty}  d_{Sas}( v_{n,m_{k}}, \hat v_n) 
&\leq & \lim_{n\to \infty} d(\pi(v_{n,m}), x_n)+  \|\hat v_{n,m_k} - \hat v_n\| < 2\epsilon.
\end{eqnarray*}
Hence  $\la v_{n,m_{k}}\ra_{m_k}$
converges as $k\to \infty$ to the Sasaki Cauchy sequence $\la \hat v_n \ra\in \bar{T^1N}$. 
\end{proof}

Clearly $V$ is an open and dense subset of  $\bar{T^1 N}$.  Let $S = \bar{T^1N}\setminus V$.
The Sasaki distance from $v$ to the singular set  $S$ is bounded above by the distance from $\pi(v)$
to $\partial N$.

\subsubsection{More (yet) on the Sasaki metric}

Condition (B) in \cite{KS}, which concerns the Riemann structure on $V$,  has three parts that require verification.  In this subsection, we establish bounds on the derivatives of the Sasaki exponential map $\bexp\colon TV\to V$, which we will then use to verify these conditions as well as later conditions on $\varphi_{t_0}$.  To control the Sasaki exponential map, we will need to control the first three derivatives of the Sasaki metric; these can be related to the higher order derivatives of the metric on $N$ via the following lemma.

\begin{lemma}\label{l=sascurvbound} There exists a cubic polynomial $C\colon \RR^3\to \RR$ such that for any Riemannian manifold $N$ and any  $v\in T^1_xN$, the Sasaki curvature tensor $R_{Sas}$ satisfies
$$
\| (R_{Sas})_v \| + \|\nabla (R_{Sas})_v \|  \leq C(\|R_x\|, \|\nabla R_x\|, \|\nabla^2 R_x\|),
$$
where $R$ is the Riemannian curvature tensor on $N$.
\end{lemma}

\begin{proof}
The sectional curvatures of the Sasaki metric on the unit tangent bundle can be computed
as follows  \cite{KowSez}.  We use the usual identification
 $T_{(x,u)}TN \cong T_x N\times T_xN$. Let $\Pi$ be a plane
in $T_{(x,u)}T^1N$, and choose an orthonormal basis $\{ (v_1,w_1), (v_2,w_2)\}$ for $\Pi$ satisfying
$\|v_i\|^2 + \|w_i\|^2 = 1$ for $i=1,2$ and $\langle v_1, v_2 \rangle = \langle w_1, w_2\rangle = 0$. Then the Sasaki sectional curvature of $\Pi$ is given by
\begin{eqnarray*}
K_{Sas}(\Pi) &=& \langle R_x(v_1,v_2)v_2,v_1\rangle + 3\langle R_x(v_1,v_2)w_2, w_1\rangle+
\|w_1\|^2 \|w_2\|^2 - \frac34 \|R_x(v_1,v_2)u\|^2 \\
&&+ \frac14\|R_x(u, w_2)v_1\|^2 + \frac14\|R_x(u, w_1)v_2\|^2+
\frac12\langle R_x(u, w_1)w_2,R_x(u, w_2)v_1\rangle\\
&& - \langle R_x(u, w_1)v_1,R_x(u, w_2)v_2\rangle+ \langle(\nabla_{v_1}R)_x(u, w_2)v_2,v_1\rangle +\langle(\nabla_{v_2}R)_x(u, w_1)v_1,v_2\rangle.
\end{eqnarray*}
The conclusion now follows from the Chain Rule and  well-known identities relating the sectional curvatures with the norm of the curvature tensor.
\end{proof}

The next lemma will be used to bound the derivative of the Sasaki exponential map.

\begin{lemma}
\label{lem:Jac2} Let $Y$ be a Riemannian manifold, and let $J$ be a Jacobi field along a geodesic $\gamma\colon [-\delta_0,\delta_0]\to Y$ satisfying $J(0) = 0$ and
$\|J'(0)\| = 1$.  Suppose that
$$\sup_{|t| < \delta_0} \|R_{\gamma(t)}\| \leq   R_0$$\
for some $R_0>1$.
Let $\epsilon\in(0,1)$ be given, and
let  $t_0 = \min\{\delta_0, \epsilon /(3 R_0)\}$.
Then for all $|t|\leq t_0$  we have
$$
(1-\epsilon)|t| \leq \|J(t)\| \leq (1+\epsilon) |t|
\quad\hbox{ and }\,
\|J'(t)\| \leq 1+\epsilon.
$$

\end{lemma}

\begin{proof}
Let $a(t) = \|J(t)\|$, and let $b(t) = \|J'(t)\|$.
Then the Cauchy-Schwarz inequality implies
$$| (a^2) ' |  = |2aa'|= |2\la J,J'\ra|  \leq  2 ab,$$    
and since $|t|< \delta_0$:
$$| (b^2)' | =  |2 bb'|=| 2 \la J',J''\ra | =  | 2 \la J',  R(\dot\gamma,
J)\dot\gamma  \ra  |  \leq   2 R_0 a b.$$
We conclude that  wherever $|a|$ and $|b|$ are not zero, we have
$|a'| \leq b$  and  $ |b'| \leq R_0 a$.

We are assuming that $a(0) = 0$ and $b(0) = 1$. Without loss of generality, assume that $|a(t)| > 0$ for $t > 0$ (otherwise, we may replace $t=0$ with a positive value of $t$ in the following argument). 
From this we obtain the integral inequality, for $t\geq 0$:
\begin{eqnarray}\label{e=abintupbound} |a'(t) |       \leq         1  +  \int_0^t |b'(s) |\, ds           \leq
 1  +  R_0 \int_{0}^t   a(s)\, ds.
\end{eqnarray}
Suppose that, for some $t_1\in (0,t_0)$  we have  
 $|a'(t)| < 1+\epsilon$ for all $t \in [0, t_1)$ and  $|a'(t_1)| = 1+\epsilon$. 
Then $a(t)  <  (1+\epsilon) t$, for all  $t \in [0, t_1)$; combined
with (\ref{e=abintupbound}), this gives that
$$|a'(t_1)| \leq 1 + R_0 \int_0^{t_1} (1+\epsilon) s \,  ds   <   1 +  \frac{R_0(1+\epsilon)}{2}t_1^2 <  1+\epsilon,$$
since $\epsilon\in(0,1)$ implies that
$$t_1^2 < t_0^2 \leq  \frac{\epsilon^2}{9 R_0^2} <  \frac{2\epsilon}{R_0(1+\epsilon)}.$$
This contradicts our assumption that
$|a'(t_1)| = 1+\epsilon$.
We conclude that $|a'(t)| < 1+\epsilon$ for all  $t \in (0,t_0)$; similarly,
$|a'(t)| < 1+\epsilon$, for all  $t \in (-t_0,0)$.
From this we conclude that $a(t)\leq (1+\epsilon) |t|$ for all 
$|t|\leq t_0$.

We now prove the lower  bound. Since $b(0)= 1$  and $|b'(t)|\leq R_0 a(t)$, for $|t|\leq t_0$
we have
$$b(t) \geq 1 -  \frac{(1+\epsilon) R_0 t^2}{2}.$$
On the other hand, we know that
$$(a^2)'' = 2b^2 - 2\la R(J,\dot\gamma)\dot\gamma, J \ra   \geq   2b^2 -
2 R_0 a^2
$$
$$> 2 \left[  \left( 1- \frac{(1+\epsilon) R_0 t^2}{2}\right)^2 - (1+\epsilon)^2 R_0 t^2  \right]
> 2[1-  2  (1+\epsilon)^2 R_0 t^2 ],
$$
(using the lower bound for $b(t)$ and upper bound of $(1+\epsilon)|t|$ for  $a(t)$).
Now, since
$$t^2\leq  t_0^2 \leq \frac{\epsilon^2}{9 R_0^2} < \frac{\epsilon^2}{2(1+\epsilon)^2 R_0},$$
 we find that
$$(a^2)''(t) > 2[1-  2  (1+\epsilon)^2 R_0 t^2 ] > 2(1-\epsilon^2).$$  
But then
$2 a(t) a'(t) = (a^2)'(t) >  2(1-\epsilon^2) |t| $, and again 
 using the upper bound on $a$, we get
$$a'(t) > \frac{(1-\epsilon^2) |t|}{a(t)} >  \frac{(1-\epsilon^2)| t|}{(1+\epsilon)| t|} = 1-\epsilon;$$
hence $a(t)> (1-\epsilon) |t|$.

Finally, since $b(0)=1$ and $|b'| \leq R_0 a \leq R_0(1+\epsilon)$, it follows that
$|b(t)|\leq 1 + |t| R_0(1+\epsilon)$, and so for $|t| < |t_0|$, we have
$|b(t)| \leq 1 + \epsilon(1+\epsilon)/3 < 1+\epsilon$. The final conclusion follows.
\end{proof}

We apply this lemma to the Sasaki exponential map $\bexp\colon TV\to V$ to obtain:

\begin{proposition}\label{p=dexpbound} There exist constants $\delta_1>0$  and $k_1 >1$ such that for every 
$v_0\in V$, if $d_{Sas}(v_0, S) < \delta_1$, then for all $v\in V$ with $d_{Sas}(v,v_0) < d_{Sas}(v_0,S)^{k_1}$: 
$$ 1 -  d_{Sas}(v_0, S) \leq   \|D_v\bexp_{v_0}^{-1}\|^{-1} \leq     \|D_\xi\bexp_{v_0} \| \leq 1 + d_{Sas}(v_0, S),$$
where $\xi=\bexp_{v_0}^{-1}(v)$.
\end{proposition}

\begin{proof}  Let $v_0\in V$ and $\xi=\bexp_{v_0}^{-1}(v)$.
 Let $\hat \xi=\frac{\xi}{\|\xi\|}$ be the unit vector in the direction of $\xi$. 
Suppose $\xi'\in T^1_{v_0} V$ is an orthogonal unit vector.  
Let $a(s,t)=(\hat \xi+s\xi')t$ be the $1$-parameter 
family of rays through the origin in $T_{v_0} V$. 
 Let $$\alpha(s,t)=\bexp_{v_0}\circ \, a(s,t)$$ be the  
$1$-parameter family of image geodesics in $V$. 
We consider the corresponding Jacobi field $J(t)$ along $\alpha(0,t)$ defined by 
$J(t)=\partial \alpha(s,t)/\partial s$ at $s=0$. 
Clearly $J(0)=0$ and $J'(0)=\xi'$.
Setting $t_1=\|\xi\|$, by the chain rule we have 
 $$\|J(t_1)\|=\|t_1 D_\xi\bexp_{v_0}(\xi')\|.$$
Thus we have to bound $\|\frac{J(t_1)}{t_1}\|$ above and below. 

By Lemma~\ref{l=sascurvbound}, the sectional curvatures of the Sasaki metric on the unit tangent bundle are bounded polynomially in terms of the absolute value of the curvature and the derivative of the curvature of the original metric. Assumption IV. gives a bound for these latter quantities, and therefore a polynomial bound on the curvatures in the Sasaki metric,  in the reciprocal of the distance to the singular set $S$. It follows that there exist $k_0>1$ and $\delta_0>0$ such that for all $v\in V$ with  $d_{Sas}(v,S)<\delta_0$,
the Sasaki curvature tensor $R_{Sas}$ satisfies
$$
\| (R_{Sas})_v\| < d_{Sas}(v,S)^{-k_0}.
$$

Let $k_1 =  k_0 + 2$.  Then there exists $\delta_1\in (0,1/3)$ such that if
$d_{Sas}(v_0,S) < \delta_1$ and
$$d_{Sas}(v_0,v)\leq d_{Sas}(v_0,S)^{k_1},$$ 
then the maximum norm $R_0$ of the Sasaki curvature tensor 
along the geodesic joining $v_0$ to $v$ also satisfies $R_0< d_{Sas}(v_0,S)^{-k_0}$. 
Lemma~\ref{lem:Jac2} implies that 
\begin{eqnarray}\label{e=epsbounds}
1 -\epsilon 
\leq \left \|\frac{J(t_1)}{t_1} \right \| < 1 + \epsilon,
\end{eqnarray}
provided 
$\epsilon > 3 R_0 |t_1| = 3 R_0 d_{Sas}(v_0,v)$.  Hence if $d_{Sas}(v_0,S) < \delta_1 $ and $d_{Sas}(v,v_0) \leq d_{Sas}(v_0,S)^{k_1}$,
then (\ref{e=epsbounds}) holds for $\epsilon = d_{Sas}(v_0,S)$,
since 
$$ 3 R_0 d_{Sas}(v_0,v) < 
3 d_{Sas}(v_0,S)^{-k_0} \cdot d_{Sas}(v_0,S)^{k_1} = 3 d_{Sas}(v_0,S)^{2}<
d_{Sas}(v_0,S)=\epsilon .$$

\end{proof}

The next proposition gives bounds on the second derivative of $\bexp$, which we will later use to verify condition (1.3) of \cite{KS}.

\begin{proposition}
\label{p=secondexp}
There exist constants $\delta_2>0$  and $ k_2 >1$ such that for every 
$v_0\in V$, if $d_{Sas}(v_0, S) < \delta_2$, then for all 
%$v_1, v_2, w_1, w_2\in V$ with $\max\{d_{Sas}(v_i,v_0), d_{Sas}(w_i, w_0} < d_{Sas}(v_0,S)^{k_2}$ for $i=1,2$, we have
$\xi_i, \eta_i\in T_{v_0} V$ with $(\xi_1,\eta_1) \neq (\xi_2,\eta_2)$ and
 $\max\{\|\xi_i\|, \|\eta_i\|\} < d_{Sas}(v_0,S)^{k_2}$ for $i=1,2$, we have:
$$ d_{Sas}(v_0, S)^{k_2}  \leq
\frac{\bd_{Sas} \left(D_{\xi_1} \bexp_{v_0}(\eta_1),  D_{\xi_2} \bexp_{v_0}(\eta_2) \right) }
{\|\xi_1 - \xi_2\| + \|\eta_1 - \eta_2\| }
\leq
d_{Sas}(v_0, S)^{-k_2}.
$$
%$$
%\max\{\|D^2_\xi \bexp_{v_0}\|,  \|D^2_v \bexp_{v_0}^{-1}\|\} < d_{Sas}(v_0, S)^{-k_2}$$
%where $\xi_i=\bexp_{v_0}^{-1}(v_i)$ and  $\eta_i=\bexp_{v_0}^{-1}(w_i)$,  for $i=1,2$.
\end{proposition}

\begin{proof} 

Suppose that $v_0\in V$ is fixed and $v_1$ lies in a neighborhood of $v_0$. Let $\xi_1 = \bexp_{v_0}^{-1}(v_1)$.  For $\xi_2 \in T_{v_0}V$, the map $D_{\xi_2}\bexp_{v_0}$ is a linear transformation between $T_{v_0}V$ and $T_{v_2}V$, where $v_2 = \bexp_{v_0}(\xi_2)$.  The Sasaki connection defines a trivialization of the bundle $TV$ in a neighborhood of the fiber over $v_1$;
in these coordinates, a vector $\eta_2\in T_{v_2} V$ is sent  to the pair $(v_2, P_{v_2,v_1}(\eta_2))\in V\times T_{v_1}V$, where
$P_{v_2,v_1}\colon T_{v_2}V\to T_{v_1}V$ is parallel translation along the unique local geodesic from $v_2$ to $v_1$.  The Sasaki metric $\bd_{Sas}$ on $TV$ is comparable in this trivializing neighborhood to the product metric on $V\times T_{v_1}V$.  In these coordinates, there is a well-defined second derivative $D^2_{\xi_1} \bexp_{v_0}\colon T_{v_0}V\times T_{v_0}V  \to T_{v_1} V$ obtained by differentiating the second component of $D_{\xi}\bexp_{v_0}$ with respect to $\xi$ and evaluating at $\xi_1$.  By the Mean Value Theorem, to prove the conclusions of the proposition, it suffices to bound
$\|D^2_{\xi} \bexp_{v_0}(\eta,\eta)\|$ from above and below, for all $\xi$ in a neghborhood of the origin in $T_{v_0}V$ and $\eta$ a unit vector perpendicular to $\xi$.

To this end, fix $v_0\in V$ and $v\in V$ in a neighborhood of $v_0$, and let $\xi = \bexp_{v_0}^{-1}(v)$.
 Let $\hat \xi=\frac{\xi}{\|\xi\|}$ be the unit vector in the direction of $\xi$, and
suppose $\eta\in T^1_{v_0} V$ is an orthogonal unit vector.  
As in the proof of Proposition~\ref{p=dexpbound}, we consider the variation of geodesics
 $$\alpha(s,t)=\bexp_{v_0}\circ \, a(s,t),$$ 
where $a(s,t)=(\hat \xi+s\eta)t$.
Define $Z, J$ and $Q$ by
$$Z(s, t) = \frac{D}{\partial t} \alpha(s,t); \, J(s,t) = \frac{D}{\partial s}\alpha(s,t);
\, Q(s,t) = \frac{D^2}{\partial s^2} \alpha(s,t) = \frac{D}{\partial s} J(s,t).
$$
The chain rule implies that
$$
Q(0,t) = D^2_{a(0,t)} \bexp_{v_0}(t \eta, t\eta),$$
and so
$$\|  D^2_{\xi} \bexp_{v_0} (\eta,\eta) \| = \frac{1}{\|\xi\|^2} \|Q(0,\|\xi\|)\|,$$
since $\xi = a(0,\|\xi\|)$.

Observe that for $s$ fixed, $J(s,\cdot)$ is a Jacobi field down the geodesic $\alpha(s,\cdot)$
and so satisfies the Jacobi equation
$$
\frac{D^2}{\partial t^2} J  =  R_{Sas}(Z,J)Z.
$$
From this, the definition of $Q$ and symmetries of the curvature tensor it follows that
$$
\frac{D^2}{\partial t^2} Q = \frac{D^2}{\partial t^2}\frac{D}{\partial s} J  =  R_{Sas}(Z,J)J' 
+  \frac{D}{\partial t} \left(R_{Sas}(Z,J)J \right) + \frac{D}{\partial s}\frac{D^2}{\partial t^2} J 
$$
$$= R_{Sas}(Z,J)J'  +  \frac{D}{\partial t}\left( R_{Sas}(Z,J)J \right)  +  \frac{D}{\partial s} R_{Sas}(Z,J)Z 
$$
$$ =  R_{Sas}(Z,J)J' +  \left(\left(\frac{D}{\partial t} R_{Sas} \right)(Z,J)J +  R_{Sas}(Z, J')J+ R_{Sas}(Z,J)J'\right),
$$
$$
 +\left(\left(\frac{D}{\partial s} R_{Sas} \right)(Z,J)Z +  R_{Sas}(J', J) Z + R_{Sas}(Z, Q)Z+ R_{Sas}(Z,J)J'\right)
$$
$$ = \left(\frac{D}{\partial t} R_{Sas} \right)(Z,J)J  + \left(\frac{D}{\partial s} R_{Sas} \right)(Z,J)Z 
 +  4  R_{Sas}(Z,J)J'  + R_{Sas}(Z, Q) Z,
$$
where $'$ denotes the derivative with respect to $t$, and we have also used the facts that $Z'=0$
and $(D/\partial s ) Z = J'$ .
Then $\|Q''(0,t)\| \leq C_1(t) + \|Q(0,t)\| C_2(t)$, where
$$C_1(t) =  \|\left(\nabla R_{Sas}\right)_{\exp_{v_0}(t\hat\xi)}\|\,(\|J(0,t)\| + \|J(0,t)\|^2) + 4  \|\left(R_{Sas}\right)_{\exp_{v_0}(t\hat\xi)}\| \|J(0,t)\|\,\|J'(0,t)\|,$$
and
$$ C_2(t) = \|\left(R_{Sas}\right)_{\exp_{v_0}(t\hat\xi)}\|.$$

Assumption IV. and  Lemma~\ref{l=sascurvbound} imply that there exists $k_0>1$ such that
$$
\max\{\|\left(R_{Sas}\right)_{v_0}\| , \|\left(\nabla R_{Sas}\right)_{v_0}\| \} < d_{Sas}(v_0,S)^{-k_0}.
$$
Fix $\delta_2\in (0,1/22)$ such that if $d_{Sas}(v_0,S) < \delta_2$, then
$$
\sup_{|t|\leq  d_{Sas}(v_0,S)^{k_0+1}} \max\{\|\left(R_{Sas}\right)_{\exp_{v_0}(t\xi)}\| , \|\left(\nabla R_{Sas}\right)_{\exp_{v_0}(t\xi)}\| \} < d_{Sas}(v_0,S)^{-k_0 - 1}.
$$
Assume that $d_{Sas}(v_0,S) < \delta_2$.
Lemma~\ref{lem:Jac2} implies that for $|t| < d_{Sas}(v_0,S)^{k_0+2}$, both
$\|J(0,t)\|$ and $\|J'(0,t)\|$ are bounded by $2$, and so 
$$C_1(t) \leq 22\, d_{Sas}(v_0,S)^{- k_0-1} <  d_{Sas}(v_0,S)^{- k_0-2},$$ and
$$C_2(t) \leq  d_{Sas}(v_0,S)^{-k_0-1} < d_{Sas}(v_0,S)^{- k_0-2}.$$

Let $q(t) = \|Q(0,t)\|$ and let $r(t) = \| Q'(0,t)\|$.
As in the proof of Lemma~\ref{lem:Jac2}, we have that
$$
| q q'| = | \la Q, Q'\ra | \leq  q r,\quad\text{ and }\,\, | r r' | \leq  |\la Q',Q''\ra | \leq  r( C_1  + q C_2) .
$$
Note that $q(0)=r(0) = 0$. 
An analysis similar to that in the proof of Lemma~\ref{lem:Jac2} (whose details we omit) shows that
for $|t| <  d_{Sas}(v_0,S)^{k_0+2}$, we have
$$q(t) \leq t^2 d_{Sas}(v_0, S)^{-k_0-2}.$$
Hence, if $\|\xi\| \leq d_{Sas}(v_0, S)^{k_0+2}$, then
$$
\|  D^2_{\xi} \bexp_{v_0} (\eta,\eta) \|  = \frac{q(\|\xi\| )}{\|\xi\| ^2}  \leq  d_{Sas}(v_0, S)^{-k_0-2}.
$$
Hence an  upper bound for the ratio in the conclusion of the proposition holds, with the exponent $k_2 = k_0+2$.  A lower bound for this ratio on the order of $d_{Sas}(v_0, S)^{-k_2}$
follows from the upper bound 
on  $\|  D^2_{v} \bexp_{v_0}\|$ we have just obtained,  the upper bounds on 
$\|D_v\bexp_{v_0}^{-1}\|$
and  $\|D_\xi \bexp_{v_0}\|$
given by Proposition~\ref{p=dexpbound}, and the fact that  for an invertible matrix-valued function
$\xi\mapsto A(\xi)$, one has 
$$
D_\xi (A^{-1}) = -A^{-1}(\xi) ( D_\xi A ) A^{-1}(\xi).
$$
The details are left to the reader.
%$$
%D^2\bexp^{-1} = - \left( D\bexp^{-1}\times D\bexp^{-1}\right) \circ \left( D^2\bexp\right) \circ \left( D\bexp^{-1} \times \, D\bexp^{-1}\right).
%$$
\end{proof}

\subsubsection{Verifying condition (B) in \cite{KS}}

For $v\in V$,
let $\operatorname{\bf inj}(v)$ denote the radius of injectivity of the Sasaki exponential map 
$\bexp_{v}\colon T_{v}V \to V$.
Since $d_{Sas}(v,w)\geq d(\pi(v),\pi(w))$, the controlled injectivity assumption V. implies
that 
$$
\operatorname{\bf inj}(v) \geq \inj(\pi(v)) \geq C d(\pi(v), \partial N)^\beta \geq  C d_{Sas}(v, S)^\beta. 
$$
This implies condition (Ba) of \cite{KS}. 
Conditions (Bb) and (Bc) in \cite{KS} follow in a straightforward way from  Proposition~\ref{p=dexpbound}.

\subsubsection{Verifying conditions (1.1) -- (1.4) of \cite{KS}}

Conditions (1.1) -- (1.4) of \cite{KS} concern the volume of the singular set $S$ and the behavior of $\varphi_{t_0}$ near $S$.  
Condition (1.1) of \cite{KS}, which concerns the volume of a neighborhood of $S$,  follows directly from Lemma~\ref{l=specialnbdvol}. Condition (1.2) of \cite{KS}, concerning the integrability of $\log_+\|D\varphi_{t_0}\|$,  follows immediately from Lemma~\ref{l=integrability}.
Condition (1.4) of \cite{KS} requires a bound on $\|D_{v_0}\varphi_{t_0}\|$
on the order of $d_{Sas}(v_0,S)^{-\beta}$, for some $\beta>0$.   This follows in a straightforward way from assumption VI. This leaves condition (1.3).

Fix $v_0\in V$, and let $\Phi = \Phi_{v_0} \colon T_{v_0}V\to T_{\varphi_{t_0}(v_0)}V$ be defined by
$$
\Phi = \bexp_{\varphi_{t_0}(v_0)}^{-1} \circ \, \varphi_{t_0} \circ \bexp_{v_0}.
$$
Condition (1.3) of \cite{KS} requires a bound on the second derivative of $\Phi$
as an inverse power of the distance to the singular set, which follows from the 
next proposition.

\begin{proposition} 
There exist constants $\delta_3>0$  and $ k_3 >1$ such that for every 
$v_0\in V$, if $d_{Sas}(v_0, S) < \delta_3$, then for all $v\in V$ with $d_{Sas}(v,v_0) < d_{Sas}(v_0,S)^{k_3}$:
$$
\|D^2_\xi \Phi_{v_0}\|  < d_{Sas}(v_0, S)^{-k_3},$$
where $\xi=\bexp_{v_0}^{-1}(v)$.
\end{proposition}

\begin{proof}
Choose constants $k_2>1$ and $\delta_2>0$ satisfying the conclusions
of Proposition~\ref{p=secondexp} and
such that if $d_{Sas}(v_0,S) < \delta_2$, then for every $|t|\leq t_0$
$$
\sup_{d_{Sas} (v,\varphi_t(v)) \leq  d_{Sas}(\varphi_t(v),S)^{k_2}} \max\{\|\left(R_{Sas}\right)_{v}\| , \|\left(\nabla R_{Sas}\right)_{v}\| \} < d_{Sas}(v_0,S)^{-k_2}.
$$
By assumption VI.,  there exist $\delta_3 < \min \{\delta_2, 1/(8n)\} $ and $k_2'>k_2$ such that for $d_{Sas}(v_0,S) < \delta_3$,
and all $|t|\leq t_0$:
$$
\max\{ \|D_{v_0}\varphi_{t}\|, \|D_{v_0}\varphi_{-t}\| \} \leq d_{Sas}(v_0,S)^{-k_2'}.
$$
Proposition~\ref{p=generalsecondderiv} implies that if $\xi, \xi' \in TV$ satisfy
$d_{Sas}(\pi_V(\xi),\pi_V(\xi')) <  d_{Sas}(v_t,S)^{k_2}$ then
$$
\bd_{Sas}( D\varphi_t(\xi),  D\varphi_t(\xi') ) \leq 
( 8n)\, d_{Sas}(v_0,S)^{-4 k_2' - 3 k_2}   \bd_{Sas}(\xi_0,\xi_1)$$
$$ \leq   d_{Sas}(v_0,S)^{-7 k_2'}.
 \bd_{Sas}(\xi_0,\xi_1).$$
To bound the norm $\|D^2_\xi \Phi_{v_0}\|$ it suffices to bound the Lipschitz constant of the
map $D_\xi\Phi_{v_0}$ in a small neighborhood of $\xi$. 
This in turn is bounded by the product of the Lipschitz constants of the
three factors  $D\bexp_{\varphi_{t_0}(v_0)}^{-1}$,  $D\varphi_{t_0}$ and  $D\bexp_{v_0}$ in the composition defining $D_\xi\Phi_{v_0}$.  

The Lipschitz constants for  $D\bexp_{v_0}$  and $D\bexp_{\varphi_{t_0}(v_0)}^{-1}$ are both bounded by
% for nearby points by the norms $\|D^2\bexp_{v_0}\|$ and $\|D^2\bexp_{\varphi_{t_0}(v_0)}^{-1}\|$, which by 
Proposition~\ref{p=secondexp} on the order of $d_{Sas}(v_0,S)^{-k_2}$. We have just shown that the Lipschitz constant for $D\varphi_{t_0}$ is bounded on the order of $d_{Sas}(v_0,S)^{-7k_2'}$. Hence the Lipschitz constant of  $D_\xi\Phi_{v_0}$ is bounded on the order of $d_{Sas}(v_0,S)^{-k_4}$, for $k_4 =2k_2 + 7k_2'$.  
%The details of this argument are left to the reader.
\end{proof}

 This completes the verification of the hypotheses in \cite{KS} implying the conclusions of Proposition~\ref{p=aclam}.

\subsection{Additional conditions in \cite{KS} implying finite, positive entropy}

The final conclusion of Theorem~\ref{t=generalergodicity} that remains to be proved 
concerns the entropy of $\varphi$.  The positivity of the entropy follows from \cite{KS} and
the hypotheses we have just verified.  Finitude of the entropy requires that an additional hypothesis -- Condition (C) -- hold.  As stated in \cite{KS}, condition (C) is the requirement that the capacity of the space $X=\bar{T^1N}$ be finite.  In fact, a slightly weaker condition is required, which is given by the following proposition. Recall that $U_\rho$, for $\rho>0$, denotes the set of $v\in T^1N$ such that $d(\pi(v), \partial N) <\rho$.
\begin{proposition} There exists  $q >1$ such that  if $\rho_0>0$ is sufficiently
small, then for
any $\rho<\rho_0$
there is a cover of $T^1N\setminus U_{\rho_0}$ by open balls of radius $\rho$, whose cardinality does not exceed $\rho^{-q}$.
\end{proposition}

\begin{proof} Proposition~\ref{p=dexpbound} implies that there exist $\delta>0$ and $k>1$ such that for $\rho_0<\delta$ and all
$v\in T^1N\setminus U_{\rho_0}$, the derivative of
the Sasaki exponential map $D_{\xi}\bexp_v$  and its inverse have norm bounded by $2$,
for all $\|\xi\| < \rho_0^{k}$.  Hence on a ball of radius $\rho_0^{k}$ in $T^1N\setminus U_{\rho_0}$, the Sasaki metric is uniformly comparable to Euclidean; in particular, the volume of a ball of radius $\rho\leq \rho_0^k$ is bounded 
%above by   $c\rho^n$  and 
below by  $c^{-1} \rho^n$, where $c>1$ is a universal constant. 

The Vitali Covering Lemma states that if $\mathcal B$ is any collection of 
balls in a metric space, then there exists a subcollection ${\mathcal B}'\subset {\mathcal B}$ such that the elements of  ${\mathcal B}'$ are pairwise disjoint, and
$$
\bigcup_{B\in {\mathcal B}'}  5B \supset  \bigcup_{B\in {\mathcal B}} B, 
$$
where $5B$ denotes the ball concentric with $B$ of $5$ times the radius.

Let $\mathcal B$ be a finite cover of the set  $T^1N\setminus U_{\rho_0}$ by metric balls of radius $\rho^k$, and let ${\mathcal B}'$ be a subcollection of disjoint balls 
supplied by the Vitali lemma. Then  the collection $\{5B : B\in {\mathcal B}'\}$ is a covering of $T^1N\setminus U_{\rho_0}$ by balls of radius $5\rho^k$.  If $\rho_0$ was chosen sufficiently small, then  $5\rho^k<\rho$, and the balls in this cover can be expanded to give a cover by balls of radius $\rho$. The cardinality of this cover equals the cardinality of ${\mathcal B}'$; this number can be bounded above using the volume:
$$
\left(\# {\mathcal B}' \right) \times  (c^{-1} \rho^{n k}) \leq \sum_{B\in {\mathcal B}'} m(B) 
= m\left(    \bigcup_{B\in {\mathcal B}'} B  \right)   \leq m(T^1N) = 1.
$$
Thus $
\# {\mathcal B}' \leq c\rho^{-nk}
$,
for all $\rho<\rho_0$.   This implies the conclusion of the proposition, with $q= nk+1$.

\end{proof}

\end{document}